\documentclass[10pt]{amsart}
\usepackage{amssymb}
\newcommand{\qsub}{{\qq^{\protect\vphantom{\mathsf{T}}}}}
\newcommand{\Bsub}{{B^{\protect\vphantom{\mathsf{T}}}}}
\newcommand{\Csub}{{C^{\protect\vphantom{\mathsf{T}}}}}
\newcommand{\qq}{{\mathbf{q}}}
\newcommand{\qtran}{{\qq^{\mathsf{T}}}}
\newcommand{\negone}{{\mathbf{-1}}}
\newcommand{\minusone}{{\mathbf{-\!1}^{\protect\vphantom{\mathsf{T}}}}}
\newcommand{\rd}[1]{\smash{ \underline{ {#1}{\vphantom{j}} } }}
\newcommand{\uy}{{\rd y}}
\newcommand{\refl}{t}
\newcommand{\mindouble}{}
\newcommand{\GL}{\mathrm{GL}}
\newcommand{\sgn}{\mathrm{sgn}}
\newcommand{\diag}{\mathrm{diag}}
\newcommand{\dcat}{\mathrm{D}}
\newcommand{\qdcat}{\mathrm{Q}}
\newcommand{\CC}{\mathbb C}

\newcommand{\tensor}{\mathop{\otimes}}
\newcommand{\utensor}{\mathop{\underline\otimes}}

\newcommand{\dash}{\nobreakdash-\hspace{0pt}}
\newcommand{\Ob}{\mathrm{Ob}}
\newcommand{\YD}[1]{{{\vphantom{X}}_{#1}^{#1}\mathcal{YD}}}
\newcommand{\lgen}{\text{$<$}}
\newcommand{\rgen}{\text{$>$}}

\newcommand{\field}{\CC}
\newcommand{\Symm}{\mathbb{S}}
\newcommand{\lcprod}{\rtimes}
\newcommand{\rcprod}{\ltimes}

\DeclareMathOperator{\id}{id}
\DeclareMathOperator{\End}{End}
\DeclareMathOperator{\im}{im}

\theoremstyle{plain}
\newtheorem{theorem}{Theorem}[section]
\newtheorem{proposition}[theorem]{Proposition}
\newtheorem{problem}[theorem]{Problem}
\newtheorem{lemma}[theorem]{Lemma}
\newtheorem{corollary}[theorem]{Corollary}

\theoremstyle{definition}
\newtheorem{remark}[theorem]{Remark}
\newtheorem*{exercise}{Exercise}
\newtheorem{definition}[theorem]{Definition}
\newtheorem{example}[theorem]{Example}

\advance\oddsidemargin by-0.5in
\advance\evensidemargin by-0.5in
\advance\textwidth by 1in

\long\def\startcmm#1\endcmm{}

\begin{document}

\title{Noncommutative Dunkl operators and braided Cherednik algebras}

\author{Yuri Bazlov}
\address{Mathematics Institute, University of Warwick, Coventry CV4 7AL, UK}
\email{y.bazlov@warwick.ac.uk}

\author{Arkady Berenstein}
\address{Department of Mathematics, University of Oregon, 
Eugene, OR 97403, USA} 
\email{arkadiy@math.uoregon.edu}
\thanks{The authors acknowledge support of the EPSRC grant
  EP/D065801/1 (Y.B.) and  NSF grant DMS-0501103~(A.B.)}
\subjclass[2000]{Primary
20G42; % Quantum groups (quantized function algebras) and their representations
%16S99,  % Rings and algebras arising under various constructions:
       % None of the above, but in this section
Secondary
%05E15, % Combinatorial problems concerning the classical groups
% 12E05% Polynomials (irreducibility, etc.)
% 14E05% Rational and birational maps
% 14M15, % Grassmannians, Schubert varieties, flag manifolds
% 14M17, % Homogeneous spaces and generalizations
% 14M99, % Special varieties: None of the above, but in this section
16S80, % Deformations of rings  
%17B37,% Quantum groups (quantized enveloping algebras) and related
% deformations
% 17B20% Simple, semisimple, reductive (super)algebras (roots)
% 17B67% Kac-Moody (super)algebras (structure and representation theory)
% 17B99. % Lie algebras and Lie superalgebras: None of the above, but in this section
% 22E46% Semisimple Lie groups and their representations
% 22E65% Infinite-dimensional Lie groups and their Lie algebras
% 20G05% Representation theory [of linear algebraic groups]
% 16S34 % Group rings, Laurent polynomial rings
% 13A99 % General commutative ring theory: None of the above, but in
% this section
 20F55% Reflection and Coxeter groups
}

\pagestyle{myheadings}
\markboth{Y.~BAZLOV and A.~BERENSTEIN}{NONCOMMUTATIVE DUNKL OPERATORS}
\begin{abstract}
We introduce \emph{braided Dunkl operators} $\rd \nabla_1,\ldots,\rd
\nabla_n$ that are acting on a $\qq$-symmetric algebra $S_\qsub(\CC^n)$
and $\qq$-commute. Generalising the approach of Etingof and Ginzburg,
we explain the $\qq$-commutation phenomenon by constructing
\emph{braided Cherednik algebras} $\rd{\mathcal H}$ for which  the above
operators form a representation.  We classify  all  braided Cherednik 
algebras using the theory of braided doubles developed in our previous paper. 
Besides ordinary rational Cherednik algebras, our
classification gives new algebras 
$\rd{\mathcal H}(W_+)$ attached to an infinite family of
subgroups of even elements in complex reflection groups, so  that the
corresponding 
braided Dunkl operators $\rd \nabla_i$ pairwise anti-commute.  
We explicitly compute these new operators in terms of braided partial
derivatives and $W_+$-divided differences. 
\end{abstract}
\maketitle
\setcounter{tocdepth}{1}
\tableofcontents

\setcounter{tocdepth}{1}
\setcounter{section}{-1}
%%%

\section{Introduction}

%Let $P=\CC[x_1,\dots,x_n]$ be the algebra 
%of polynomial functions on the affine space $\Aff^n$. 

In 1989, Dunkl introduced  the remarkable operators 
$\nabla_1,\dots,\nabla_n$ that act on the polynomial algebra
$\CC[x_1,\dots,x_n]$  by
%\colon P \to P$, $i=1,\ldots,n$, given by
$$
\nabla_i  = \frac{\partial}{\partial x_i} + c \sum_{j\ne i} 
\frac{1}{x_i-x_j}(1-s_{ij})\ ,
$$ 
where $s_{ij}$ is the automorphism of $\CC[x_1,\dots,x_n]$ switching
the variables 
$x_i$ and $x_j$; together, the $s_{ij}$ generate the symmetric group
$\Symm_n$.  
Dunkl operators are a very special 
deformation (with parameter $c\in \CC$) 
of the partial derivatives $\frac{\partial}{\partial x_i}$,
because they commute: $\nabla_i \nabla_j = \nabla_j \nabla_i$ for all
$i$, $j$. This striking fact, 
originally proved in \cite{Dun} by a  lengthy direct computation, 
is interpreted in quantum mechanical terms as the integrability of the
Calogero-Moser system, and  algebraically via  
rational Cherednik algebras introduced in 2002 by Etingof and
Ginzburg \cite{EG}. A family of commuting Dunkl operators is attached 
not only to $\Symm_n$ but to any finite group generated by (complex)
reflections.

The motivating idea behind the present paper is to introduce analogues
of Dunkl operators on 
noncommutative deformations of the  symmetric algebra $S(V)$, 
e.g., on the $\qq$-symmetric algebra 
$$
S_\qsub(V)= \CC\langle x_1, \dots, x_n \mid x_i x_j = q_{ij} x_j x_i\
\text{for}\ i<j\rangle
$$
where $V$ is a $\CC$-vector space with basis $x_1,\ldots,x_n$ 
and $\qq=(q_{ij})$ is a complex $n\times n$ matrix.

Assume that $q_{ij}q_{ji}=q_{ii}=1$ for all $i,j$ and 
define braided partial derivatives $\rd\partial_i\colon 
S_\qsub(V)\to S_\qsub(V)$ by
$$
\rd \partial_i(x_1^{a_1}\cdots x_i^{a_i}\cdots
x_n^{a_n})=a_i\, q_{1,i}^{a_1}\cdots q_{i-1,i}^{a_{i-1}} \, x_1^{a_1}
\cdots x_i^{a_i-1}\cdots x_n^{a_n} \ .
$$ 
Clearly, $\rd \partial_i\rd \partial_j=q_{ij}\rd \partial_j\rd
\partial_i$ and $\rd \partial_i x_j=q_{ji}\rd \partial_i
x_j+\delta_{ij}$ for all $i,j$.

It turns out that if $q_{ij}=-1$ for all $i\ne j$ 
(denote such matrix $\qq$ by $\negone$), 
then  we can introduce the following deformations $\rd
\nabla_i\colon S_\minusone(V)\to S_\minusone(V)$ of the braided
partial derivatives $\rd \partial_i$. Let
${\mathcal C'}\subset {\mathcal C}$ be finite (cyclic) subgroups of
$\CC^\times$ such that $-1\in {\mathcal C}$ (i.e., $\mathcal C$ is
of even order), and let $c\colon {\mathcal C'}\to \CC$ be a function
$\varepsilon'\mapsto c_{\varepsilon'}$. For each $i=1,\ldots,n$
define the operator $\rd \nabla_i$ on the skew-field of fractions of
$S_\minusone(V)$ as follows:
\begin{equation}
\label{eq:braided Dunkl anticommutative}
\displaystyle{\rd \nabla_i=\rd \partial_i+c_1\sum_{j\ne i,\,\varepsilon \in {\mathcal C}}
\frac{x_i+\varepsilon x_j}{x_i^2-\varepsilon^2
  x_j^2}(1-\sigma_{ij}^{(\varepsilon)}) 
+\sum_{\varepsilon'\in {\mathcal
    C}'\setminus\{1\}}\frac{c_{\varepsilon'}}{1-\varepsilon'} \cdot
\frac{1}{x_i}(1-t_i^{(\varepsilon')})} \,,
\end{equation}
%for $i=1,\ldots,n$, 
where
$t_i^{(\varepsilon)}$, $\sigma_{ij}^{(\varepsilon)}$  are algebra
automorphisms of $S_\minusone(V)$ defined by  
$$
t_i^{(\varepsilon)}(x_k)=\begin{cases}
x_k & \text{if $k\ne i$}\\
\varepsilon x_i& \text{if $k=i$}\\
\end{cases}\ ;\qquad
%for $i$ and $\varepsilon\in \CC^\times$; ~ $
\sigma_{ij}^{(\varepsilon)}(x_k)=
\begin{cases}
x_k & \text{if $k\notin \{i,j\}$}\\
\varepsilon x_j& \text{if $k=i$}\\
-\varepsilon^{-1} x_i& \text{if $k=j$}\\  
\end{cases} 
\qquad
\text{for all $i\ne j$, $\varepsilon\in \CC^\times$.}
$$
We refer to these new operators $\rd \nabla_i$ as
\emph{noncommutative} or \emph{braided} Dunkl operators due to the following 
\begin{theorem} 
\label{th:braided dunkl} 
The operators $\rd \nabla_1,\ldots,\rd \nabla_n$ preserve 
$S_\minusone(V)$ and satisfy
$\rd \nabla_i\rd \nabla_j+\rd \nabla_j\rd \nabla_i=0$ for all $i\ne j$. 

 \end{theorem}

Therefore, the operators $\rd \nabla_i$ may be viewed as an
anti-commutative analogue of a Calogero-Moser integrable system.  
We will elaborate on 
the connection with integrable systems in a separate paper.

We  prove Theorem \ref{th:braided dunkl}
% in Section \ref{subsect:Braided Dunkl operators} 
by means of {\it braided Cherednik algebras} which we introduce
(following the logic of  \cite{EG}) as the main tool for establishing
anti- and $q$-commutation relations between 
operators~\eqref{eq:braided Dunkl anticommutative} and their generalisations. 

Namely, let $W_{{\mathcal C}, {\mathcal C}'}$ be the subgroup of
$\GL(V)$ generated by all $\sigma_{ij}^{(\varepsilon)}$,
$\varepsilon\in {\mathcal C}$ and $t_i^{(\varepsilon')}$,
$\varepsilon'\in {\mathcal C}'$. Define $\rd  {\mathcal
  H}_c(W_{{\mathcal C}, {\mathcal C}'})$ to be the subalgebra of
$\End_\CC(S_\minusone(V))$ generated by $W_{{\mathcal C}, {\mathcal C}'}$,
$y_i=\rd\nabla_i$, and operators of multiplication by $x_i$.

%The following is our first main result.

\begin{theorem} 
\label{th:PBW anti-commutating}

$(a)$
In the algebra $\rd  {\mathcal H}_c(W_{{\mathcal
    C}, {\mathcal C}'})$, the generators satisfy:  
\begin{itemize}
\item[$(i)$] $x_ix_j+x_jx_i=y_iy_j+y_j y_i=0$ for all $i\ne j$;
\item[$(ii)$] $wx_iw^{-1}=w(x_i), w y_iw^{-1}=w(y_i)$ for all $w\in W_{{\mathcal C}, {\mathcal C}'}$, $i=1,\ldots,n$;
\item[$(iii)$] 
$y_jx_i+x_iy_j=c_1\sum\limits_{\varepsilon\in {\mathcal C}} \varepsilon\sigma_{ij}^{(\varepsilon)}$ for all $i\ne j$, and 
\item[]    
                $y_ix_i-x_iy_i=1+ c_1\sum\limits_{j\ne i,\ \varepsilon\in {\mathcal C}} \sigma_{ij}^{(\varepsilon)}+\sum\limits_{\varepsilon'\in {\mathcal C}'\setminus\{1\}} c_{\varepsilon'}t_i^{(\varepsilon')}$ for $i=1,\ldots,n$.
\end{itemize}              

$(b)$ As an associative algebra, $\rd  {\mathcal H}_c(W_{{\mathcal C}, {\mathcal C}'})$ is determined by relations $(i)$--$(iii)$ and admits a triangular decomposition 
$$\rd  {\mathcal H}_c(W_{{\mathcal C}, {\mathcal
    C}'})=S_\minusone(V)\tensor\CC W_{{\mathcal C}, {\mathcal
    C}'}\tensor S_{\minusone}(V^*)\,.$$
\end{theorem}

We prove Theorem \ref{th:PBW anti-commutating} in Section
\ref{subsect:sigma}. In what follows we will refer to $\rd{\mathcal
  H}_c(W_{{\mathcal C}, {\mathcal C}'})$ as a {\it negative braided
  Cherednik algebra}. Using the Verma module
$M_{\mathit{triv}}=S_\minusone(V)$  for 
$\rd{\mathcal H}_c(W_{{\mathcal C}, {\mathcal C}'})$, induced from the
trivial representation of $W_{{\mathcal C}, {\mathcal C}'}$, we obtain
the following refinement of Theorem \ref{th:braided dunkl} (to be
proved in Section \ref{subsect:Braided Dunkl operators}). 

\begin{theorem} 
\label{th:braided Verma} 
The generators $y_i$ of the negative braided Cherednik algebra $\rd
{\mathcal H}_c(W_{{\mathcal C}, {\mathcal C}'})$ act on its Verma
module $M_{\mathit{triv}}=S_\minusone(V)$ by braided Dunkl operators
\eqref{eq:braided Dunkl anticommutative}.  
\end{theorem}

\begin{remark} 
\label{rem:degenerate}
In fact, if one drops the constant term $1$ in the
  second relation of Theorem \ref{th:PBW anti-commutating}$(a)(iii)$
  and at the same time drops the braided derivative $\rd
  \partial_i$ in \eqref{eq:braided Dunkl anticommutative}, then one
  obtains a \emph{degenerate} braided Cherednik algebra for which
  Theorems \ref{th:PBW anti-commutating} and \ref{th:braided Verma}
  are also valid. The latter algebra admits a finite-dimensional
  quotient by the $W_{{\mathcal C}, {\mathcal C}'}$-invariant ideals
  of $S_\qsub(V)$ and $S_\qsub(V^*)$, which is an analogue of restricted
  Cherednik algebras; see \cite{Go}. 
\end{remark}

\begin{remark} 

In a series of papers \cite{KW1,KW2} Ta Khongsap and
  Weiqiang Wang have discovered a different class of algebras with
  triangular decomposition and anti\dash commuting generators.  This similarity inspired us to start a new project \cite{BB3}
  where constructions of \cite{KW1,KW2} are uniformly treated in the context of
  braided doubles (developed in \cite{BB} and in Section~\ref{sect:qd}
  of the present paper).

\end{remark}

%One can classify 
The groups $W_{{\mathcal C}, {\mathcal C}'}$ are
classified in terms of  the infinite family of \emph{classical} 
complex reflection groups 
$G(m,p,n) = \Symm_n\rcprod T(m,p,n)$, $m\ge 1$, $p\,|\,m$, where 
$\Symm_n\subset \GL(V)$ is the group of permutation matrices and 
$T(m,p,n) = \{ \diag(\varepsilon_1,\dots,\varepsilon_n) \mid
\varepsilon_i^m=1 \ \forall i, \ (\varepsilon_1\dots
\varepsilon_n)^{m/p}=1 \}$.
%For any finite subgroup $G\subset \GL(V)$, let $d$ be the least
%positive integer such that $\det(g)^d\in\{\pm 1\}$ for all $g\in G$.
%Denote
%$$
%    \sgn=\det\nolimits^d\colon G \to \{\pm 1\}, 
%\qquad
%    G_+ = \ker\ \sgn,
%$$
%and call $g\in G_+$ even elements of $G$. 
It is not difficult to show that     
\begin{itemize}
\item[---] if $|\mathcal C'|$ is even, 
$W_{{\mathcal C}, {\mathcal C}'} = G(m,p,n)$ with $m=|\mathcal C|$,
  $p=|\mathcal C / \mathcal  C'|$;
\item[---] if $|\mathcal C'|$ is odd, 
$W_{{\mathcal C}, {\mathcal C}'} = G(m,p,n)_+$, a subgroup of index
  $2$ in $G(m,p,n)$ with $m=|\mathcal C|$,
  $p=|\mathcal  C/\pm \mathcal C'|$, consisting of $g$ such that 
$\det(g)^{|\mathcal C'|}=1$.
\end{itemize}
%(In particular, one gets the usual notion of even elements in
%classical Coxeter groups.)
%
(Note however, that the generators $\sigma_{ij}^{(\varepsilon)}$ of
$W_{{\mathcal C}, {\mathcal C}'}$ are of order $4$
and are not complex reflections.)
The smallest example of such group in rank $n$ is  
$W_{\{\pm 1\},\{1\}}=B_n^+$ of even elements in a Coxeter
group of type $B_n$; see Example~\ref{ex:bn+}.

%
%\begin{proposition} Let ${\mathcal C}'\subset  {\mathcal C}$ be arbitrary finite cyclic subgroups of $\CC^\times$. Denote $m:=|{\mathcal C}|$ and $p:=|{\mathcal C}/{\mathcal C}'|$, so that  $|{\mathcal C}'|=\frac{m}{p}$. Then
%$$W_{{\mathcal C}, {\mathcal C}'}=\{w\in G(m,p,n):\det(w)\in {\mathcal C}'\} \ .$$
%In particular, if $\frac{m}{p}$ is even, then $W_{{\mathcal C}, {\mathcal C}'}=G(m,p,n)$ and if $\frac{m}{p}$ is odd, then $W_{{\mathcal C}, {\mathcal C}'}\subsetneq G(m,p,n)$.
%
%
%\end{proposition}     

Having been inspired by the construction of the above negative braided
Cherednik algebras corresponding to the matrix $\qq=\negone$, 
we formulated and solved the following problem.

\begin{problem}   
\label{prob:braided Cherednik algebras} Let $\qq=(q_{ij})$ be a
complex $n\times n$ matrix with $q_{ij}q_{ji}=q_{ii}=1$
and $V$ be a vector space with basis $x_1,\ldots,x_n$. Find all finite
groups $W\subset \GL(V)$ acting on $S_\qsub(V)$ by automorphisms and
all algebras $\rd  {\mathcal H}$ generated by $W$, $V$, and $V^*$ such
that:  

(a) $\rd  {\mathcal H}$ admits a triangular decomposition $\rd
{\mathcal H}=S_\qsub(V)\tensor\CC W\tensor S_\qsub(V^*)$, where
$S_\qsub(V)\tensor\CC W$ and $\CC W\tensor S_\qsub(V^*)$ are semidirect product algebras;

(b) $y_j x_i-q_{ij} x_i y_j\in \CC W$  for all $i,j$, 
where   $y_1,\ldots,y_n$ is the  basis of $V^*$ dual to $x_1,\ldots,x_n$.
\end{problem}

Clearly,  ordinary rational Cherednik algebras $H_c(W)$ associated
to complex reflection groups $W\subset \GL(V)$ and the above
negative algebras 
$\rd  {\mathcal H}_c(W_{{\mathcal C}, {\mathcal C}'})$ 
solve Problem \ref{prob:braided Cherednik algebras} for special
examples of the matrix $\qq$.
We refer to solutions of Problem \ref{prob:braided Cherednik
  algebras}  as \emph{braided Cherednik algebras}.

Let us briefly describe how to construct all braided Cherednik
algebras out of the above building blocks (see Section
\ref{sect:Classification} for details). Informally speaking, we prove
that each braided Cherednik algebra is a braided product (which, once
again, justifies the name) of rational Cherednik algebras and the
negative algebras  $\rd  {\mathcal H}_c(W_{{\mathcal C}, {\mathcal C}'})$.

More precisely,  let each of $\mathcal H(W_1),\dots,\mathcal H(W_m)$
be either a rational Cherednik algebra of an irreducible complex
reflection group $W_k$ (one of the groups $G(m,p,n)$ and $G_4,\dots,G_{37}$
in Shephard-Todd's classification \cite{ST}) or a negative braided Cherednik
algebra of $W_k=W_{{\mathcal C}, {\mathcal C}'}$. 
One has $\mathcal H(W_k)\cong S(V_k)\tensor \CC W_k \tensor
S(V_k^*)$ or $S_\minusone(V^k) \tensor \CC W_k \tensor
S_\minusone(V_k^*)$, respectively. Choose $m(m-1)/2$
arbitrary parameters $r_{kl}\in 
\CC^\times$, $1\le k<l\le m$. 
Define $\mathcal H$ to be the algebra generated by 
all $\mathcal H(W_k)$ subject to the relations
$$
   xx'=r_{kl} x'x, \quad yy' = r_{kl} y'y, \quad xy' = r_{kl}^{-1}
   y'x, \qquad yx' = r_{kl}^{-1}  x'y
$$
for $x\in V_k$, $y\in V_k^*$, $x'\in V_l$, $y'\in V^*_l$, and 
the relation that $W_k$ commute with $\mathcal H(W_l)$ for $k\ne l$.  
In Section~5, we prove 
\begin{theorem} 
\label{thm:brprod}
(a) $\mathcal H$ factorises as a tensor product of its subalgebras, 
\begin{equation}
\label{eq:braided factorization1}
\mathcal H = \mathcal H(W_1)\tensor \dots \tensor \mathcal H(W_m),
\end{equation}
and is a braided Cherednik algebra of the group $W_1\times \dots
\times W_m$;

(b) each braided Cherednik algebra of a group $W$ is isomorphic to one
of such algebras $\mathcal H$ (under a simple assumption of minimality
of $W$, see Section \ref{subsect:Braided Cherednik algebras: the main
  structural theorem}). 
\end{theorem}

The braided Cherednik algebra $\mathcal H$
so constructed has factorisation $S_\qsub(V)\tensor \CC W \tensor
S_\qsub(V^*)$ with $V=V_1\oplus \dots \oplus V_m$. Here $\qq$ is the
matrix consisting of $m^2$ blocks $M_{kl}$ (of size $\dim V_k \times
\dim V_l$), $k,l=1,\dots,m$. The 
block $M_{kk}$  has all off\dash diagonal entries equal to $1$
(respectively 
$-1$) if $\mathcal H(W_k)$
is a rational (respectively negative braided) Cherednik algebra. 
The block $M_{kl}$ has all entries equal to $r_{kl}$ if $k<l$ or to
$r_{lk}^{-1}$ if $k>l$. 

The proof of the theorem  is based on the theory of \emph{braided
  doubles} developed in our previous paper
\cite{BB}. Results specific to \emph{quadratic} braided
  doubles over group algebras are given in Section \ref{sect:qd}. 
Using these results, in Section~\ref{sect:q_cher} we introduce and
classify 
$\qq$-Cherednik algebras as specific quadratic doubles with
triangular factorisation 
$S_\qsub(V)\tensor \CC \widetilde W\tensor S_{\qtran}(V^*)$ 
where $\widetilde W$ is a certain Abelian extension of $W$ depending
on $\qq$ 
(and $\qtran$ stands for the transposed matrix). 
Based on this classification and the \emph{braided reduction} introduced
in Section \ref{subsect:Braided reduction}, we prove that 
each braided Cherednik algebra is naturally a subalgebra of one
of the $\qq$\dash Cherednik algebras (Proposition \ref{prop:brred}), and, on
the other hand, that each braided Cherednik algebra naturally admits a
factorisation \eqref{eq:braided factorization1} in an appropriate
braided tensor category (Theorem \ref{thm:main}). 

Let us illustrate our construction of braided Cherednik algebras and
the corresponding braided Dunkl operators for any finite Abelian group
$W$, i.e., 
$W={\mathcal C}_1\times {\mathcal C}_2\times \cdots \times
{\mathcal C}_n$, where each ${\mathcal C}_i$ is a finite (cyclic)
subgroup of $\CC^\times$, and an arbitrary $n\times n$ matrix $\qq$ as
above.  Define  
the braided Dunkl operators  $\rd \nabla_1,\ldots,\rd \nabla_n\colon
S_\qsub(V)\to S_\qsub(V)$ attached to $W$ by 
\begin{equation}
\label{eq:dunkl factored}
\rd \nabla_i =\rd \partial_i +\sum_{\varepsilon\in {\mathcal
    C}_i\setminus \{1\}}\frac{c_{i,\varepsilon}}{1-\varepsilon}\cdot
\frac{1}{x_i}(1-t_i^{(\varepsilon)}).
\end{equation}
\begin{proposition} 

\label{pr:dunkl factored}
The braided Dunkl operators $\rd \nabla_1, \ldots, \rd \nabla_n$  satisfy for all $i,j$:
$$
\rd  \nabla_i x_j-q_{ji} x_j\rd  \nabla_i =
\delta_{ij}(1+\sum_{\varepsilon\in {\mathcal C}_i\setminus
  \{1\}}c_{i,\varepsilon}t_i^{(\varepsilon)}),
\qquad
\rd \nabla_i\rd \nabla_j = q_{ij}\rd \nabla_j\rd \nabla_i \ .
$$
\end{proposition}

In general, braided Dunkl operators attached to a direct product
$W_1\times \dots\times W_m$ of groups are obtained from
Theorem~\ref{thm:brprod}. For each group $W_k$, one writes down either 
commuting Dunkl operators for complex reflection groups \cite{DO} or
anticommuting operators $\rd\nabla_i$ as above. The differential parts
of these operators now become braided derivatives $\rd\partial_i$ 
of $S_\qsub(V)$, an altogether they form a complete list of braided
Dunkl operators for $W$.

Our next result deals with universal embedding of each
braided Cherednik algebra in a \emph{modified Heisenberg double}.
This embedding is crucial in the proof of Theorem~\ref{th:braided
  dunkl} and Theorem~\ref{th:braided Verma}. Besides that, it leads to
new interesting representations of braided Cherednik algebras (see
below). 

A modified Heisenberg double $\mathcal A$ 
is an algebra attached to two Yetter\dash
Drinfeld module structures on the same module $Y$ over a group $W$. It
has triangular decomposition 
$\mathcal A = U^- \tensor \CC W \tensor U^+$, where 
%\begin{equation}
%\label{eq:U-U+}
$$
U^-=T(Y)/\langle \wedge^2_{\Psi_1} Y 
\cap \wedge^2_{\Psi_2}Y\rangle,
\qquad
U^+=T(Y^*)/\langle \wedge^2_{\Psi_1^*} Y^* 
\cap \wedge^2_{\Psi_2^*}Y^*\rangle,
$$
%\end{equation}
where $\wedge^2_\Psi X = \ker(\id_{X\tensor X}+\Psi)$ for $\Psi\in
\End(X\tensor X)$ on any vector space $X$, and $\Psi_1$, $\Psi_2$ are
braidings on $Y$ induced by the two Yetter\dash Drinfeld structures
over $W$. 
%The relevant terminology and machinery is introduced in
%Section~\ref{sect:yd}; in particular, 
%the existence of such algebra with triangular decomposition
%follows by applying a general operation $\diamond$ to two Heisenberg
%quadratic doubles. 

For a braided Cherednik algebra $\rd  {\mathcal H}(W)=S_\qsub(V)\tensor \CC  W
\tensor S_\qsub(V)$, there is an extension 
$\widetilde W=W\cdot \Gamma\subset \GL(V)$ of $W$ 
by means of an Abelian group $\Gamma\subset
\GL(V)$. To this data 
we associate a \emph{$\qq$\dash reflections module} 
$Y$ over $\widetilde W$ with two Yetter\dash
Drinfeld structures, hence a modified Heisenberg double 
$\mathcal A(\widetilde W) = U^-\tensor \CC \widetilde W \tensor U^+$. 
\begin{theorem} 
\label{th:emb0}
In the above setup, there exists an injective algebra homomorphism
$$
\varphi \colon \rd  {\mathcal H}(W) \to \mathcal A(\widetilde W)
$$
such that $\varphi|_W$ is the natural inclusion of $W$ in $\widetilde
W$, $\varphi(V)\subset Y$ and $\varphi(V^*) \subset \Gamma \cdot
Y^*$. 
\end{theorem}

The embedding $\varphi\colon\rd  {\mathcal H}(W)\hookrightarrow {\mathcal
  A}(\widetilde W)$ generalises our earlier result
\cite[Theorem
  7.26]{BB}, where we constructed such  embeddings for all rational
Cherednik algebras. 
This way we can obtain new
representations of $\rd  {\mathcal H}(W)$ in ${\mathcal A}(\widetilde W)$
 or in the Verma-type ${\mathcal A}(\widetilde W)$\dash module $U^-$. 

The quadratic algebra $U^-$ arising from
Theorem~\ref{th:emb0} is itself of great interest. In \cite{BK}, 
Anatol Kirillov and the first author show  that    
when $\rd  {\mathcal H}(W)$ is a rational Cherednik algebra,  
the defining relations in $U^-$ are generalised 
classical Yang\dash Baxter equations. 
In particular, if $W=\Symm_n$, $U^-$ coincides with
the triangular enveloping algebra 
$U(\mathrm{tr}_n)$ of Bartholdi-Enriquez-Etingof-Rains \cite{BEER}; moreover, $U(\mathrm{tr}_n)$ surjects onto the Fomin-Kirillov quadratic algebra $\mathcal{E}_n$ from \cite{FK}, which is relevant for embeddings of rational Cherednik algebras (see \cite[Example 7.24]{BB}).
%(see e.g.,\ \cite[Examples 7.23 and 7.27]{BB} for more details). 
%
It is also quite surprising that when $\rd  {\mathcal H}(W)$ is a
negative braided Cherednik algebra from Theorem~\ref{th:PBW anti-commutating}
with $W=B_n^+$, then the image $\varphi(S_{\minusone}(V))$ which is a
subalgebra of $U^-$ by Theorem~\ref{th:emb0}, coincides 
with what Majid called the \emph{algebra of flat
  connections} with constant coefficients 
in the noncommutative differential geometry of the
symmetric group \cite{Mnoncomm}.

To conclude the Introduction, we list 
relevant open problems and new directions of study.

\subsection*{Degenerate $q$-Hecke algebras} Here, the problem is two-fold:

\noindent\textbf{Problem 1.} (a) Given a $\qq$-symmetric algebra $S_\qq(\tilde
V)$, find all finite groups $W\subset GL(V)$ such that the $W$-action
on $V$ extends to the $W$-action on  $S_\qq(\tilde V)$ by algebra
automorphisms.  

(b) For each such $W$, find all flat deformations of the semidirect
product algebra $S_\qq(\tilde V)\rtimes \CC W$.  

Here, we solve Problem 1(a) in the case when all $q_{ij}\ne 1$ for 
$i\ne j$ and under the assumption that $W$ also acts on  $S_\qq(\tilde
V)^*$ by algebra automorphisms (Section \ref{subsect:preserves}). In fact,
the above groups $W_{{\mathcal C}, {\mathcal C}'}$ form the most
important class of solutions(when all $q_{ij}=-1$ for $j\ne i$). 

Each $\qq$-Cherednik algebra and braided Cherednik algebra is a
solution to Problem 1(b) in the case when $\tilde V=V\oplus
V^*$. In our forthcoming paper \cite{BB3} we construct more solutions
to the problem. 

\subsection*{Representations of braided Cherednik algebras}
Similarly to the ordinary (rational) Cherednik algebras, 
one defines the category
${\mathcal O}$ for each braided Cherednik algebra $\rd  {\mathcal
  H}$. The following natural problem emerges:

\noindent\textbf{Problem 2.} For each braided Cherednik algebra $\rd  {\mathcal
  H}=\rd  {\mathcal H}_c$, describe the category ${\mathcal O}$. In
particular, find all values of parameters $c$ such that ${\mathcal O}$
contains finite-dimensional objects.  
  
Even though  ${\mathcal O}$ is not a tensor category, in addition to
the Verma modules, it contains a number of interesting objects:
$U^-\otimes \rho$, where  $U^-$ is the ``generalised $r$-matrix
algebra'' from Theorem \ref{th:emb0} and $\rho$ is any representation
of $W$. As we mentioned above, if $W=S_n$ it is known from \cite{BK}
and \cite{BEER} that the quadratic algebra $U^-$ is Koszul. We expect
this phenomenon to persist in general, therefore, having an $\rd
{\mathcal H}$-module structure on $U^-\otimes \rho$ and on $U^-$
itself is beneficial for understanding this quadratic algebra. 

We plan to study finite\dash dimensional quotients of the Verma
module $M_{triv}=S_\qq(V)$ for $\rd  {\mathcal H}$ in a separate
paper. We expect that for negative braided Cherednik algebras the
answer can be given along the lines of \cite{BEG} and \cite{VV}. And,
according to Remark \ref{rem:degenerate}, the degenerate version of
$\rd  {\mathcal H}$ has a number of finite-dimensional modules that
can be studied along the lines of \cite{Go}. 

However,  when $\rd  {\mathcal H}$ is a braided tensor product (of
  negative braided or ordinary Cherednik algebras) as in
  \eqref{eq:braided factorization1}, the representation category of
  $\rd  {\mathcal H}$ is not at all determined by those of the tensor
  factors. For instance, by varying the matrix  $\qq$ and parameters
  $c_{i,\varepsilon}$ in \eqref{eq:dunkl factored} and Proposition
  \ref{pr:dunkl factored}, one can expect new interesting submodules
  of the Verma module $S_\qsub(V)$ even when $W=\Symm_2\times
  \Symm_2\times \cdots \times 
  \Symm_2$. Another ``degree of freedom'' 
in representations of such factored $\rd
       {\mathcal H}$ is a choice of the field of definition ${\mathbb 
	 K}\subset \CC$ containing all $q_{ij}$, say, ${\mathbb
	 K}={\mathbb Q}(q_{ij}|i,j=1,\ldots,n)$. Then, in the
       assumption that all $q_{ij}$ are roots of unity, i.e.,
       ${\mathbb K}$ is a cyclotomic extension of ${\mathbb Q}$, there
       exist finite\dash dimensional quotients ${\mathcal B}_\qsub$ of
       $S_\qsub(V)$ and we expect that some of these ${\mathcal
	 B}_\qsub$  are, in fact, representations of $\rd {\mathcal
  H}$. It
       follows from the famous Merkurjev-Suslin theorem that
       essentially all central simple algebras over ${\mathbb K}$ are  
%tensor products of symbol algebras $\{a,b\}_m=\langle x_1,x_2|
%x_1x_2=qx_2x_1, x_1^m=a,x_2^m=b\rangle$ for $q$ being a primitive
%$m$-th root of unity in ${\mathbb K}^\times$ 
%and $a,b\in {\mathbb K}^\times$, i.e., essentially all central simple algebras over ${\mathbb K}$ are 
simple finite-dimensional quotients of various $S_\qsub(V)$ so that an
$\rd {\mathcal H}$-module structure on them would be of interest in
$K$-theory.

%%%%%%%%%%%%%%%%%%%%%%%%%%%%%%%%%%%%%%%%%%%%%%%%%%%%%%%%%%%%%%%%

%\medskip

\subsection*{Acknowledgments}
 
We thank the organisers of the Workshop on Cherednik algebras
at ICMS, Edinburgh, where this paper was started. 
The second author (A.B.) expresses his gratitude to the 
Mathematics Institute, University of Warwick for hospitality and support  
during his visit in Summer 2007. 
The authors are grateful to Weiqiang Wang for bringing to our
attention the remarkable joint papers \cite{KW1,KW2} with  Ta Khongsap.

%%%
\section{Quadratic doubles}
\label{sect:qd}

In this section we introduce quadratic doubles as  a sub-class
of \emph{braided doubles} over bialgebras (introduced and
studied in our earlier paper \cite{BB}) and %we give an overview of quadratic doubles 
%\ref{sect:qd} and \ref{sect:yd} 
%quadratic doubles over group algebras, together with 
present new results related specifically to quadratic doubles
over group algebras. It is an open question, if (and how) results such as 
Proposition \ref{prop:operations} and Theorem \ref{thm:perfect}
can be extended to doubles over arbitrary bialgebras or Hopf
algebras. 

\subsection{Triangular decomposition and braided doubles}

Triangular decomposition of an associative algebra is defined as follows. 

\begin{definition}
An algebra $A$ admits  \emph{triangular decomposition}
$A=U^-\tensor U^0 \tensor U^+$, if $U^0$, $U^\pm$ are subalgebras in
$A$ such that  $U^-U^0$ and $U^0 U^+$ are also subalgebras in $A$,
and the vector space map 
$$ 
U^-\tensor U^0 \tensor U^+\to A, \qquad
u^-\tensor u^0 \tensor u^+ \mapsto u^- u^0 u^+\in A, 
$$
is bijective.  
\end{definition}

Let $V$ be a finite\dash dimensional module over a group $W$.  
To a linear map $\beta\colon V^*\tensor V \to \CC W$ and two 
subspaces
$R^-\subset T^{>0}(V)$, $R^+ \subset T^{>0}(V^*)$
we associate the algebra 
$$ 
A_\beta(R^-,R^+) = \frac{ T(V\oplus V^*)\lcprod \CC W }{\lgen \, R^-,\ R^+,\
\{f\tensor v-v\tensor f-\beta(f\tensor v)
\mid f\in V^*, v\in V\}\, \rgen}^{\vphantom{T}}\  .
$$
Here the symbol $\lcprod$ is used to denote  a semidirect product. 
(If $A$ is a a $W$\dash module
algebra, $A\lcprod \CC W$ is the algebra with underlying vector space
$A\tensor \CC W$ and multiplication $(a\tensor w)(a'\tensor
w')=aw(a')\tensor ww'$, where $a,a'\in A$ and $w,w'\in W$.)
The angular brackets $\lgen \ \rgen$
denote the two-sided ideal with given generators. 

\begin{definition}

The algebra $A_\beta(R^-,R^+)$ is a \emph{braided double}, if it
has triangular decomposition
$$
A_\beta(R^-,R^+) \ \cong \ T(V) / \lgen R^-\rgen \, \tensor \,\CC W \,\tensor\, T(V^*)/\lgen R^+ \rgen.
$$
\end{definition}

\begin{remark}
%1. 
The algebra $A_\beta(R^-,R^+)$ may either be a proper 
%is only a
quotient of the vector space on the right, or 
%
%2. Suppose that 
%$A_\beta(R^-,R^+)$ 
be a braided double. In the latter case, $\lgen R^-\rgen$ 
%and $\lgen R^+ \rgen$ are 
is automatically a $W$\dash invariant ideal in the tensor algebra
$T(V)$, and the subalgebra $T(V) / \lgen R^-\rgen \tensor \CC W$ 
%$T(V^*)/\lgen
%   R^+ \rgen \tensor \CC W$ are 
is isomorphic to 
%the respective semidirect products  
$T(V) / \lgen R^-\rgen \lcprod \CC W$; similarly for $R^+$.
%,  $\CC W
%   \rcprod T(V^*)/\lgen
%   R^+ \rgen$.
\end{remark}
  
To understand braided doubles, one would like to study the locus of
parameters $(\beta,R^-,R^+)$ such that the algebra 
$A_\beta(R^-,R^+)$ has
triangular decomposition. The first major step is to determine for
which $\beta$ braided doubles 
of the form $A_\beta(R^-,R^+)$ exist. 
We say that $\beta\colon V^*\tensor V \to \CC W$ is a \emph{$W$\dash
equivariant map} if $\beta$ is a 
$W$\dash homomorphism with respect to the standard  
diagonal $W$\dash action on $V^*\tensor V$ and the $W$\dash
action on $\CC W$ by \emph{conjugation}.
\begin{theorem}
\label{thm:equivar}
Let $V$ be a finite\dash dimensional module over a group $W$ and 
 $\beta\colon V^*\tensor V \to \CC W$ be a linear map. 
The algebra $A_\beta(0,0)$ is a braided double, if
 and only if $\beta$ is a $W$\dash equivariant map. 
\end{theorem}
\begin{proof}
To prove the `only if' part, 
pick any $f\in V^*$, $v\in V$ and $w\in W$. 
Using the
%following 
relations $wvw^{-1} = w(v)$ and $wfw^{-1}= w(f)$ in the algebra
%hold in the smash product 
$T(V\oplus V^*)\lcprod\CC W$, write
%$$
%     wvw^{-1} = w(v) , \qquad wfw^{-1}= w(f)\,.
%$$ 
%Write 
$\beta(w(f),w(v))=(wfw^{-1})(wvw^{-1})-(wvw^{-1})(wfw^{-1})$,
which is equal to $w(fv-vf)w^{-1}=w\beta(f,v)w^{-1}$. 
Thus, %the relation 
$\beta(w(f),w(v))= w\beta(f,v)w^{-1}$ holds in 
$A_\beta(0,0)$. Both sides of this relation lie in the group
algebra $\CC W$ which embeds injectively
in $A_\beta(0,0)$ because of the triangular decomposition of
$A_\beta(0,0)$. Hence the relation holds in $\CC W$ and $\beta$ is
$W$\dash equivariant. The (more difficult) `if' part is proved in
\cite[Theorem 3.3]{BB}; the 
key point here is that  the $W$\dash equivariance may be interpreted
as the Yetter\dash Drinfeld condition for modules
over a group algebra. 
\end{proof}

A braided double of the form $A_\beta(0,0)$ is called a \emph{free
  braided double} and denoted $\widetilde A_\beta$.
%From 
The proof of the Theorem implies 
%it is easy to deduce the following
that
%\begin{corollary}
if $A_\beta(R^-,R^+)$ is a braided double, then $\beta$ is $W$\dash
equivariant and $A_\beta(R^-,R^+)$ is a quotient of $\widetilde
A_\beta$. The quotient map in question is a morphism in the category
of braided doubles:

\begin{definition}[The category $\dcat_W$]
Denote by $\dcat_\beta(V)$ the set of braided doubles of the form
$A_\beta(R^-,R^+)$. 
We introduce the category $\dcat_W$ such that 
%, the objects
%of which are braided doubles over $\CC W$:
$$
\Ob\ \dcat_W = \bigcup_{V,\beta} \dcat_\beta(V),
$$
where the union is taken over all finite\dash dimensional $W$\dash modules
$V$ and all $W$\dash equivariant maps $\beta\colon V^*\tensor V \to
\CC W$. If $A\in \dcat_\beta(U)$ and $B\in
\dcat_\gamma(V)$, a morphism $\varphi\colon A\to B$ in $\dcat_W$ 
is an algebra map
such that $\varphi(U)\subset V$, $\varphi(U^*)\subset V^*$ and
$\varphi|_W=\id_W$. 
\end{definition}

Clearly, $\varphi$ is uniquely determined by the two
$W$\dash module maps $\mu=\varphi|_U$ and $\nu=\varphi|_{U^*}$.
However, not every 
pair of $W$\dash module maps $U\xrightarrow{\mu}
V$, $U^*\xrightarrow{\nu}V^*$  extends to an algebra homomorphism
$A\to B$. 
For example, zero maps $\mu=\nu=0$ do not extend to a morphism 
between $A\in \dcat_\beta(U)$ and $B\in
\dcat_\gamma(V)$, unless $\beta=0$.

%\end{corollary}
%
%It also follows from the $W$\dash equivariance of $\beta$ that 

Observe also that a braided double $A_\beta(R^-,R^+)$ is a $W$\dash
module algebra  (where the action of $W$ on generators $w\in W$ of
$A_\beta(R^-,R^+)$ is by conjugation).

Using Lemma 4.4 \cite{BB}, one obtains a way to construct braided
doubles in terms  of $\widetilde A_\beta$: 

\begin{proposition}
\label{prop:way}
 Let $R^-\subset T^{>0}(V)$ and $R^+\subset T^{>0}(V^*)$ be $W$\dash
 submodules such that $[R^+,V]=[V^*,R^-]=0$ in the free braided double
 $\widetilde  A_\beta$. Then $A_\beta(R^-,R^+)$ is  a braided double. 
\qed
\end{proposition}
\begin{remark}

%%In a braided double $A=A_\beta(R^-,R^+)$, one can vary $R^-$ leaving
%%$\lgen R^-\rgen$ fixed (same for $R^+$), thus hoping to obtain an
%%``optimal'' presentation of $A$ with $R^\pm$ satisfying the zero
%%commutator condition in .
Based on this result, it is natural to expect that all braided doubles can be obtained as quotients of free braided doubles by the zero 
commutator condition in Proposition~\ref{prop:way}.
For example, enveloping algebras
$U(\mathfrak g)$ and $U_q(\mathfrak g)$ of a semisimple Lie algebra
$\mathfrak g$ have such presentation (with $R^\pm$ being
the Serre relations). Finding such ``optimal'' presentation for braided
Heisenberg doubles \cite[5.3]{BB} would imply interesting
results on the structure of Nichols algebras. 
This optimal presentation is available for the main object of this section -- quadratic doubles  
(see Theorem~\ref{thm:relations} below).
\end{remark}

%\begin{proof}
%Since the ideals $\lgen R^-\rgen \subset T(V)$ and $\lgen R^+\rgen
%\subset T(V^*)$  are $W$\dash invariant, the claim follows directly
%from \cite[Lemma 4.4]{BB}. 
%\end{proof}

%Observe also that if $\beta$ is $W$\dash equivariant, the group $W$
%acts covariantly on any algebra $A_\beta(R^-,R^+)$, 
%where the action of $W$ on generators $w\in W$ of the algebra is 
%by conjugation. 
%Theorem~\ref{thm:equivar} thus implies  
%\begin{corollary}
%\label{cor:mod-alg}
%Let $W$ be a group. 
%Any braided double over $\field W$ is a $W$\dash module algebra. 
%\qed
%\end{corollary}

\subsection{Quadratic doubles}

%Most relevant for us is the case when $R^\pm$ are
%spaces of quadratic relations, i.e., 
A braided double $A_\beta(R^-,R^+)$ in
$\dcat_\beta(V)$ is called a \emph{quadratic double}, if
$R^-\subset V\tensor V$ and
$R^+\subset V^*\tensor V^*$.  
%if $R^\pm$ are
%quadratic. 
Our original motivating example of this is rational
Cherednik algebra; free braided doubles
are quadratric too. 
%(On the other hand, important braided
%doubles such as enveloping algebras
%$U(\mathfrak g)$ and $U_q(\mathfrak g)$ of a semisimple Lie algebra
%$\mathfrak g$, are not quadratic.
%See \cite{BB} for more about that line of research.) 
We denote by $\qdcat_\beta(V)$ the set of quadratic doubles in
$\dcat_\beta(V)$ and by $\qdcat_W$ the category of quadratic doubles
over $\CC W$ (a full subcategory of $\dcat_W$). 

%The next Theorem is a practical tool which in many cases allows us to
%determine possible quadratic relations $R^\pm$ in a double of the  form
%$A_\beta(R^-,R^+)$ with given $\beta$. 
%Here and below, square brackets denote the usual
%commutator $[a,b]=ab-ba$ in an associative algebra. 
 
\begin{theorem}
\label{thm:relations}
Let $\beta\colon V^*\tensor V\to \CC W$ be a $W$\dash equivariant
map. Then 
$A_\beta(R^-,R^+)$ is a quadratic double if and only if
$R^-\subset V\tensor V$, $R^+\subset V^*\tensor V^*$ are 
$W$\dash submodules and 
$$
     [R^+,V]=0, \qquad [V^*,R^-]=0 \qquad \text{in the free double}\ \widetilde A_\beta\,.
$$
\end{theorem}
\begin{proof}
If 
%the algebra 
$A_\beta(R^-,R^+)$ is a quadratic double, 
%it has triangular
%decomposition  of the form $A_\beta(R^-,R^+)\cong T(V)/\lgen R^- \rgen
%\tensor \CC W \tensor T(V^*)/\lgen R^+\rgen$. It follows that 
the ideal $\lgen R^- \rgen$ of $T(V)$ is $W$\dash invariant, 
hence so is its quadratic part $R^-$. The same %argument 
applies to $R^+$. Furthermore, the relations
in the free double $  \widetilde A_\beta$ imply that 
the commutator $[V^*,R^-]$ is a subspace of $V \tensor \CC W$ which
must obviously be in the kernel of the quotient map $\widetilde
A_\beta \to A_\beta(R^-,R^+)$. The quotient map has no kernel in degrees less
than $2$ with respect to generators from $V$, thus $[V^*,R^-]=0$ in 
$\widetilde A_\beta$. The same argument applies to $[R^+,V]$ and thus
establishes the `only if' statement. The `if' statement follows by
Proposition~\ref{prop:way}.  
\end{proof}

We denote 
$$
R^-_{\beta\, \mathrm{max}}= \{ r^- \in V\tensor V \mid [V^*,r^-]=0 \
\text{in} \ \widetilde A_\beta\}\,,
\qquad 
R^+_{\beta\, \mathrm{max}} = \{r^+\in V^*\tensor V^* \mid [r^+,V]=0\ \text{in} \ \widetilde A_\beta\}. 
$$
It easy to see that as long as
$\beta\colon V^*\tensor V \to \CC W$ is $W$\dash equivariant (as
above), $R^\pm_{\beta\, \mathrm{max}}$ are $W$\dash invariant subspaces. This
observation is useful in the following

\begin{corollary}
Let $\beta\colon V^*\tensor V \to \CC W$ be a $W$\dash equivariant
map. The algebras 
$\widetilde A_\beta := A_\beta(0,0)$ and $\mindouble A_\beta:=
A_\beta(R^-_{\beta\, \mathrm{max}},R^+_{\beta\, \mathrm{max}})$  are
quadratic doubles in $\qdcat_\beta(V)$. 
For any quadratic double $A\in \qdcat_\beta(V)$ there are  quotient maps 
$\widetilde A_\beta \twoheadrightarrow  A  \twoheadrightarrow \mindouble
A_\beta$
in $\qdcat_W$.
\qed
\end{corollary}
\begin{definition}
The quadratic double
$
\mindouble A_\beta \ \cong\  T(V)/\lgen R^-_{\beta\, \mathrm{max}} \rgen 
             \, \tensor\, \CC W \,\tensor\, 
              T(V^*)/\lgen R^+_{\beta\, \mathrm{max}}\rgen
$ 
in $\qdcat_\beta(V)$ is called 
the \emph{minimal quadratic double} with parameter
$\beta\in\{W$\dash equivariant maps $V^*\tensor V \to \CC W\}$.  
\end{definition}
%Among 
%quadratic doubles in 
Of the objects  $\qdcat_\beta(V)$, it is the minimal quadratic
double $\mindouble A_\beta$ that is most 
interesting algebraically.
The quadratic  relations in $\mindouble A_\beta$  
are given implicitly as kernels of certain linear operators (see
Lemma~\ref{lem:T} below). The
central problem in the theory of quadratic doubles is two-fold:
\begin{problem}
\label{prob:qd1}
Let $V$ be a finite\dash dimensional module for a group $W$. 
Given a $W$\dash equivariant map $\beta \colon V^*\tensor V \to \CC 
W$, define the algebra  $\mindouble A_\beta$ explicitly by
generators and relations (i.e., find an explicit description of 
$R^\pm_{\beta\, \mathrm{max}}$).  
\end{problem}
\begin{problem}
\label{prob:qd2} 
Given $W$\dash submodules $R^-\subset V\tensor V$ and 
   $R^+\subset V^*\tensor V^*$, find all maps $\beta\colon V^*\tensor
   V \to \CC W$ such that $A_\beta(R^-,R^+)$ is a quadratic double.  
\end{problem}

Problem~\ref{prob:qd2}
is in fact a deformation question.
Regard $A_\beta(R^-,R^+)$ as a deformation, with
parameter $\beta$, of the algebra $A_0(R^-,R^+)$; the latter 
is a quadratic double by Theorem~\ref{thm:relations}. One needs 
to find the values of $\beta$ for which the deformation is flat (the
flatness locus).

\begin{example}
\label{ex:cherEG}
When $W\subset \GL(V)$, $R^-=\wedge^2 V\subset V\tensor V$ and 
$R^+=\wedge^2 V^*\subset V^*\tensor V^*$, the solution to 
Problem~\ref{prob:qd2} is given by rational Cherednik
algebras $A_\beta(R^-,R^+)$ with 
$$
\beta(\xi\tensor v) = \langle v, \xi \rangle 
  + \sum_s c_s \langle v, \alpha_s^\vee \rangle \langle \alpha_s, \xi
   \rangle s\, 
$$
for $\xi\in V^*$, $v\in V$; cf.\ \cite{EG}. The sum is taken over all
complex reflections $s\in W$, the parameters 
$c=\{c_s\}$, $c_s\in \CC$ satisfy $c_{wsw^{-1}}=c_s$ for all $w\in W$, and
$\alpha_s\in V$, $\alpha_s^\vee\in V^*$ 
is the root-coroot pair for the complex reflection $s$, meaning that
$s(v)=v-\langle v, \alpha_s^\vee\rangle \alpha_s$ for all $v\in V$.
Here $\langle v, \xi \rangle$ can be any
$W$\dash invariant pairing between $V$ and $V^*$.
If it is the standard evaluation pairing, denote the corresponding
rational Cherednik algebra by $H_c(W)$, whereas if $\langle v, \xi
\rangle=0$, denote the corresponding algebra by $H_{0,c}(W)$.
\end{example}

\subsection{Operations $\diamond$ and $\star$ on quadratic doubles}

Recall that the parameter $\beta$ in a quadratic double
$A_\beta(R^-,R^+)$ belongs to the space of $W$\dash equivariant linear
maps from $V^*\tensor V$ to $\CC W$.  
Let us now observe that this parameter space has the structure of an
algebra. Write $\beta$ in the form
$$
   \beta(f\tensor v) = \sum_{w\in W} \langle L_w(v), f\rangle w,
\qquad f\in V^*, \ v\in V,
$$
where $L_w\in \End(V)$ are zero for all but finitely many $w\in
W$. We identify $\beta$ with the element $\sum_{w\in W} \delta_w \tensor L_w$
of the algebra 
$$
     (\CC(W)_0 \tensor \End(V))^W. 
$$
Here $\CC(W)_0$ is the algebra, with respect to pointwise
multiplication, of complex\dash valued functions on $W$
with finite support. It is spanned by delta\dash functions $\delta_w$,
$w\in W$. The action of $W$ on $\CC(W)_0$ is by conjugation:
$w(\delta_\sigma)=\delta_{w\sigma w^{-1}}$. The $\tensor$ is the
standard tensor product of algebras, where the tensorands commute. 
The algebra $(\CC(W)_0 \tensor \End(V))^W$ of parameters contains
an identity if and only if the group $W$ is finite. 

Let  $\beta=\sum_w \delta_w \tensor L_w$,
$\gamma = \sum_w \delta_w \tensor M_w$ be elements of $(\CC(W)_0
\tensor \End(V))^W$. We observe that their sum and product in the algebra of
parameters are rewritten as linear maps from $V^*\tensor V$ to $\CC W$
as follows:
$$
   (\beta+\gamma)(f\tensor v) =  \sum_{w\in W} \langle (L_w+M_w)(v),
   f\rangle w, 
   \qquad
   (\beta\gamma)(f\tensor v) =  \sum_{w\in W} \langle (L_w M_w)(v),
   f\rangle w\,. 
$$
Let 
$$
     \qdcat(V) := \bigsqcup_{\beta\in (\CC(W)_0 \tensor \End(V))^W
     } \qdcat_\beta(V)
$$
be the set of all quadratic doubles of the $W$\dash module $V$. 
We will now see how the above sum and product can be ``lifted'' from
the algebra of 
parameters to  $\qdcat(V)$, to yield two
operations, $\diamond$ and $\star$. 

\begin{definition}
\label{def:stardiamond}
Let $A=A_\beta(R^-,R^+)$ and $B=A_\gamma(S^-,S^+)$ be quadratic doubles 
in $\qdcat(V)$. Denote
$$
A\diamond B = A_{\beta+\gamma}(R^-\cap S^-,R^+\cap S^+), 
\qquad
A\star B = A_{\beta\gamma}(S^-,R^+).
$$
\end{definition}
\begin{proposition}
\label{prop:operations}
If $A$ and $B$ are quadratic doubles in $\qdcat(V)$, then $A\diamond
B$ and $A\star B$ are also quadratic doubles in $\qdcat(V)$. 
\end{proposition}
The proposition will follow from a technical
\begin{lemma}
\label{lem:T}
Let $\beta=\sum_w \delta_w \tensor L_w$ be an element of the parameter
algebra $(\CC_0(W)\tensor \End(V))^W$. 
The quadratic relations in the minimal quadratic double $A_\beta$
are given by 
$$
R^-_{\beta\, \mathrm{max}}=\cap_{w\in W}\ker T_{w,\beta}^-
, 
\qquad 
R^+_{\beta\, \mathrm{max}} =\cap_{w\in W}\ker T_{w,\beta}^+ 
,
$$
where $T_{w,\beta}^-\in \End(V\tensor V)$ and $T_{w,\beta}^+\in
\End(V^*\tensor V^*)$ are defined by 
$$
  T_{w,\beta}^-(u\tensor v) = (L_w\tensor \id_V)(u\tensor w(v)+v\tensor u), 
\qquad
  T_{w,\beta}^+(f\tensor g) = (\id_{V^*}\tensor L_w^*)(w^{-1}(f)\tensor g +
  g\tensor f).
$$
\end{lemma}
\begin{proof}
Recall that $R^-_{\beta\, \mathrm{max}}$ is defined, following
Theorem~\ref{thm:relations}, as the space of quadratic tensors in
$V\tensor V$ that commute, in the free double $\widetilde A_\beta$,
with all elements of $V^*$. 
By the Leibniz rule, the commutator of $f\in V^*$ with $u\tensor v\in
V\tensor V$ in $\widetilde A_\beta$ is 
\begin{align*}
   \sum_{w\in W} \langle L_w(u),f\rangle w v + u\langle
   L_w(v),f\rangle w 
   & = \sum_{w\in W} \Bigl( \langle L_w(u),f\rangle w(v) + 
     \langle L_w(v),f\rangle u \Bigr) \tensor w 
\\
   & =  \sum_{w\in W} (\langle \cdot,f\rangle \tensor \id_V)
   T_{w,\beta}^-(u\tensor v)\tensor w,  
\end{align*}
hence indeed $R^-_{\beta\, \mathrm{max}}=\cap_{w\in W}\ker T^-_{w,\beta}$. The argument for
$R^-_{\beta\, \mathrm{max}}$ is similar. 
\end{proof}
\begin{proof}[Proof of Proposition~\ref{prop:operations}]
To establish that $A\diamond B$ is a quadratic double, we need to show
that $R^\pm\cap S^\pm$ is a $W$\dash submodule of 
$R^\pm_{\beta+\gamma\, \mathrm{max}}$ and to apply
Theorem~\ref{thm:relations}. 
But clearly $T_{w,\beta+\gamma}^\pm=T_{w,\beta}^\pm+T_{w,\gamma}^\pm$, 
thus $\ker T_{w,\beta+\gamma}^\pm$ contains 
the intersection of $\ker T_{w,\beta}^\pm$ and $\ker T_{w,\gamma}^\pm$,
which contains $R^\pm \cap S^\pm$. The latter is a $W$\dash
submodule as an intersection of $W$\dash submodules.  
In a similar fashion, to show that $A\star B$ is a quadratic double,
we need to check that $S^-\subset\ker T_{w,\beta\gamma}^-$ and 
$R^+\subset\ker T_{w,\beta\gamma}^+$ for all $w\in W$. 
Write $\beta=\sum_{w\in W} \delta_w \tensor L_w$ and
$\gamma=\sum_{w\in W} \delta_w \tensor M_w$. 
Observe that 
$ T_{w,\beta\gamma}^- = (L_w\tensor \id_V)T_{w,\gamma}^-$, therefore
the kernel of $ T_{w,\beta\gamma}^-$ contains that of
$T_{w,\gamma}^-$, which contains $S^-$. Furthermore, 
$ T_{w,\beta\gamma}^+ = (\id_{V^*}\tensor M_w^*)T_{w,\beta}^+$, hence
its kernel contains that of $T_{w,\beta}^+$, which contains $R^+$. 
\end{proof}

\begin{remark}
The two operations 
$\diamond, \star\colon \qdcat(V)\times \qdcat(V) \to \qdcat(V)$ 
satisfy the following axioms:
\begin{itemize}
\item[]
$(A\diamond B)\diamond C = A\diamond (B\diamond C)$, 
\quad
$(A\star B)\star C = A\star (B\star C)$;
\item[]
$A\diamond B=B\diamond A$;
\item[]
$A_0 \diamond A = A \diamond A_0 = A$;
\item[]
$A\star(B\diamond C)=(A\star B)\diamond (A\star C)$,
\quad
$(A\diamond B)\star C = (A\star C) \diamond (B\star C)$,
\end{itemize} 
where $A_0\cong V\tensor \CC W \tensor V^*$ is the minimal quadratic
double corresponding to $\beta=0$ (the ``smallest possible'' quadratic
double). 
This is a subset of the semiring axioms, however, note that there 
is no zero or identity element with respect to $\star$. 

Warning: the operations $\diamond$ and $\star$ do not preserve the
minimality of quadratic doubles: $A\diamond B$ and $A\star B$ may not
be minimal even if $A$, $B$ are both minimal. 
\end{remark}

We will now see  how the operation $\diamond$ 
``behaves'' with respect to morphisms in $\qdcat_W$.

\begin{proposition}
\label{prop:diamond}
Let $U$, $V$ be two finite\dash dimensional $W$\dash modules, and
 assume that $A,B\in \qdcat(U)$, $A',B'\in \qdcat(V)$ are quadratic
 doubles. If a pair $U\xrightarrow{\mu} V$, $U^* \xrightarrow{\nu}
 V^*$  of $W$\dash module maps 
extends to
a morphism $\varphi\colon A\to A'$ and to a morphism $\psi\colon B\to
 B'$ in $\qdcat_W$, 
then the same maps $\mu$, $\nu$ extend to a morphism 
$$
\varphi\diamond\psi \colon A\diamond B \to A'\diamond B'.
$$  
\end{proposition}
\begin{proof}
Let $A=A_\beta(R^-,R^+)$, $B=A_\gamma(S^-,S^+)$ 
where $\beta,\gamma\colon U^*\tensor U \to \CC W$ are
$W$\dash equivariant maps, $R^-,S^-\subset U\tensor U$ and
$R^+,S^+\subset U^*\tensor U^*$.  Let  
$A'=A_{\beta'}({R'}^-,{R'}^+)$, $B'=A_{\gamma'}({S'}^-,{S'}^+)$,
similarly to $A$, $B$ but with $U$ replaced with $V$.
Looking at the relations in quadratic doubles, we conclude 
that $\mu$, $\nu$ extend to algebra homomorphisms $A\to A'$, $B\to B'$ 
if and only if 
$\beta=\beta'\circ(\nu\tensor \mu)$,
$\gamma=\gamma'\circ(\nu\tensor \mu)$ and
$$
(\mu\tensor \mu)R^-\subset {R'}^-, \quad (\mu\tensor \mu)S^-\subset
{S'}^-, 
\quad
(\nu\tensor \nu)R^+\subset {R'}^+, 
\quad (\nu\tensor \nu)S^+\subset {S'}^+.
$$
But then $\beta+\gamma=(\beta'+\gamma')\circ(\nu\tensor \mu)$ and 
$(\mu\tensor \mu)(R^-\cap S^-)\subset ({R'}^-\cap {S'}^-)$, similarly
for $R^+\cap S^+$. Thus, $\mu$, $\nu$ extend to a morphism $A\diamond
B \to A'\diamond B'$, which we denote by $\varphi\diamond\psi$.
\end{proof}

%\section{Yetter-Drinfeld modules and Heisenberg quadratic doubles}
%\section{Heisenberg quadratic doubles and braided reduction}
%\label{sect:yd}

\subsection{Yetter-Drinfeld modules}

Yetter\dash Drinfeld modules over $W$ provide a family of 
deformation parameters $\beta$, for which the minimal quadratic
doubles $\mindouble A_\beta$ have a nice description and are in a sense
universal, as many quadratic doubles can be realised as their subalgebras 
(see Theorem~\ref{thm:perfect} below). 
Let us recall the definition of a Yetter\dash Drinfeld module. When
the group $W$ is finite, it is the same as a module over the Hopf
algebra $D(W)$, the Drinfeld quantum double of $W$. 

\begin{definition}
A \emph{Yetter-Drinfeld module} for a group $W$ is a $W$\dash module
$Y$ with a grading $Y=\bigoplus\limits_{w\in W} Y_w$, such that  
$\sigma(Y_w)=Y_{\sigma w \sigma^{-1}}$ for all $w,\sigma\in W$. 
\end{definition}

Whenever $Y$ is a Yetter\dash Drinfeld (YD-) module over $W$, we
denote by $|y|$ the $W$\dash degree of homogeneous $y\in Y$. When 
the notation $|\cdot|$ is used in formulas, 
extension from homogeneous elements to all
elements of $Y$ by linearity is implied. 
For example, the Yetter\dash Drinfeld axiom may be written as
$|w(y)|=w|y|w^{-1}$.

Clearly, if $Y$ is finite\dash dimensional, the dual module $Y^*$ is a 
YD-module via 
$Y^*=\bigoplus\limits_{w\in W} (Y^*)_w$ with  
$(Y^*)_w=\mathrm{Hom}_\field(Y_{w^{-1}},\field)$.  
%For a YD-module $Y$ d
Define the linear map 
$\beta_Y\colon Y^*\tensor Y\to \field W$ by 
$$
   \beta_Y(f\tensor v) = 
             \langle v,f\rangle   |v|.
$$
%The following 
It is straightforward to verify that 
%
%\begin{proposition}
the map $\beta_Y$ is $W$\dash equivariant. 
%\qed
%\end{proposition}
%
%It follows that 

\subsection{Heisenberg quadraic doubles}

To each finite\dash dimensional Yetter\dash Drinfeld module $Y$ over $W$ is
therefore associated a minimal quadratic  double 
$\mindouble A_Y :=\ A_{\beta_Y}$,
%\cong T(Y)/\lgen R^-_{\beta_Y\, \mathrm{max}}
%\rgen\tensor \field W 
%\tensor T(Y^*)/\lgen R^+_{\beta_Y\,\mathrm{max}}\rgen\,,
%$$ 
referred to as the \emph{Heisenberg quadratic double} of $Y$.

To describe 
%the algebras $T(Y)/\lgen R^-_{\beta_Y\,
%  \mathrm{max}}\rgen $ and $T(Y^*)/\lgen
%R^+_{\beta_Y\,\mathrm{max}}\rgen$, 
Heisenberg quadratic doubles more explicitly,
recall that the linear map 
$$
\Psi_Y\colon Y\tensor Y \to Y\tensor Y,
\qquad
\Psi_Y(y \tensor z) = |y|(z) \tensor y
$$
is a braiding on $Y$, i.e., a solution to the braid equation; see
\cite[Section 5]{BB}. Viewing $Y^*\tensor Y^*$ as a dual space
to $Y\tensor Y$, denote by  $\Psi_Y^*$ the adjoint map to $\Psi_Y$. 
(This braiding on $Y^*$ is not the same as the braiding $\Psi_{Y^*}$ given by the YD-module
structure on $Y^*$; the two are related via 
$\Psi_Y^*=\tau\circ \Psi_{Y^*} \circ \tau$, where $\tau(x\tensor
y)=y\tensor x$ is the trivial braiding.) Furthermore, 
any braiding $\Psi\in \End(V\tensor V)$ on a vector space $V$ gives
rise to a braided analogue of the symmetric algebra of $V$:
$$
S(V,\Psi) = T(V) / \, \lgen \ker(\id_{V\tensor V} + \Psi) \rgen\, ,
$$
of which $S(V)$ is a particular case corresponding to $\Psi=\tau$.
%with $\tau(x\tensor y)=y\tensor x$ the trivial braiding.
Theorem~5.4 in \cite{BB} implies
\begin{proposition}
\label{prop:heisenberg}
%Let $Y$ be a finite\dash dimensional Yetter\dash Drinfeld module for a
%group $W$.
The 
%minimal
Heisenberg
 quadratic double $\mindouble A_Y$ has
triangular decomposition  
$$
\mindouble A_Y = S(Y,\Psi_Y) \tensor \field W \tensor S (Y^*,\Psi_Y^*).
\qed
$$
\end{proposition}
%We will refer to the minimal quadratic double $\mindouble A_Y$, associated to
%a Yetter\dash Drinfeld module $Y$, as the \emph{Heisenberg quadratic double}.
%For future reference, we state the presentation of $\mindouble
%A_Y$ in 
%
%\begin{lemma}
%\label{lem:Heisenberg}
%In the above notation, the algebra $\mindouble A_Y$ is 
%generated by $Y$, $W$, $Y^*$ subject to the relations 
%$$
%  wxw^{-1}=w(x),\quad
%  f  v - v f = 
%  \langle v,f\rangle |v|, \quad
%  \ker(\id+\Psi_Y), \quad \ker(\id+\Psi_Y^*),
%$$
%where $w\in W$, $x\in (Y\text{ or }Y^*)$, $f\in
%Y^*$, $v\in Y$.
%\qed
%\end{lemma}
%
%
%\subsection{Universality of Heisenberg quadratic doubles}
%
%Here comes 
The crucial property of Heisenberg quadratic doubles is given in  

\begin{theorem}
\label{thm:perfect}
For any finite\dash dimensional $W$\dash module $V$ and any 
two quadratic doubles  $A$, $B$ in $\qdcat(V)$, there exists a
finite\dash dimensional Yetter\dash Drinfeld module $Y$ over $W$ and a
morphism   $A\star B\to A_Y$ in $\qdcat_W$.  
\end{theorem}
%\begin{remark}
%Given a quadratic double $A$, one would ideally like to embed $A$ in a
%Heisenberg quadratic double $A_Y$. In a number of cases, this can be
%achieved by representing 
%$A$ as $B\star C$ for some quadratic doubles $B$, $C$, applying the
%Theorem and using certain minimality properties of $A$ which guarantee
%that the morphism is injective. (One can show that $A$ being a minimal
%quadratic double is a necessary condition for injectivity.) 
%In the present paper, we use this method
%for rational and braided Cherednik algebras.   
%\end{remark}

\begin{proof}%[Proof of Theorem~\ref{thm:perfect}]
Let $A=A_\beta(R^-,R^+)$ and $B=A_\gamma(S^-,S^+)$, where $\beta$,
$\gamma$ are $W$\dash equivariant maps from $V^*\tensor V$ to $\CC W$ given
by 
$
   \beta(f\tensor v) = \sum_{w\in W} \langle L_w(v), f\rangle w
$, 
%\qquad
$      \gamma(f\tensor v) = \sum_{w\in W} \langle M_w(v), f\rangle w
$
%,
%\qquad 
%f\in V^*, \ v\in V,
%$$ 
with $L_w,M_w\in \End(V)$. The 
%following 
finite subset
$%$
   E = \{g\in W \mid L_g\ne 0\ \text{or}\ M_g\ne 0\}
$
%Clearly,  
of $W$ is conjugation\dash invariant by 
the $W$\dash equivariance of $\beta$, $\gamma$.
% implies that
%$E$ is a conjugation\dash invariant subset of $W$. 
Denote by $\CC E$
the linear span of $E$ in $\CC W$. 
We introduce the space $Y$ equipped  with 
$W$\dash action and $W$\dash grading  by
$$
     Y = \CC E \tensor V,
\qquad
     w(g\tensor v)=wgw^{-1}\tensor w(v),
\qquad
     |g\tensor v|=g
$$
for all $g\in E$, $v\in V$, $w\in W$.
It is easy to see that $Y$ is a Yetter\dash Drinfeld module for $W$. 
The dual Yetter\dash Drinfeld 
module $Y^*$ can also be described explicitly:
$
  Y^* = \CC E^{-1} \tensor V^*$, 
%\qquad
$      w(h\tensor f)=whw^{-1}\tensor w(f)$,
%\qquad
$     |h\tensor f|=h 
$
for all $h\in E^{-1}=\{g^{-1}\mid g\in E\}$, $f\in V^*$ and $w\in
W$. One checks that $
 \langle g\tensor v, h\tensor f\rangle = \delta_{g^{-1},h}  \langle
   v,f\rangle$
is a pairing between $Y$ and $Y^*$ that indeed makes $Y^*$ the YD
module dual to $Y$.
%Let the 
%$W$\dash module homomorphisms 
The maps $\mu\colon V\to Y$, $\nu \colon
V^* \to Y^*$ given by 
$$
     \mu(v) = \sum_{w\in W} w\tensor M_w(v), 
\qquad
     \nu(f) = \sum_{w\in W} w^{-1} \tensor L_w^*(f) 
$$
%The maps $\mu$, $\nu$ 
are $W$\dash module homomorphisms because
$\gamma$, $\beta$ are $W$\dash equivariant. 

%In the rest of the proof we 
It remains to show that $\mu$, $\nu$ extend to a
morphism between the quadratic doubles $A\star B=A_{\beta\gamma}(S^-,R^+)$
and $A_Y$. 
As in the proof of Proposition~\ref{prop:diamond}, it is enough to
show that 
$$
 \beta\gamma=\beta_Y\circ(\nu\tensor\mu), 
\qquad
 (\mu\tensor \mu)S^- \subset \ker(\id_{Y\tensor Y}+\Psi_Y),
\qquad
 (\nu\tensor\nu)R^+ \subset \ker(\id_{Y^*\tensor Y^*}+\Psi_Y^*).
$$
Since  $\beta_Y(h\tensor f \tensor g \tensor v)=\delta_{g,h^{-1}} \langle
v,f\rangle g$ where $g,h^{-1}\in E$, $f\in V^*$, $v\in V$, 
one has $\beta_Y(\nu(f)\tensor\mu(v))=\sum_{w\in E} \langle M_w(v),
L_w^*(f)\rangle w = \sum_w \langle (L_w M_w)(v),
f\rangle w =(\beta\gamma)(f\tensor v)$ as required. 
%Next, we have to check that $(\id+\Psi_Y)(\mu\tensor \mu)S^-=0$. 
%Observe that $(\mu\tensor \mu)(u\tensor v)=\sum_{w,x\in W} w\tensor
%M_w(u)\tensor x \tensor M_x(v)$, so that
%\begin{align*}
%  (\id+\Psi_Y)(\mu\tensor \mu)(u\tensor v) & = 
%  \sum_{w,x\in W} w\tensor M_w(u)\tensor x \tensor M_x(v)+
%       wxw^{-1}\tensor w(M_x(v)) \tensor w \tensor M_w(u)   
%\\
%&= \sum_{w\in W} 
%\mu(u) \tensor w \tensor M_w(v) + w(\mu(v)) \tensor w \tensor M_w(u). 
%\end{align*}
%Given that $w(\mu(v))=\mu(w(v))$, the latter expression is obtained by
%cyclic permutation of tensor factors in 
%$$
%   \sum_{w\in W} w\tensor (\id_V \tensor \mu)(M_w\tensor
%   \id_V)(v\tensor u + u\tensor w(v)) = 
%   \sum_{w\in W} w\tensor (\id_V \tensor \mu) T^-_{w,\gamma}(u\tensor v),
%$$
%where the maps $T^-_{w,\gamma}$ are as introduced in Lemma~\ref{lem:T}. 
%It follows that the kernel of $(\id+\Psi_Y)(\mu\tensor \mu)$
%in $V\tensor V$ contains $\cap_{w\in W}\ker T^-_{w,\gamma}$, which by
%Lemma~\ref{lem:T} contains $S^-$.
%
%In the same way it is shown that 
%$(\id_{Y^*\tensor Y^*}+\Psi_Y^*)(\nu\tensor \nu)(R^+)=0$. 
The remaining two equalities are established by applying
Lemma~\ref{lem:T} (similarly to the proof of
\cite[Theorem 6.9]{BB}). 
The Theorem is proved.
\end{proof}
%\begin{remark}
%\label{rem:submodule}
%Instead of the Yetter\dash Drinfeld module $Y$ constructed in the
%proof of the Theorem, its YD submodule spanned by $w\tensor M_w(v)$
%and $w\tensor L_w(v)$, $w\in W$, $v\in V$, can be used. The argument
%remains the same.
%\end{remark}

\startcmm

\subsection{Rational Cherednik algebras and Heisenberg
  quadratic doubles}
\label{subsect:H0}

We will now apply Theorem~\ref{thm:perfect} to rational Cherednik
algebras. Let $W$ be a finite subgroup of $\GL(V)$ and $S$ be the set of
complex reflections in $W$. Recall the rational Cherednik algebra
$H_c(W)$ (Example~\ref{ex:rca}), and consider its 
``degenerate'' version
$$
H_{0,c}(W):=A_{\beta'}(\wedge^2 V, \wedge^2 V^*)\cong S(V)\tensor
\CC W \tensor S(V^*), 
\qquad
\beta'(\xi\tensor v)=\sum_{s\in S} c_s \langle v,\alpha_s^\vee\rangle \langle
\alpha_s, \xi\rangle,
$$
where the commutator $\beta'$ lacks the $\langle v,\xi\rangle$ term. 
The maps $L'_s\in \End(V)$, $s\in S$ given by  
$L'_s(v)=c_s \langle v,\alpha_s^\vee \rangle \alpha_s$ are idempotents
up to a scalar factor: one has $L_s'L_s' = c_s \langle
\alpha_s,\alpha_s^\vee\rangle L_s'$. Therefore, 
putting $c_0=\{\langle
\alpha_s,\alpha_s^\vee\rangle^{-1}: s\in S\}$ we get
$$
         H_{0,c}(W) = H_{0,c}(W) \star H_{0,c_0}(W).  
$$  
Theorem~\ref{thm:perfect} and Remark~\ref{rem:submodule} now suggest
to consider the vector space
$$
Y_{\mathit{refl}}=\bigoplus\limits_{s\in S} \field\cdot [s]
$$
with the basis $\{[s] \mid s\in S\}$ labelled by complex reflections in $W$.
(This is a YD submodule of $\CC S \tensor V$ where $[s]$ stands for
$s\tensor \alpha_s$.) 
The space $Y_{\mathit{refl}}$ is a Yetter\dash Drinfeld module
for $W$, via the $W$\dash action  
$$
     \sigma([s])=\chi(\sigma,s)\cdot [\sigma s\sigma^{-1}] \ ,
$$ 
where $\chi\colon W\times S\to \field^\times$ is a function such that 
$\sigma(\alpha_s)=\chi(\sigma,s)\alpha_{\sigma s \sigma^{-1}}$. 
The grading on $Y_{\mathit{refl}}$ is by $\big|[s]\big|=s$. 
The dual basis of $Y_{\mathit{refl}}^*$ will be written as 
$\{[s]^* \mid s\in S\}$. (In fact, $Y_{\mathit{refl}}^*$ is
a YD submodule of $\CC S \tensor V^*$, where $[s]^*$ is 
$s^{-1}\tensor \langle \alpha_s,\alpha_s^\vee\rangle^{-1} \alpha_s^\vee$.) From the proof of
Theorem~\ref{thm:perfect} we deduce 
\begin{theorem}
\label{thm:emb0}
The $W$\dash equivariant maps 
$\mu_c\colon V \to  Y_{\mathit{refl}}$ and $\nu \colon V^* \to
Y_{\mathit{refl}}^*$ given by 
$$
\mu_c(v)=\sum_{s\in S}c_s \langle v, \alpha_s^\vee \rangle  [s],
\qquad
\nu(\xi)=\sum_{s\in S} \langle \alpha_s, \xi \rangle  [s]^*
$$
extend to an algebra map
$$
H_{0,c}(W) \to
 \mindouble A_{Y_{\mathit{refl}}} = 
 S(Y_{\mathit{refl}},\Psi_{Y_{\mathit{refl}}}) \tensor 
 \field W \tensor 
 S(Y_{\mathit{refl}}^*,\Psi^*_{Y_{\mathit{refl}}}).
\qed
$$
\end{theorem}

Theorem~\ref{thm:emb0} provides a realisation only of degenerate
Cherednik algebras $H_{0,c}(W)$ but not of $H_c(W)$, 
and the map $H_{0,c}(W)\to \mindouble
A_{Y_{\mathit{refl}}}$ is not injective. In fact, the algebras 
$S(Y_{\mathit{refl}},\Psi_{Y_{\mathit{refl}}})$ and $S(Y^*_{\mathit{refl}},\Psi_{Y_{\mathit{refl}}}^*)$ may even be finite\dash
dimensional. A good example is $W=\Symm_n$ where these algebras are
isomorphic to the Fomin\dash Kirillov quadratic algebra $\mathcal
E_n$, finite\dash dimensional for $n\le 5$ \cite{FK}.

There are two ways to correct this situation. First, one can construct
a bigger Yetter\dash Drinfeld module $V\oplus Y_{\mathit{refl}}$
by adding to $Y_{\mathit{refl}}$ a copy of $W$\dash module
$V$ with trivial $W$\dash grading: $V=V_1$. The maps $\mu_c$, $\nu$
from Theorem~\ref{thm:emb0} are replaced by  
${\mu_c}_{\mathrm{new}}=\id_V \oplus \mu_c\colon V \to V\oplus
Y_{\mathit{refl}}$  and 
${\nu}_{\mathrm{new}}=\id_{V^*} \oplus \nu \colon V^* \to 
V^*\oplus Y^*_{\mathit{refl}}$. By Theorem~\ref{thm:perfect}, 
they give rise to an embedding of
$H_c(W)$ in the Heisenberg quadratic double $\mindouble A_{V\oplus
  Y_{\mathit{refl}}}$.

Nevertheless, our preferred second way to deal with algebras $H_c(W)$  
is to ``modify'' the Heisenberg quadratic double of
$Y_{\mathit{refl}}$ as follows. We will use the fact that
$$
   H_c(W) = H_{0,c}(W) \diamond H_0(W),
$$
where, of course, $H_0(W)=\mathcal D(V)\lcprod \CC W$, the algebra
$\mathcal D(V)$ being the Weyl algebra (Example~\ref{ex:weyl}). For a Yetter\dash
Drinfeld module $Y$ over $W$, let
$$
\mathcal D_t (Y)\lcprod\CC W \cong S(Y) \tensor \CC W \tensor 
S(Y^*), \qquad t\in \CC
$$
be the quadratic double over $W$ with defining relation
$fy-yf=t\langle y,f\rangle$ for $y\in Y$, $f\in Y^*$. It is isomorphic
to the Weyl algebra $\mathcal D(Y)$ of $Y$. The following Lemma
is straightforward:

\begin{lemma}
\label{lem:value}
The above maps $\mu_c\colon V \to Y_{\mathit{refl}}$ and $\nu \colon V^* \to
Y^*_{\mathit{refl}}$ extend to a 
morphism $H_0(W) \to \mathcal D_t(Y_{\mathit{refl}})\lcprod \CC W$ of quadratic
doubles, if and 
only if $t\sum_{s\in S}c_s \langle v,\alpha_s^\vee\rangle \langle
\alpha_s,\xi \rangle =\langle v,\xi\rangle$ for all $v\in V^*$,
$\xi\in V^*$. 
\qed
\end{lemma} 

The embedding theorem is easiest to state for rational Cherednik
algebras of irreducible groups $W$ (i.e., such that $V$ is an
irreducible $W$\dash module):

\begin{theorem}
\label{thm:emb}
In the above notation, let $W$ be an irreducible subgroup of $\GL(V)$. 
For any tuple $c$ of parameters outside a certain hyperplane, 
there exists $t\in \CC$ 
such that the maps $\mu_c\colon V \to Y_{\mathit{refl}}$, $\nu
\colon V^* \to Y^*_{\mathit{refl}}$ extend to an embedding 
$$
    H_c(W) \hookrightarrow A^t_{Y_{\mathit{refl}}} :=
    A_{Y_{\mathit{refl}}} \diamond ( D_t (Y_{\mathit{refl}})\lcprod \CC W)
$$
of quadratic doubles. 
\end{theorem}
\begin{proof}
Let $C\subset S$ be a conjugacy class.
The sum $\sum_{s\in C}\langle v,\alpha_s^\vee\rangle \langle
\alpha_s,\xi \rangle$ is a nonzero $W$\dash invariant pairing between
$v\in V$ and
$\xi \in V^*$, hence it is proportional to $\langle v,\xi\rangle$ if $W$ is
irreducible. Therefore, the sum in Lemma~\ref{lem:value} represents a
pairing between $V$ and $V^*$ proportional to $\langle \cdot, \cdot
\rangle$. The coefficient of proportionality is non\dash zero unless $c$ belongs to a certain exceptional
hyperplane in the space of parameters, and is made equal to $1$ by
 choosing $t$ appropriately. 
Then by Lemma~\ref{lem:value} the maps $\mu_c$ and $\nu$ extend to a
morphism $H_0(W) \to \mathcal D_t (Y_{\mathit{refl}})\lcprod \CC W$, and thus by
Theorem~\ref{thm:emb0} and Proposition~\ref{prop:diamond} they extend
to a morphism from $H_c(W)$ to $A^t_{Y_{\mathit{refl}}}$ in the
category of quadratic doubles. 

It remains to argue that every such morphism with $H_c(W)$ as the
source is injective. This is guaranteed by the fact that $H_c(W)$ is a
minimal braided double --- a condition, stronger than the minimality
of quadratic doubles, established in \cite[Proposition 7.11]{BB}.
\end{proof}
\begin{remark}
Both the irreducibility of $W$ and the condition on $c$ are not
essential and may be eliminated by replacing the single parameter $t$ 
with conjugation\dash
invariant parameters
$t_s$, $s\in S$ ($t_s=0$ is allowed). 
The commutation relation in $A^t_{Y_{\mathit{refl}}}$
becomes $[r]^*[s]-[s][r]^*=t_s\delta_{r,s}$. 
\end{remark}
Of course, information provided by Theorem~\ref{thm:emb} is incomplete
without a description of the algebra $A^t_Y$ that contains the
rational Cherednik algebra of $W$. 
For a braiding $\Psi$ on a space $Y$, 
introduce the following quadratic algebra, which is an extension of
both the symmetric algebra $S(Y)$ and the algebra $S(Y,\Psi)$:
$$
S'(Y,\Psi) = T(Y)/ \, \lgen \,\ker(\id+\Psi) \cap \wedge^2 Y\,
\rgen.
$$
\begin{proposition}
For a finite\dash dimensional Yetter\dash Drinfeld module $Y$ over
$W$ and for $t\in\field$, 
the quadratic double $A^t_Y=A_Y \diamond (\mathcal D_t(Y)\lcprod \CC W)$ is of the form 
$$ 
    A^t_Y = S'(Y,\Psi_Y) \tensor \field W \tensor S'(Y^*,\Psi_Y^*),
$$
and is generated by $Y$, $W$, $Y^*$ subject to the following relations which
depend on $t$:
$$
wxw^{-1}=w(x),\quad
  f  y - y f = 
  \langle y,f\rangle(|y|+t\cdot 1), \quad
  \ker(\id+\Psi_Y) \cap \wedge^2 Y, 
  \quad 
  \ker(\id+\Psi_Y^*) \cap \wedge^2 Y^* 
$$
with $w\in W$, $x\in (Y\text{ or }Y^*)$, $f\in
Y^*$, $y\in Y$.
\qed
\end{proposition}

\begin{remark}
\label{rem:BK}
To better understand the algebras $\mindouble A_{Y_{\mathit{refl}}}$ and
$A^t_{Y_{\mathit{refl}}}$, it is instructive to analyse them in the case when 
$W$ is a  symmetric or, more generally, Coxeter group. 
In this case, the algebra
$S'(Y_{\mathit{refl}},\Psi_{Y_{\mathit{refl}}})$ is studied by
Kirillov and the first author in \cite{BK}. 
When $W=\Symm_n$, the algebra  $S({Y_{\mathit{refl}}},\Psi_{Y_{\mathit{refl}}})$ is the Fomin\dash Kirillov quadratic
algebra $\mathcal E_n$ related to the Schubert calculus \cite{FK}, 
and the algebra 
$S'(Y_{\mathit{refl}},\Psi_{Y_{\mathit{refl}}})$ 
coincides with the triangular enveloping algebra
$U(\mathrm{tr}_n)$ introduced by
Bar\-tholdi-En\-ri\-quez-Etin\-gof-Rains 
in \cite{BEER}.
See \cite[Examples 7.23 and 7.27]{BB} for more details.
\end{remark}
\endcmm

\subsection{Generalised Dunkl operators}
\label{subsect:revisited}

We now look at the Dunkl operators
in the context of quadratic doubles and propose their
generalisation.

Let $V$ be a finite\dash dimensional module over a group $W$ and $A=U^-\tensor \field W \tensor U^+$ be a quadratic double in $\qdcat(V)$. 
%By definition of a quadratic double, 
The subalgebra
$\CC W\tensor U^+$ of $A$ 
%is isomorphic to the semidirect product $\CC W \rcprod U^+$ and 
has one\dash dimensional trivial
representation $\mathit{triv}$, where $w\in W$ acts by $1$ and $V$ acts by
$0$. One has the induced representation of $A$:
$$
    \mathrm{Ind}_{\CC W\tensor U^+}^A (\mathit{triv})  \quad 
    \cong U^-\ \text{as vector spaces.}
$$
Denote the action of $\xi\in V^*$ on $U^-$ via this representation by
$\partial_\xi$. The operators $\partial_\xi \colon U^- \to U^-$
are of degree $-1$ with respect to the grading in $U^-$.

In Heisenberg quadratic doubles, the operators $\partial_\xi$ are
called \emph{braided derivatives}. They can be computed 
%in an efficient way 
as follows. 
Suppose that $Y=\oplus_{w\in W}Y_w$ is a Yetter\dash Drinfeld
module for $W$. Let $\{y_w^i \in Y_w \mid w\in W, \ i=1,\dots,\dim
Y_w\}$ be a basis of $Y$ compatible with the $W$\dash grading, and let 
$f^i_w\in (Y^*)_{w^{-1}}$ form the dual basis of $Y^*$. One shows that 
the braided derivatives 
$$
\partial^i_w = \partial_{f^i_w}\colon 
S(Y,\Psi_Y) \to S(Y,\Psi_Y)
$$ 
satisfy, and are determined by, the following
properties:
\begin{itemize}
\item[$(i)$] $\partial^i_w y^j_\sigma =
  \delta_{i,j}\delta_{w,\sigma}$;
\item[$(ii)$] ($w$\dash twisted Leibniz rule) 
$\partial^i_w (ab) = (\partial^i_w a) w(b) + a(\partial^i_w b)$ for all
  $a,b\in S(Y,\Psi_Y) $. 
\end{itemize}
Let us now consider a quadratic double 
$A=U^-\tensor \field W \tensor U^+$ in $\mathcal \qdcat(V)$ which is not
Heisenberg. For emphasis, we will now write $\nabla_\xi$ instead of
$\partial_\xi$ in $A$ and call $\nabla_\xi$ \emph{generalised Dunkl
  operators}. 

Suppose that $A$ embeds in a Heisenberg quadratic double $A_Y$ for
some  Yetter\dash Drinfeld module $Y$, and let 
$V\xrightarrow{\mu} Y$, $V^* \xrightarrow{\nu} Y^*$ be a pair of maps
which give rise to such embedding. Put
$$
     \alpha^i_w := \nu^*(y^i_w) \in V,
\qquad
    {\alpha^\vee}^i_w := \mu^*(f^i_w)\in V^*,
\qquad
    i=1,\dots,\dim Y_w.
$$
The vectors $\alpha^i_w$, respectively ${\alpha^\vee}^i_w$, are an analogue
of roots, respectively coroots, of the group $W$. Initial
data for this generalised root system is a $W$\dash module $V$
equipped with $W$\dash homomorphism $\beta$ from $V^*\tensor V$ to the
adjoint representation of $W$. Strictly speaking, the construction
depends on a particular ``quantisation'' $Y$ of $(V,\beta)$, but one hopes
that the ``root system'' has a reasonable uniqueness property; we do
not consider this question here.
It is not difficult to see that the generalised Dunkl operators are
expressed as
$$
   \nabla_\xi = \sum_{w\in W}\sum_{i=1}^{\dim Y_w}
    \langle \alpha^i_w,\xi\rangle \bar\partial^i_w,
$$
with $\bar\partial^i_w\colon U^-\to U^-$ defined by its two properties,
$$
  \bar\partial^i_w(v)=\langle v,{\alpha^\vee}^i_w \rangle, \quad v\in
  V;
  \qquad
  \bar\partial^i_w\quad \text{satisfies the $w$-twisted Leibniz rule.}
$$
%
%The Dunkl operators for reflection groups are a model case of
%this. Explicitly, let $V$ be the reflection representation of
%a complex reflection group $W$. Constructions in and before
%Theorem~\ref{thm:emb} imply the following analysis of Dunkl's
%operators $\nabla_\xi$. For $w=1\in W$ choose
%$\alpha_1^i$, $i=1,\dots,\dim V$ to 
%be any basis of $V^*$, and let $\{{\alpha^\vee}_1^i\}$ be the dual basis
%    of $V^*$. 
%For $w=s\in S$ put $\alpha^1_w=\alpha_s$ to
%be the usual root of $s$, and let ${\alpha^\vee}^1_s =
%c_s \alpha^\vee_s$.
%This explains the logic behind the formula for
%Dunkl operators, quoted in Example~\ref{ex:rca}.
%
%\begin{remark}
%In the case of a complex reflection group, the
%operators $\bar\partial^i_w$ are written as directional derivatives
%(for $w=1$) or as Demazure divided difference operators
%$\frac{1}{\alpha_s}(1-s)$ for $w=s\in S$. The classical Dunkl operators 
%are then treated as elements of the ring $\mathcal
%D(V_{\mathit{reg}})\lcprod \CC W$, where $V_{\mathit{reg}}$ is the variety
%$\{x\in V \mid s(x)\ne x \ \forall s\in
%S\}$, and $\mathcal D(\cdot)$ stands for the algebra of polynomial
%differential operators. %
%
%Although this property is not expected to have an analogue for
%$\nabla_\xi$ in an arbitrary quadratic double, a ``noncommutative'' 
%version of it is present in formula (1) for the braided Dunkl
%operators $\rd\nabla_i$, given in the Introduction. 
%\end{remark}
%
\subsection{Braided reduction of quadratic doubles}
\label{subsect:Braided reduction}

To conclude this Section, we describe a method of obtaining 
a wider class of algebras with triangular decomposition over $\CC
W$ as subalgebras in quadratic doubles. Recall that if $Y$ is a
finite\dash dimensional Yetter\dash Drinfeld module over $W$, then $\Psi_Y$ is
the braiding $y\tensor z\mapsto |y|(z)\tensor y$ on $Y$, 
and $\tau$ is the trivial
braiding. The spaces $Y^*$, $Y^*\tensor Y^*$ etc.\ are also YD-modules, 
and the $W$\dash grading on $Y^*\tensor Y^*$ is by $|f\tensor g|=|f||g|$. 

\begin{proposition}
\label{prop:twist}
Assume that a $W$\dash module $V$ has $W$\dash grading which makes $V$
a Yetter\dash Drinfeld module. Let $A=A_\beta(R^-,R^+)$ be a quadratic
double in $\qdcat_\beta(V)$, such that $\tau(R^+)$ is a $W$-graded subspace of
$V^*\tensor V^*$. Let
$$
\rd A = \frac{T(V\oplus V^*)\lcprod \CC W}{\lgen \, R^-, \ \rd R^+,\ 
\{ \rd{[\theta , v ]} - |\theta|\beta(\theta\tensor v) 
: \theta\in V^*, \ v\in V\} \,\rgen}\ , 
$$
where $\rd{[f,v]} = f\tensor v - |f|(v)\tensor f$ is the braided
commutator between $V^*$ and $V$, and 
$\rd R^+= (\Psi_{V^*}\circ \tau)R^+$. Then there exists an algebra isomorphism
$t\colon \rd A \to A$, 
given on generators 
of $\rd A$ by 
$$
   t|_V=\id_V, \qquad t|_W=\id_W, \qquad t(\theta)=|\theta|\cdot
   \theta, \quad \theta\in V^*.
$$
\end{proposition}
\begin{proof}
First, we have to show that $t$ maps relations in $\rd A$ to relations
in $A$.
For $\theta\in V^*$ we have $t(w\theta w^{-1}-w(\theta))=w\cdot |\theta|\theta \cdot
w^{-1}-|w(\theta)|w(\theta)$. This is a relation in $A$, as
$|w(\theta)|=w|\theta|w^{-1}$ by the Yetter\dash Drinfeld condition on
$V^*$. Furthermore, $t(\rd{[\theta , v ]}) = |\theta|\theta v
-|\theta|(v)\cdot |\theta|\theta = |\theta|(\theta v - v\theta)$ which
in $A$ coincides with $|\theta|\beta(\theta,v)$. 

Now observe that $\Psi_{V^*}(\tau(\theta\tensor
\kappa))=|\kappa|(\theta)\tensor \kappa$ for $\theta\tensor\kappa\in
V^*\tensor V^*$. This is mapped by $t$ to $|\kappa||\theta||\kappa|^{-1}\cdot
|\kappa|(\theta)\cdot |\kappa|\kappa = |\kappa||\theta| \cdot
\theta\kappa$. Hence if $a\in R^+$ is such that $\tau(a)$ is
$W$\dash homogeneous (such $a$ span $R^+$), then
$t(\Psi_{V^*}\tau(a))=|\tau(a)|a$. Thus, $t(\rd R^+)\subset \CC
W\tensor R^+$. It remains to note that the relations $wvw^{-1}-w(v)$
and $R^-$ in $\rd A$ are mapped by $t$ to exactly the same relations
in $A$. We conclude that $t$ is a map of algebras. 

In the same fashion one shows that the map $t^{-1}$, given on
generators of $A$ by $t^{-1}|_V=\id_V$, $t^{-1}|_W=\id_W$,
$t^{-1}(f)=|f|^{-1} \cdot f$, $f\in V^*$, is an algebra homomorphism
from $A$ to $\rd A$. As $tt^{-1}$ and $t^{-1}t$ are identity on
generators, $t^{-1}$ is the inverse of $t$. 
\end{proof}
\begin{remark}
It is easy to deduce from the Proposition that 
the algebra $\rd A$, given by its presentation, has
triangular decomposition $\rd A=T(V)/\lgen R^-\rgen \tensor \CC W
\tensor T(V^*)/\lgen \rd R^+ \rgen$. 
\end{remark}
\begin{definition}
In the above notation, assume that for all $\theta\in V^*$, $v\in V$
the braided commutator $\rd{[\theta,v]}$ in $\rd A$  lies in 
$\CC W'$ for some 
subgroup $W'$ of $W$. The algebra $A=\rd A$ has
subalgebra $\rd A' \cong  T(V)/\lgen R^-\rgen \tensor \CC W'
\tensor T(V^*)/\lgen \rd R^+ \rgen$. We call $\rd A'$ a \emph{braided
  reduction} of $A$. 
\end{definition}
\begin{example}[Braided Weyl algebra]
\label{ex:braidedweyl}
To show how the braided reduction works, we consider the ``extreme''
example which is $A=\mindouble A_V$, the Heisenberg quadratic double of a
Yetter\dash Drinfeld module $V$ over $W$. 
We compute the braided commutator in $\rd A$ of $\theta\in V^*$ and
$v\in V$:
$$
    \rd{[\theta,v]}=|\theta|\beta_V(\theta,v) = 
  |\theta|\cdot \langle v,\theta\rangle\cdot |v| = \langle
    v,\theta\rangle\cdot 1,      
$$ 
as $\langle v, \theta\rangle\ne 0$ for 
$W$\dash homogeneous $\theta$, $v$ only if
$|\theta|=|v|^{-1}$. We thus have a braided reduction   
$\rd A'\cong  T(V)/\lgen R^-\rgen \tensor \CC \cdot 1
\tensor T(V^*)/\lgen \rd R^+ \rgen$ of $A$. 
Furthermore, using 
$\Psi_{V^*}=\tau\circ \Psi_V^*\circ \tau$
%Remark~\ref{rem:adjbraiding} 
we find 
$\rd R^+= (\Psi_{V^*}\tau)\ker(\id+\tau \Psi_{V^*}\tau) = 
\ker(\id+\Psi_{V^*})$. Hence
$$
 \rd A_V := \rd A' \cong S(V,\Psi_V) \tensor S(V^*,\Psi_{V^*})
$$
with defining commutation relation 
$\theta v-|\theta|(v)\theta = \langle v,\theta\rangle\cdot 1$ between $\theta\in V^*$ and
$v\in V$. 
We have 
$\mindouble A_V \cong \rd A_V \lcprod \CC W$ as algebras.
The algebra $\rd A_V$  is a particular case of 
Majid's braided Weyl algebra \cite{Mcalculus}, hence
\begin{definition}
\label{def:braidedweyl}
$\rd A_V$ is called the \emph{braided Weyl algebra} of
the Yetter\dash Drinfeld module $V$.
\end{definition} 
\end{example}

\section{$\qq$-Cherednik algebras}
\label{sect:q_cher}

In this Section we introduce $\qq$\dash Cherednik algebras. 
They are quadratic doubles, which allows us to use the
methods of
%developed in 
Section~\ref{sect:qd}. 
%and \ref{sect:yd}. 
On the other
hand, results about $\qq$\dash Cherednik algebras will be translated
to  braided Cherednik algebras, obtained from $\qq$\dash Cherednik
algebras by braided reduction.

%\subsection{The braiding $\tau_\qsub$}
\subsection{The $\qq$-polynomial algebra}

Recall that a rational Cherednik algebra of a finite linear group $W\subset
\GL(V)$ is a deformation of the semidirect product 
$\mathcal D(V)\lcprod \CC W$, where  $\mathcal D(V)\cong
S(V)\tensor S(V^*)$  
is the Weyl algebra of the space $V$. 
Our 
%general  
aim is to replace the polynomial algebra $S(V)$ with
its $q$\dash analogue.
% $S_\qsub(V)$. 
%noncommutative deformation, the 
%Koszul algebra $S(V,\Psi)=T(V)/\lgen \id_{V\tensor V} + \ker \Psi \rgen$
%associated to a braiding 
%$\Psi\in \End (V\tensor V)$ of Hecke type (see Remark~\ref{rem:Hecke}).
%In the present paper we only consider
%braidings $\Psi$ which meet the conditions in
Throughout, 
%$V$ will be a vector space of dimension $n\ge 1$ and 
%the symbol $\qq$ will stand for an $n\times n$ matrix of deformation  
%parameters $q_{ij}$ 
%such that $q_{ii}=1$, $q_{ij}q_{ji}=1$ for all $i,j=1,\dots,n$. 
\startcmm
\begin{definition}
\label{def:type}
(1) 
A braiding $\Psi$ on a space $V$ is \emph{of group type},
if $V$ has a Yetter\dash Drinfeld module structure over some group
$\Gamma$ which induces the braiding $\Psi$.

(2) $\Psi$ is a \emph{flat deformation} of the trivial braiding $\tau$,
if $S(V,\Psi)\cong S(V)$ as graded vector spaces. 
\end{definition}

The group type condition \cite{T}
allows us to use results from Sections
\ref{sect:qd} and \ref{sect:yd}; if it is relaxed, more general
methods for Hopf algebras such as in \cite{BB} are to be applied.  
Condition (2) restricts the
class of algebras we are constructing and working with; relaxing it
to, say, certain types growth may lead to new classes of quadratic
doubles, to be studied elsewhere. 

A family of flat deformations of $\tau$ of group type is constructed in
the following 
\begin{lemma}
\label{lem:tau}
Let $x_1,\dots,x_n$ be a basis of $V$ and 
let $\qq=(q_{ij})_{i,j=1}^n$ be a complex matrix such that
$q_{ii}=1$, $q_{ij}q_{ji}=1$ for all $i,j=1,\dots,n$. The map
$$
\tau_\qsub\colon V\tensor V \to V \tensor V, \qquad
\tau_\qsub(x_i \tensor x_j) = q_{ij} x_j \tensor x_i 
$$
is a braiding on $V$ of group type which is a flat deformation of
$\tau$, and  $\tau_\qsub^2=\id$.
\end{lemma}
\begin{proof}
The braid equation for $\tau_\qsub$ applied to $x_i
\tensor x_j \tensor x_k$ reduces to a tautological equality 
$q_{ij}q_{ik}q_{jk}=q_{jk}q_{ik}q_{ij}$.
That $\tau_\qsub^2=\id$, follows from $q_{ij}q_{ji}=1$. 
To show that $\tau_\qsub$ is of group type, consider the 
elements  $\gamma_1,\dots,\gamma_n\in  \GL(V)$ and 
the group $\Gamma_\qsub \subset\GL(V) $ defined by
$$
    \gamma_i(x_j) = q_{ij} x_j, \quad i,j=1,\dots,n;
\qquad
    \Gamma_\qsub=\langle \gamma_1, \dots,\gamma_n\rangle\, . 
$$
(Clearly, $\Gamma_\qsub$ is finite if and only if all $q_{ij}$ are
    roots of unity.)
The space $V$ is $\Gamma_\qsub$\dash graded via $|x_i|=\gamma_i$. 
This grading and the natural action of $\Gamma_\qsub$ on $V$ make $V$ a
Yetter\dash Drinfeld module over $\Gamma_\qsub$. Indeed, as $\Gamma_\qsub$
is Abelian (a subgroup of a maximal torus of $\GL(V)$), the
Yetter\dash Drinfeld condition reads 
$|\gamma_i(x_j)|=|x_j|$, which is true. It is easy to see
that the induced braiding is $\tau_\qsub$. Finally, the space of
quadratic relations in the
algebra $S(V,\tau_\qsub)$ is
$$
   \wedge^2_\qsub V = \mathrm{span}(x_i\tensor x_j - q_{ij} x_j \tensor
   x_i \mid 1\le i,j\le n) = \ker(\id+\tau_\qsub)
   \quad \subset V\tensor V,
$$
thus $S(V,\tau_\qsub)$ is the algebra of ``$q$\dash polynomials'' which
is a flat deformation of $S(V)$.
\end{proof}
Note that $\tau_\qsub$ is an involution, which is the simplest case of a
Hecke type braiding. 
It turns out that the $\tau_\qsub$ exhaust all suitable braidings:
\begin{proposition}
All  flat deformations of $\tau$ of Hecke type and of
group type are $\tau_\qsub$ given by Lemma~\ref{lem:tau}.   
\end{proposition}
\begin{proof}
Suppose that $V$ is a Yetter\dash Drinfeld module over a group
$\Gamma$, and put $\Psi=\Psi_V$. The Hecke condition for $\Psi$ reads
$\Psi^2+(1-\lambda) \Psi=\lambda \id$ with $\lambda\ne 0$. 
This is equivalent to the equation
$$
   aba^{-1}(x)\tensor a(y) + (1-\lambda)a(y)\tensor x = \lambda
   x\tensor y
\qquad\text{for all}\ a,b\in \Gamma,\ x\in V_a,\ y\in V_b\,.
$$ 
The terms of this equation lie respectively in
$V_{abab^{-1}a^{-1}}\tensor V_{aba^{-1}}$, $V_{aba^{-1}}\tensor V_a$
and $V_a\tensor V_b$. These subspaces of $V\tensor V$ form a direct
sum unless $ab=ba$. We conclude that $\{a\in \Gamma \mid V_a\ne 0\}$ commute. 

Now note that, assuming $V_a\ne 0$ and $V_b\ne 0$, 
the above equation may only hold when $b=a$ or $\lambda=1$.
If $\lambda=1$, we have $b(x)\tensor a(y) = x\tensor y$, which means
that the $\Gamma$\dash degree of a non\dash zero  $\Gamma$\dash
homogeneous element of $V$ 
acts on each $\Gamma$\dash homogeneous subspace of $V$ by
a scalar. 

Choose a basis $x_1,\dots,x_n$ of $V$ compatible with the
$\Gamma$\dash grading, and denote $\gamma_i=|x_i|$; then one has 
$\gamma_i(x_j)=q_{ij}x_j$ for some $q_{ij}\in
\CC^\times$, and $\Psi=\tau_\qsub$ with $\qq=(q_{ij})$. The algebra
$S(V,\Psi)$ is a flat deformation of $S(V)$ if and only if
$q_{ii}=q_{ij}q_{ji}=1$ for all $i$, $j$.

If $\lambda\ne 1$, the grading on $V$ must be given by $V=V_a$ for
some $a\in \Gamma$. 
Substitution $y=a^{-1}(x)$ gives $\bigl(a(x)+(1-\lambda)x\bigr)
\tensor x = \lambda x\tensor a^{-1}(x)$, hence $a^{-1}(x)$ is
proportional to $x$ for all $x\in V$. It follows that $a$ acts on $V$
as a scalar. The argument in the previous paragraph then applies
and leads to a contradiction with $\lambda\ne 1$.
\end{proof}

\subsection{$\qq$-polynomials and the $\qq$-Heisenberg double} 
\label{subsect:qpoly} 

From now on, let the $n\times n$ matrix
$$
 \qq = (q_{ij})_{i,j=1}^n,  
\qquad q_{ii}=q_{ij}q_{ji}=1 \quad \forall i,j
$$
be fixed, and 

\endcmm
$
   V = \mathrm{span}\, (x_1,\dots,x_n)
$
will be a space spanned by $n$ independent variables, and 
the symbol $\qq$ will stand for an $n\times n$ matrix of deformation  
parameters $q_{ij}$ 
such that $q_{ii}=1$, $q_{ij}q_{ji}=1$ for all $i,j=1,\dots,n$. 
The $\qq$\dash polynomial algebra
$$
S_\qsub(V) := \CC \langle \, x_1, \dots, x_n\, |
\, x_i x_j =
 q_{ij} x_j x_i \, \rangle
$$
is a flat deformation of the symmetric algebra $S(V)$. 
The space of quadratic relations in $S_\qsub(V)$ is
$$
\wedge^2_\qsub V = \mathrm{span}(x_i\tensor x_j - q_{ij} x_j \tensor
   x_i \mid 1\le i,j\le n) 
%= \ker(\id+\tau_\qsub)
   \quad \subset V\tensor V,
$$
the $\qq$\dash exterior square of $V$.
%Introduce the group
%$$
%  \Gamma_\qsub \subset \GL(V), \qquad \Gamma_\qsub=\langle \gamma_1,
%  \dots,\gamma_n\rangle, \qquad \gamma_i(x_j) = q_{ij} x_j
%$$
%as in the proof of Lemma~\ref{lem:tau}.
%We will be working with 
%the Koszul algebra 
%$$
%S_\qsub(V) := \CC \langle \, x_1, \dots, x_n\, |
%\, x_i x_j =
% q_{ij} x_j x_i \, \rangle,
%$$
%of $\qq$\dash polynomials, otherwise written as $S(V,\tau_\qsub)$. 
For future use, we denote by $y_1,\dots,y_n$ the basis of $V^*$ dual
to $\{x_i\}$, so that 
$
   V^* = \mathrm{span}\, (y_1,\dots,y_n)
$.
Furthermore, we introduce the Abelian group
$$
  \Gamma_\qsub \subset \GL(V), \qquad \Gamma_\qsub=\langle \gamma_1,
  \dots,\gamma_n\rangle, \qquad \gamma_i(x_j) = q_{ij} x_j,
$$
and observe that $V$ is a Yetter\dash Drinfeld module over
$\Gamma_\qsub$, via the natural action of $\Gamma_\qsub$ and the
grading 
$$
    |x_i|=\gamma_i. 
$$
This gives rise to the braiding $\tau_\qsub$ on $V$, defined by
$\tau_\qsub(x_i \tensor x_j) = q_{ij} x_j \tensor x_i$. 
The algebra $S_\qsub(V)$ coincides with the braided symmetric algebra 
$S(V,\tau_\qsub)$. 

\subsection{The $\qq$-Heisenberg double}

Our next step is to introduce a $\qq$\dash differential calculus 
via a $\qq$\dash analogue of the Weyl algebra $\mathcal D(V)$. 
We have two candidates for the role of such
$\qq$\dash analogue. One is the Heisenberg quadratic double 
$A_V$, associated to $V$ as a Yetter\dash Drinfeld module over
the group $\Gamma_\qsub$.
% (see the proof of Lemma~\ref{lem:tau}). 
The other candidate is the braided Weyl algebra $\rd A_V$ 
(Definition~\ref{def:braidedweyl}) of $V$, and this will be relevant
for braided Cherednik algebras later. The two are related by 
braided reduction, as described in Section~\ref{sect:qd}. 

In this Section we construct $\qq$\dash Cherednik algebras as
deformations of the Heisenberg quadratic double $A_V$. 
The presentation of $A_V$ is given in 

\begin{proposition}
\label{prop:q-Heisenberg}
Let $V$ be viewed as a Yetter\dash Drinfeld module over the group
$\Gamma_\qsub$ as above.
The Heisenberg quadratic double $A_\qsub:=A_V$ 
has the  triangular decomposition 
$$
    A_\qsub \cong S_\qsub(V) \tensor \CC \Gamma_\qsub \tensor 
                S_\qtran(V^*)\,,
$$
where $\qtran$ is the transpose of the matrix $\qq$, so that 
$S_\qtran(V^*)$ is generated by $V^*$ subject to relations 
$y_i y_j = q_{ji}y_j y_i$. 
The commutation relation between $x_i$ and $y_j$ is
$
      y_j x_i - x_i y_j = \delta_{i,j} \gamma_i \in \CC \Gamma_\qsub
$.
\end{proposition}
\begin{proof} 
%The triangular decomposition of $A_\qsub$ is a particular case of
%The generators and the relations in
%$S(V,\tau_\qsub)=S_\qsub(V)$ were
%  given in the proof of Lemma~\ref{lem:tau}. The $y_i$ are
%the generators of the algebra  $S(V^*,\tau_\qsub^*)$. 
%Based on Remark~\ref{rem:adjbraiding}, 
The braiding 
$\tau_\qsub^*$ on $V^*$ is computed as $\tau_\qsub^*(y_i\tensor
y_j)=q_{ji} y_j \tensor y_i$ (note the order of the indices). One
thus obtains the relations $y_iy_j-q_{ji}y_j y_i$ in
$S(V^*,\tau_\qsub^*)$ as the kernel of  
$\id+\tau_\qsub^*$.
% (see Lemma~\ref{lem:Heisenberg}).  
%Finally, by Lemma~\ref{lem:Heisenberg} the commutator of $y_j$ and $x_i$ is
%$\langle x_i, y_j\rangle |x_i|=\delta_{i,j} \gamma_i$. 
The claim then follows from Proposition~\ref{prop:heisenberg} and the
definition of the map $\beta_Y$ in Section~\ref{sect:qd}. 
\end{proof}

It is now natural to look for a $\qq$\dash analogue of
rational Cherednik algebras among quadratic doubles 
 with triangular decomposition 
$$
A\cong S_\qsub(V) \tensor \CC W \tensor S_\qtran(V^*) ,
$$
where $W$ is a subgroup of $\GL(V)$.
Our next step is to determine what the group $W$ can be.

\subsection{The subgroup of $\GL(V)$ that preserves $\wedge^2_\qsub V$
and $\wedge^2_\qtran V^*$}
\label{subsect:preserves}

Suppose that $W$ is a subgroup of $\GL(V)$ 
such  that there exists a
quadratic double of the form $S_\qsub(V) \tensor \CC W
\tensor S_\qtran(V^*)$. By Theorem~\ref{thm:relations},  
$$
    W(\wedge_\qsub^2 V) = \wedge_\qsub^2 V\,, 
    \qquad
    W(\wedge_\qtran^2 V^*) = \wedge_\qtran^2 V^*\,.
$$
In the case $q_{ij}=1$ $\forall i,j$  (the commutative case) these
conditions are vacuous; but they are not so in general. 
The group $W$ must be a subgroup of
$$
        N(\qq) = \{ w\in \GL(V) \mid w(\wedge_\qsub^2 V) =
        \wedge_\qsub^2 V, \ w(\wedge_\qtran^2 V^*) = \wedge_\qtran^2 V^* \}\,.
$$
To describe $N(\qq)$, we denote 
$$
\Symm(\qq) = \{\sigma\in \Symm_n \mid
q_{\sigma(i)\sigma(j)}=q_{ij} \ \forall i,j\}
$$
and view $\Symm(\qq)$ as a subgroup of $\GL(V)$ acting on $V$ by
permutations of the basis $\{x_i\}$. 
Recall the grading $V=\oplus_{\gamma\in \Gamma_\qsub} V_\gamma$ given by
$|x_i|=\gamma_i$ and 
observe that the component $V_{\gamma_i}$ of $V$ is spanned by 
$\{x_j \mid q_{jk}=q_{ik}\ \text{for all indices}\ k\}$. 
Denote 
$$
  L(\qq) = \{g\in \GL(V) \mid g(V_\gamma)=V_\gamma\ \text{for all}\ 
  \gamma\in\Gamma_\qsub\}.
$$
Clearly, the group $\Symm(\qq)$ normalises $L(\qq)$, therefore
$\Symm(\qq)\cdot L(\qq)$ is a subgroup in $\GL(V)$. 
\begin{proposition}
\label{prop:autom}
1. $N(\qq)= \Symm(\qq)\cdot L(\qq)$. 

2. $N(\qq)$ is the stabiliser of the set $\{\gamma_i | i=1,\dots,n\}$
   in $\GL(V)$ and is the normaliser of $\Gamma_\qsub$ in $\GL(V)$.
\end{proposition}
The Proposition will follow from two elementary lemmas.
\begin{lemma}
\label{lem:elem1}
Let $c_i$, $c_i'$ $(i=1,\dots,n)$ and $q$ be scalars, and let
$x=\sum_i c_i x_i$, $x'=\sum_i c_i' x_i$. Then 
$x\tensor x'-qx'\tensor x\in \wedge^2_\qsub V$ if and only if 
$(1-qq_{ij})c_i' c_j = (q-q_{ij})c_i c_j'$ for all indices $i, j$. 
\end{lemma}
\begin{proof}
Recall that $\wedge^2_\qsub V=\ker(\id+\tau_\qsub)$.
Applying $\id+\tau_\qsub$ to $x\tensor x'-qx'\tensor x$ and equating
the coefficient of $x_j\tensor x_i$ with zero gives the desired identity. 
\end{proof}
\begin{lemma}
\label{lem:elem2}
Let $A$ be an endomorphism of $V$ acting by $Ax_i = \sum_j A_i^j
x_j$. Then $A$ preserves $\wedge^2_\qsub V$ and the adjoint $A^*$ of
$A$ preserves $\wedge^2_\qtran V^*$, if and only if 
$$
   (q_{kl}-q_{ij})A_k^i A_l^j=0
$$ 
for all indices $i,j,k,l$. 
\end{lemma}
\begin{proof}
The condition that $A$ preserves $\wedge^2_\qsub V$ is equivalent to 
$Ax_k\tensor Ax_l - q_{kl} Ax_l \tensor Ax_k\in \wedge^2_\qsub V$ for
all $k$, $l$. By Lemma~\ref{lem:elem1} this is the same as 
$(1-q_{kl}q_{ij})A_l^i A_k^j = (q_{kl}-q_{ij})A_k^i A_l^j$ for all
$i,j,k,l$. The matrix of $A^*$ with respect to the basis $\{y_i\}$
dual to $\{x_i\}$ is the transpose of the matrix of $A$,
and $\qtran$ is the transpose of $\qq$, therefore in $V^*$ we obtain the
condition 
$(1-q_{lk} q_{ji})A_i^l A_j^k = (q_{lk}-q_{ji})A_i^k A_j^l$. 
Swapping the indices $i$, $l$, as well as $j$, $k$, we get 
$(1-q_{kl}q_{ij})A_l^i A_k^j = -(q_{kl}-q_{ij})A_k^i
A_l^j$. Hence both sides of this equation are zero. 
\end{proof}
\begin{proof}[Proof of Proposition \ref{prop:autom}]
1. It is clear that $\Symm(\qq)$ and $L(\qq)$ both preserve
$\wedge^2_\qsub V$ and $\wedge^2_\qtran V^*$, thus it is enough
to show that  $N(\qq)\subset \Symm(\qq) \cdot L(\qq)$. 
Let $w\in N(\qq)$. By definition of the action of $\GL(V)$ on $V^*$,
the action of $w^{-1}$ on $V^*$ is by the adjoint $w^*$ of $w$, 
thus $w^*$ preserves $\wedge^2_\qtran V^*$.
By Lemma~\ref{lem:elem2}, 
$(q_{kl}-q_{ij})w_k^i w_l^j =0$ for all indices $i,j,k,l$, where 
$w_k^i$ are entries of the matrix of $w$ in the basis $\{x_i\}$. 
By invertibility of
$w$, there exists a permutation $\sigma\in \Symm_n$ such that
$w_i^{\sigma(i)}\ne 0$ for all indices $i$. For any pair $i,j$ of
indices one has the relation $(q_{ij} - q_{\sigma(i)\sigma(j)})
w_i^{\sigma(i)} w_j^{\sigma(j)}=0$, hence $q_{ij} =
q_{\sigma(i)\sigma(j)}$ and $\sigma\in \Symm(\qq)$.

We are left to prove that the matrix $g:=\sigma^{-1} w$, with entries
also satisfying the equation in Lemma~\ref{lem:elem2}  
and with $g_i^i\ne 0$ for
all $i$, is in $L(\qq)$; equivalently, that 
$\gamma_i\ne \gamma_j\in \Gamma_\qsub$
implies $g_i^j=0$. Indeed, find $l$ such that $q_{il}\ne q_{jl}$. 
The relation $(q_{il}-q_{jl})g_i^j g_l^l=0$ implies that 
$g_i^j=0$ as required.

2. For $g\in \GL(V)$, $g\,\Gamma_\qsub\, g^{-1}=\Gamma_\qsub$ if and only
   if $g$ permutes the simultaneous eigenspaces of
   $\Gamma_\qsub$. These are the same as simultaneous eigenspaces of
   $\gamma_i$, i.e., the subspaces $V_{\gamma_j}$ of $V$. 
   It is obvious that such $g$ are precisely elements of $\Symm(\qq)L(\qq)$. 
\end{proof}

\begin{corollary}
\label{cor:yd}
Let $W$ be a subgroup of $\GL(V)$ that contains $\Gamma_\qsub$. 
The group $W$ preserves $\wedge^2_\qsub V$ and $\wedge^2_\qtran V^*$  
if and only if $V$ is a Yetter\dash
Drinfeld module via the $W$\dash action on $V$  
and the $W$\dash grading by $|x_i|=\gamma_i$. 
\qed
\end{corollary}  
\begin{remark}
An element $w\in \GL(V)$ stabilises $\wedge^2_\qtran V^*$, if and only
if $w$ stabilises the $\qq$\dash symmetric square $S_\qsub^2
(V):=\mathrm{span}_{i,j} (x_i\tensor x_j + q_{ij} x_j \tensor x_i)$ of
$V$. This is because $S_\qsub^2(V)$ is the orthogonal complement  of
$\wedge^2_\qtran V^*$ with respect to the standard pairing. Note that
$S^2_\qsub(V)$, $\wedge^2_\qsub V$ are the eigenspaces of the involutive
braiding $\tau_\qsub$ on $V\tensor V$. Therefore, 
$N(\qq)$ is the centraliser of $\tau_\qsub$ in $\GL(V)$.
\end{remark}

\subsection{$\qq$-Cherednik algebras}

To obtain a nice classification of deformations
of the Heisenberg quadratic double $A_\qsub$, we impose an extra
nondegeneracy condition:

\begin{definition}
A quadratic double $A\cong T(V)/\lgen R^-\rgen \tensor \CC W \tensor
T(V^*)/\lgen R^+ \rgen$ is called \emph{non\dash degenerate}, 
if  the commutator map $[\cdot,\cdot]\colon V^* \times V \to \CC W$
has no non\dash trivial kernels in~$V^*$ and in~$V$. 
\end{definition}

\begin{definition}
A \emph{$\qq$\dash Cherednik algebra} is a non\dash degenerate
quadratic double  with triangular decomposition $S_\qsub(V)\tensor \CC W \tensor
S_\qtran(V^*)$, where $W$ is a subgroup of $\GL(V)$ (not necessarily
finite).   
\end{definition}

%\begin{remark}
%Our motivation for the nondegeneracy requirement is as follows. 
%In the commutative case ($q_{ij}=1$ for all $i$, $j$), 
%\cite[Theorem 1.3]{EG} guarantees that if $V$ is an irreducible
%$W$\dash module, quadratic doubles of the form $S(V)\tensor \CC W
%\tensor S(V^*)$ are rational Cherednik algebras $H_{t,c}(W)$. 
%Irreducibility of $W$ implies that
%doubles are non\dash degenerate, because the left and right kernels of
%the $W$\dash equivariant commutator $[\cdot,\cdot]$ are 
%$W$\dash submodules of $V$, $V^*$.
%Moreover, rational Cherednik algebras $H_{t,c}(W)$ with $t\ne 0$ are
%non\dash degenerate in the sense of our definition. Thus,
%nondegeneracy is a good way to single out quadratic doubles that are to
%be $\qq$\dash analogues of rational Cherednik algebras. 
%\end{remark}

In the next Proposition, we keep the notation for $\qq$, $x_i$, $y_i$,
$V$ and $W$. 
\begin{proposition}
\label{prop:qcomm}
A $\qq$\dash Cherednik algebra is generated by 
$x_1,\dots,x_n\in V$, $w\in W$ and $y_1,\dots,y_n\in V^*$
subject to relations
\begin{itemize}
\item
$x_i x_j = q_{ij} x_j x_i$,\quad $y_i y_j = q_{ji} y_j y_i$,\quad 
$     w x_i w^{-1} = w(x_i) \in V$,\quad $w y_i w^{-1} = w(y_i) \in V^*$,
and
\item    $y_j x_i - x_i y_j = \sum_{w\in W} \langle L_w(x_i), y_j
    \rangle w$ for some $L_w\in \End(V)$.
\end{itemize}
The maps $L_w$ are such that $\cap_{w\in W}\ker L_w = 0$, 
$\cap_{w\in W}\ker L_w^* = 0$, and satisfy
\begin{align*}
    g(L_w(g^{-1}(x)))  = L_{gwg^{-1}}(x) 
& \qquad (\text{$W$\dash equivariance});
\\
     (x_i - q_{ij} w(x_i)) \tensor L_w(x_j) = (q_{ij}x_j
     -w(x_j))\tensor L_w(x_i) & 
\\
     (y_i - q_{ji} w(y_i)) \tensor L_w^*(y_j) = (q_{ji}y_j
     -w(y_j))\tensor L_w^*(y_i)
&
\qquad(\text{$\qq$\dash commutativity equations}) 
\end{align*}
for all $g,w\in W$ and all indices $i$, $j$.
Conversely, an algebra with the above presentation, with $W\le \GL(V)$
centralising $\tau_\qsub$ and $L_w$ subject
to the above conditions, is a $\qq$\dash Cherednik algebra. 
\end{proposition}

\begin{proof}
The defining relations follow from the definition of a $\qq$\dash
Cherednik algebra as a quadratic double, 
while $\cap_{w\in W}\ker L_w = 0$, $\cap_{w\in W}\ker L_w^* = 0$ is
precisely the nondegeneracy condition. 
Furthermore, the $W$\dash equivariance condition in the Proposition is
the same as $W$\dash equivariance of the commutator, as required by
Theorem~\ref{thm:equivar}. It remains to show that the $\qq$\dash
commutativity equations are equivalent to the conditions $[R^+,V]=0$
and $[V^*,R^-]=0$ in Theorem~\ref{thm:relations}, where
$R^+=\wedge^2_\qsub V^*$ and $R^-=\wedge^2_\qsub V$. To analyse 
the commutator $[V^*, R^-]$, write
\begin{align*}
[y,x_i x_j -q_{ij}x_j x_i] &= x_i [y,x_j] - q_{ij} [y,x_j] x_i +
[y,x_i]x_j-q_{ij} x_j[y,x_i]
\\
&= \sum_{w\in W}
\bigl(x_i w - q_{ij} w x_i \bigr) \langle L_w(x_j), y\rangle 
-\bigl( q_{ij} x_j w - w x_j\bigr)\langle L_w(x_i)), y\rangle 
\\
&=\sum_{w\in W}
(x_i  - q_{ij} w(x_i))\cdot w \cdot \langle L_w(x_j), y\rangle 
-( q_{ij} x_j - w(x_j)))\cdot w\cdot\langle L_w(x_i)), y\rangle, 
\end{align*}
which vanishes for all $y\in V^*$ if and only if the first $\qq$\dash
commutativity equation holds. 
Similarly, $[R^+,V]=0$ is equivalent to the second $\qq$\dash
commutativity equation. 
\end{proof}

\subsection{The block structure of the matrix $\qq$}

The structure of the subgroup $\Symm(\qq)$ of $\Symm_n$ and its action
on the space $V$ may be complicated, depending  on
the combinatorics of the matrix~$\qq$. We will soon show, however, that 
the $\qq$\dash commutativity equations in Proposition~\ref{prop:qcomm}
imply that only the part of $\Symm(\qq)$ generated by transpositions
actually matters for $\qq$\dash Cherednik algebras. This leads to $V$ and
the matrix $\qq$ being split into blocks; let us formally introduce
this block structure. 

\begin{definition}
We say that indices $i,j\in \{1,\dots,n\}$ are in the same block (with
respect to the matrix $\qq$), if 
$$
   q_{ik}=q_{jk} \quad \text{for all}\ k\ne i,j; 
   \qquad 
   q_{ij}=\pm 1.
$$
\end{definition} 
\begin{lemma}
\label{lem:blocks}
Being in the same block is an equivalence relation on the set 
$\{1,\dots,n\}$ of indices. An equivalence class $B$ (a block
of indices)  can be of one of the following two types: 
\begin{itemize}
\item[---]
positive  block: $q_{ij}=1$\phantom{$-$} for all $i,j\in B$;
\item[---]
negative block: $q_{ij}=-1$ for all $i,j\in B$, $i\ne j$, where $|B|>1$.
\end{itemize}
\end{lemma} 
\begin{proof}
Let us write $i\sim_+ j$, respectively $i\sim_- j$, if $i$, $j$ are indices
such that $q_{ik}=q_{jk}$ for any $k\ne i,j$ and $q_{ij}=1$ (respectively
$q_{ij}=-1$). We need to check that the relation $\sim \  =\ \sim_+\cup
\sim_-$ is an equivalence relation. 
Note that $\sim_+$ is an equivalence relation, because 
$i\sim_+ j$ means that rows $i$ and $j$ of the matrix $\qq$ are
identical. 
Hence $\sim$ is reflexive, and is symmetric as both $\sim_+$ and
$\sim_-$ are. Moreover, $\sim_+$ is transitive, therefore it 
remains to check that $a\sim_- b\sim c$ implies $a\sim c$. If $c=a$
or $c=b$, we are done, otherwise $a\sim_- b$ implies $q_{ac}=q_{bc}$
and $b\sim c$ implies $q_{ba}=q_{ca}$. Since $q_{ba}=-1$, we have 
$q_{ac}=q_{bc}=-1$ and $q_{ab}=q_{cb}=-1$. Finally, for any $k\ne
a,b,c$ we have $q_{ak}=q_{bk}=q_{ck}$. Thus, $a\sim_- c$. 
\end{proof}
\begin{corollary}[Block structure of the matrix $\qq$]
\label{cor:block}
Let the matrix $\qq$ be given.
The index set $\{1,\dots,n\}$ is split into disjoint 
blocks. To each pair $B$, $C$ of blocks there
is associated a complex number $q_{B,C}=q_{C,B}^{-1}\in \CC$ such that 
$$ 
   q_{ij} = q_{B,C} \quad\text{whenever}\ i\ne j,\ \ i\in B, \ j\in C.
$$
In particular, $q_{B,B}$ is $1$ or $-1$, depending on whether the block $B$ is
positive or negative. 
\qed
\end{corollary}

Let $B$ be a block of indices. Introduce the following subspaces:
$$
  V_B = \mathrm{span}\, (x_i \mid i\in B)\, \quad \subset V,
\qquad
   V^*_B = \mathrm{span}\, (y_i \mid i\in B)\, \quad \subset V^*, 
$$
and let $\gamma_{B}\in N(\qq)\subset \GL(V)$ be such that
$$
   \gamma_\Bsub|_{V_C}=q_{B,C} \id_{V_C}
$$
for any block $C$, where the scalars $q_{B,C}$ are as introduced in
Corollary~\ref{cor:block}.

\subsection{$\qq$-Cherednik algebras: the structure theorem}

Proposition~\ref{prop:qcomm} gives the relations in a $\qq$\dash
Cherednik algebra explicitly, except the most important one --- 
the commutation relation between $V^*$ and $V$. It turns out that, similar
to rational Cherednik algebras, the
commutator is expressed in terms of complex reflections in the group
$W$,  but premultiplied with elements $\gamma_B$ as an extra ingredient. 
For reference, we need a list of 
complex reflections in $\GL(V)$ that preserve the relations in the
algebras $S_\qsub(V)$ and $S_\qtran(V^*)$. 

\begin{lemma}
\label{lem:compref}
Let $s\in \GL(V)$ be a complex reflection (not necessarily of finite
order) and 
$\alpha_s\in V$, $\alpha_s^\vee \in V^*$ be the
root\dash coroot pair for $s$. 
If $s\in N(\qq)$, then:

$(1)$ There is a block $B\subset \{1,\dots,n\}$ of indices,
such that 
$\alpha_s\in V_{B}$ and $\alpha_s^\vee\in V^*_{B}$.

$(2)$ If $B$ is a positive block, $s$ is an arbitrary complex
reflection in the space $V_{B}$. 

$(3)$ If $B$ is negative, 
$s$ must be of the form  $\refl_i^{(\eta)}$ ($\eta\ne 0,1$) 
or $(ij)\refl_i^{(\varepsilon)}\refl_j^{(\varepsilon^{-1})}$.
Here $(ij)$ permutes variables $x_i$ and $x_j$ with $i,j\in B$, and 
$\refl_i^{(\varepsilon)}$ multiplies the variable $x_i$ by 
$\varepsilon\in \CC^\times$, leaving the rest of the variables intact. 
\end{lemma}
\begin{proof}
$(1)$ 
By Proposition~\ref{prop:autom}, $s=\overline s \cdot g$ where
$\overline s$ is a permutation of indices such that $q_{\overline
  s(i)\overline s(j)} = q_{ij}$, and $g$ preserves all
$\Gamma_\qsub$\dash graded components $V_{\gamma_k}$ of $V$. 
It follows that 
$(1-s)V_{\gamma_k} \subset V_{\gamma_k}+V_{\gamma_{\overline s(k)}}$. 
But then $\dim(1-s)V=1$ implies that there are at most two indices 
$k$ such that $(1-s)V_{\gamma_k}\ne 0$. If there is only one such
index $k$, let $B$ be the block of indices containing $k$.
Otherwise, there are two such indices $i$, $j$, and
necessarily $\overline s$ is the permutation $(i j)$. 
One has $q_{ji}=q_{ij}$ (hence $q_{ij}=\pm 1$) and
$q_{ia}=q_{ja}$ for all $a\ne i,j$, thus $i$, $j$ belong to the same
block; let $B$ be the block which contains $i$, $j$. 
In either case, $(1-s)V_{B'}=0$ for $B'\ne B$ and
$(1-s)V_{B}\subset V_{B}$, which implies $\alpha_s\in V_B$ and
$\alpha_s^\vee\in V^*_B$. 
 
Part $(2)$ is clear, as $V_{B}$ is of the form $V_{\gamma_k}$ if
$B$ is a positive block, thus any  complex reflection $s$ in $V_B$
has decomposition $\id\cdot s\in \Symm(\qq)\cdot L(\qq)$ and 
hence commutes with $\tau_\qsub$. 

Finally, if $B$ is a negative
block, then by Lemma~\ref{lem:blocks} 
$V_B = \oplus \{ V_{\gamma_i} \mid i\in B\}$ is a direct
sum of one\dash dimensional $\Gamma_\qsub$\dash graded components. 
By Proposition \ref{prop:autom}, $s$ must act imprimitively and permute
these $1$\dash dimensional subspaces. All such imprimitive complex
reflections are listed in $(3)$, cf.\ \cite[Section 3]{DO}.
\end{proof}

The following Theorem completes the description of 
 the structure of $\qq$\dash Cherednik algebras. 

\begin{theorem}
\label{thm:structure}
Let $A \cong S_\qsub(V)\tensor \CC W \tensor
S_\qtran(V^*)$ be 
a $\qq$\dash Cherednik algebra. Then the commutator of $y\in V^*$ and
$x\in V$ in $A$ is of the form
$$
  y x  - x y = \sum_{\mathrm{blocks}\ B}\gamma_\Bsub\cdot 
                   \Bigl( \, (x, y)_B\cdot 1 + 
                   \sum_s  c_s \langle x,\alpha_s^\vee\rangle
\langle \alpha_s,y\rangle s \, \Bigr) , 
$$
where the sum is taken over complex reflections $s$ which commute with
the braiding $\tau_\qsub$ and such that $\alpha_s\in V_B$, $\alpha_s^\vee\in V^*_B$
and $\gamma_\Bsub s \in W$. The pairing $(\cdot , \cdot)_B$ between
$V$ and $V^*$ is such that
$(x_i,y_j)_B=0$ unless $i,j\in B$, and is so chosen, together with
the constants $c_s$, as to make the commutator $W$\dash equivariant and
non\dash degenerate.  
\end{theorem}
\begin{corollary}
\label{cor:zerocomm}
In particular, in a $\qq$\dash Cherednik algebra one has $y_j x_i -
x_i y_j=0$ unless $i$, $j$ are in the same block of indices with
respect to $\qq$.
\end{corollary}
\begin{proof}[Proof of the Theorem]
We write the commutator as $[y,x]=\sum_{w\in W} \langle L_w(x),
y\rangle w$ 
with $L_w\in \End(V)$. 
%By Proposition~\ref{prop:autom} each $w\in W$ 
%decomposes as $w=\overline w \cdot g$, with $\overline w\in \Symm_\qsub$
%and $g\in L(\qq)$. The permutation $\overline w$ may not be unique,
%and is defined up to a permutation of indices which does not change
%the $\Gamma_\qsub$\dash grading. 
It is enough to show that if the map
$L_w$ is non\dash zero, then either
\begin{itemize}
\item[$(a)$]  $w=\gamma_\Bsub$ for some block $B$,  
$L_w(V_{B'})=0$ for blocks $B'\ne B$, and $L_w(V_B)\subset V_B$;\quad
or
\item[$(b)$]
$w=\gamma_\Bsub s$ for a complex reflection $s$ such that $\alpha_s\in
V_B$, $\alpha_s^\vee\in V^*_B$, and $L_w(x) = \mathrm{const}\cdot 
\langle x,\alpha_s^\vee\rangle \alpha_s$. 
\end{itemize}

\subsection*{Case 1: $w$ preserves each $\Gamma_\qsub$-homogeneous
  component $V_{\gamma_i}$ of $V$} 

Find an index $i$ such that 
$L_w(x_i)\ne 0$. 
For an index $j$ such that $\gamma_j\ne \gamma_i$,  the vectors
$x_i-q_{ij}w(x_i)\in V_{\gamma_i}$ and $q_{ij} x_j - w(x_j)\in
V_{\gamma_j}$ cannot be nonzero and proportional, therefore, both sides of the
$\qq$\dash commutativity equation in Proposition~\ref{prop:qcomm} must
be zero.  
It follows that $q_{ij} x_j - w(x_j)=L_w(x_j)=0$. Hence
$w|_{V_{\gamma_j}}=\gamma_i|_{V_{\gamma_j}}$ and $L_w|_{V_{\gamma_j}}=0$  on all
$\Gamma_\qsub$\dash homogeneous
components $V_{\gamma_j}$ of $V$ such that $\gamma_j\ne \gamma_i$.
Similarly, $L_w^*$ vanishes on $\Gamma_\qsub$\dash homogeneous
components of $V^*$ other than $(V^*)_{\gamma_i^{-1}}$, which means
that $L_w(V_{\gamma_i})\subset V_{\gamma_i}$. 

Furthermore, Proposition~\ref{prop:qcomm} implies that for $x,x'\in
V_{\gamma_i}$ 
\begin{equation*}
\tag{$*$}
\label{eq:case1}
  (x-w(x))\tensor L_w(x') = (x'-w(x'))\tensor L_w(x).   
\end{equation*}
It is easy to see that this tensor equation may hold only either if $w=\id$
on $V_{\gamma_i}$, or if $\dim (1-w) V_{\gamma_i} = \dim
L_w(V_{\gamma_i})=1$. In the former case, $w=\gamma_i$. If $i$ belongs to a
positive block $B$, one has $V_B = V_{\gamma_i}$ and therefore
$w=\gamma_\Bsub$, so that option $(a)$ holds. If $i$ belongs to a negative block $B$, then $\dim V_{\gamma_i}=1$, therefore
$L_w(x)=\mathrm{const}\cdot \langle x, y_i\rangle x_i$.  The element $w=\gamma_i$
decomposes as $\gamma_\Bsub \refl_i^{(-1)}$, and $\refl_i^{(-1)}$ is a complex
reflection on $ V_{\gamma_i}=\CC x_i$ with root\dash coroot pair
$x_i$, $2y_i$, so that option $(b)$ holds. 

In the case $\dim (1-w) V_{\gamma_i} = \dim
L_w(V_{\gamma_i})=1$, the element $w$ is necessarily $\gamma_i s$,
where $s$ is a complex reflection on $V_{\gamma_i}$. Write
$L_w(x)=\langle x, \alpha\rangle \beta$ with $\alpha\in
V_{\gamma_i^{-1}}^*$ and $\beta\in V_{\gamma_i}$. By
(\ref{eq:case1}), $L_w(x)$ vanishes on $\ker(1-s)$, therefore
$\alpha=\mathrm{const}\cdot \alpha^\vee_s$. Moreover,
Proposition~\ref{prop:qcomm} implies the 
equation $(y-s(y))\tensor L_w^*(y') = (y'-s(y'))\tensor L_w^*(y)$ for
$y,y'\in (V^*)_{\gamma_i^{-1}}$, so that $\beta=\mathrm{const}\cdot
\alpha_s$. If $i$ belongs to a
positive block $B$, we have $V_B=V_{\gamma_i}$, and option $(a)$
holds. If  $i$ belongs to a
negative block $B$ so that $\dim V_{\gamma_i}=1$, 
then $s=\refl_i^{(\eta)}$ for some 
$\eta\ne 1$, and $w=\gamma_i\refl_i^{(\eta)}= \gamma_\Bsub \refl_i^{(-\eta)}$. 
No matter what $\eta$ is, the root and the coroot of 
$\refl_i^{(-\eta)}$ are proportional to $x_i$ and $y_i$, respectively,
hence option $(b)$ still holds.

\subsection*{Case 2: there exist indices $i$, $j$ such that
  $\gamma_i\ne \gamma_j$ and $w(V_{\gamma_i})=V_{\gamma_j}$}
By Proposition~\ref{prop:qcomm},
$$
    (x_i -q_{ik} w(x_i)) \tensor L_w(x_k) = (q_{ik} x_k -w(x_k))
    \tensor L_w(x_i)\qquad
\text{for all}\ k.
$$
Note that $x_i -q_{ik} w(x_i)$ cannot be zero, because $x_i\in
V_{\gamma_i}$, $w(x_i)\in V_{\gamma_j}$ and $ V_{\gamma_i} \cap
V_{\gamma_j}=0$. 
Therefore $L_w(x_i)\ne 0$, as otherwise the commutativity equation
would imply that $L_w(x_k)=0$ for all $k$.
Now observe that for any $x\in V_{\gamma_j}$, $x\ne 0$ one has 
\begin{equation*}
\label{eq:case2}
\tag{$**$}
      (x_i -q_{ij} w(x_i)) \tensor L_w(x) = (q_{ij} x -w(x))
    \tensor L_w(x_i)
\end{equation*}
and $q_{ij}x -w(x)\ne 0$ because $w(V_{\gamma_j})\cap V_{\gamma_j}=0$.
It follows that $L_w(x)$ is proportional to $L_w(x_i)$ for any $x\in
V_{\gamma_j}$, thus $\dim L_w(V_{\gamma_j})=1$. 

Now if the $\dim V_{\gamma_j}$ is greater than $1$, the map
$L_w$ must have a kernel in $V_{\gamma_j}$. Pick $0\ne x\in V_{\gamma_j}$ such
that $L_w(x)=0$. Substituting $x$ in (\ref{eq:case2}) 
leads to a contradiction, as the left-hand side
of the equation becomes zero while the right-hand side does
not. Therefore $\dim V_{\gamma_i} = \dim V_{\gamma_j} = 1$ 
and $w(x_i)=\varepsilon x_j$ for some $\varepsilon \in \CC^\times$.

By (\ref{eq:case2}),   $x_i - q_{ij} \varepsilon x_j$ is
proportional to $q_{ij} x_j - w(x_j)$. It follows that
$w(x_j)=\varepsilon^{-1}x_i$. 
%and that $\overline w$ permutes the indices $i$ and
%$k$.
Hence for $l\ne i,j$ the subspace $w(V_{\gamma_l})$ has zero
intersection with $V_{\gamma_i}$ and with $V_{\gamma_j}$, thus
the vector $q_{il} x_l -w(x_l)$ cannot coincide with $x_i -q_{il}
w(x_i)$ up to a non\dash zero factor. 
Equation (\ref{eq:case2})  therefore forces $q_{il} x_l -w(x_l)=0$ and
$L_w(x_l)=0$  for each $l\ne i,j$. A similar equation in $V^*$ forces 
$L_w^*(y_l)=0$  for  $l\ne i,j$, hence $L_w(\CC x_i + \CC x_j) \subset
\CC x_i + \CC x_j$. 

Furthermore, Proposition~\ref{prop:autom} implies that
$q_{ik}=q_{ki}$ and that 
$q_{il}=q_{jl}$ for all $l\ne i,j$. As $\gamma_i\ne \gamma_j$, one has 
$q_{ij}=-1$. This means that $i$ and $j$ belong to a negative
block $B$, and $w$ acts as $\gamma_\Bsub$ on each $V_{\gamma_l}$ with
$l\ne i,j$. We have the decomposition
$$
     w = \gamma_\Bsub \cdot (ij)\refl_i^{(\varepsilon)}\refl_j^{(\varepsilon^{-1})}.
$$ 
Now (\ref{eq:case2}) reads 
$$
      (x_i + \varepsilon x_j) \tensor L_w(x_j) = - (\varepsilon^{-1} x_i +
      x_j ) \tensor L_w(x_i)\,,
$$
so that $L_w(x_i) = -\varepsilon L_w(x_j)$. Moreover, the $W$\dash
equivariance condition in Proposition~\ref{prop:qcomm} implies that
$L_w$ commutes with $w$, whence $w(L_w(x_i))=\varepsilon L_w(x_j) =
-L_w(x_i)$, thus $L_w(x_i) = x_i - \varepsilon x_j$. It follows that 
$$
     L_w(x) = \mathrm{const}\cdot 
     \langle x, \ y_i - \varepsilon^{-1} y_j \rangle (x_i - \varepsilon x_j)
$$
and option $(b)$ holds. 
\end{proof}

\begin{corollary}
\label{cor:structure}
Let $W$ be a  subgroup of $\GL(V)$ centralising the braiding
$\tau_\qsub$. 
An algebra given by generators and relations
from Proposition~\ref{prop:qcomm} and the commutation relation from
Theorem~\ref{thm:structure} is a $\qq$\dash Cherednik algebra, if the
commutator is $W$\dash equivariant and non\dash degenerate. 
\end{corollary}
\begin{proof}
Indeed, we checked in the proof of Theorem~\ref{thm:structure} 
that the $\qq$\dash commutativity equations in
Proposition~\ref{prop:qcomm} are
satisfied. 
\end{proof}

%%%%%%%%%%%%%%%%%%%%%%%%%%%%%%%%%%%%%%%%%%%%%%

\section{Braided Cherednik algebras}

In this Section, we introduce braided Cherednik algebras. 
%and prove that they are 
%braided reductions of $\qq$\dash Cherednik
%algebras. However, the main result of this Section is that all braided
%Cherednik algebras are 
%obtained as braided tensor products of algebras
%attached to irreducible groups. 
Besides the well-known 
rational Cherednik algebras of Etingof and Ginzburg, irreducible groups
give rise to a  new class of negative braided Cherednik algebras. 

\subsection{The $\qq$-Weyl  algebra $\protect\rd A_\qsub$} 

Recall from Section~\ref{sect:q_cher} that the classical 
Weyl algebra of polynomial differential operators on the space
$V$ admits two possible
$\qq$\dash versions. One of them is the Heisenberg quadratic double
$A_\qsub \cong S_\qsub(V)\tensor \CC \Gamma_\qsub
\tensor S_\qtran(V^*)$ 
%of the Yetter\dash Drinfeld module $V=\mathrm{span}(x_1,\dots,x_n)$ 
over the group $\Gamma_\qsub$; 
we introduced $\qq$\dash Cherednik algebras 
as deformations of this. The other is the braided Weyl algebra
of $V$, obtained from $A_\qsub$ via braided reduction. 
We denote it by $\rd A_\qsub$ and will now review it in more detail. 
Note the appearance of the $\qq$\dash
symmetric algebra 
$$
   S_\qsub(V^*) := \CC\langle \uy_1,\dots,\uy_n \mid \uy_i \uy_j =
   q_{ij} \uy_j \uy_i \rangle\,,
$$
which is not the same as $S_\qtran(V^*)$ used in the previous
Section; in fact, $ S_\qsub(V^*) \cong S_\qtran(V^*)^{\mathrm{op}}$.

\begin{proposition}
Let $V$ be viewed as a Yetter\dash Drinfeld module over the group
$\Gamma_\qsub$ as above. The braided Weyl algebra 
$\rd A_\qsub := \rd A_V$ decomposes as 
$$
       \rd A_\qsub \cong S_\qsub(V)\tensor 
                         S_\qsub(V^*) \,,
$$
where $\uy_1\dots,\uy_n$ are a basis of $V^*$ dual to $\{x_i\}$, 
and the commutation relation between $\uy_j$ and $x_i$ is given by 
$$
      \uy_j x_i - q_{ij} x_i \uy_j = \delta_{i,j}\,.
$$
\end{proposition}
\begin{proof}
Follows immediately from Proposition~\ref{prop:q-Heisenberg} and
Proposition~\ref{prop:twist}. Alternatively, can be deduced from
Example~\ref{ex:braidedweyl}. 
\end{proof}

%\subsection{Yetter-Drinfeld modules over $\Gamma_\qsub$}
%\begin{remark}

We can view the space $V\oplus V^*$ as a Yetter\dash
Drinfeld module over the group $\Gamma_\qsub$ (a direct sum of two YD
modules) and denote the resulting braiding on $V\oplus V^*$ again by
$\tau_\qsub$. Then one has the braided commutator 
$$
    [a,b]_\qsub := a\tensor b - \tau_\qsub(a\tensor b), \qquad a,b\in
    V\oplus V^*.
$$
The $\Gamma_\qsub$\dash grading on $V\oplus V^*$ is given by 
$|x_i|=\gamma_i$, $|\uy_i|=\gamma_i^{-1}$, and recall that
$\gamma_i(x_j)=q_{ij} x_j$, $\gamma_i(\uy_j) = q_{ji} \uy_j$. 
Hence, the $\qq$\dash commutator is explicitly written as 
$$
   [x_i,x_j]_\qsub = x_i \tensor x_j - q_{ij} x_j \tensor x_i,
\qquad
   [\uy_i,\uy_j]_\qsub = \uy_i \tensor \uy_j - q_{ij} \uy_j \tensor \uy_i, 
\qquad
   [\uy_j,x_i]_\qsub = \uy_j \tensor x_i - q_{ij} x_i \tensor \uy_j.
$$
Let $\omega(a,b)$ be the skew\dash symmetric bilinear form 
on $V\oplus V^*$ uniquely determined by $\omega(x,x')=\omega(y,y')=0$,
$\omega(x,y)=\langle x,y\rangle$ for $x,x'\in V$, $y,y'\in V^*$. The
$\qq$\dash Weyl algebra can be 
defined as 
$$
    \rd A_\qsub = T(V\oplus V^*) / \lgen  [a,b]_\qsub - \omega(a,b) \rgen.
$$ 
Moreover, any subgroup $W\le \GL(V)$ which preserves the $\qq$\dash
deformed exterior squares $\wedge_\qsub^2 V$ and $\wedge_\qtran^2 V^*$,
will also preserve $\wedge_\qsub^2 V^* = \tau(\wedge_\qtran^2 V^*)$ and
centralise the braiding $\tau_\qsub$ on $V\oplus V^*$. Trivially, $W$
preserves the form $\omega$, therefore the $\qq$\dash Weyl algebra
$\rd A_\qsub$ will be a $W$\dash module algebra. 

\subsection{Braided Cherednik algebras}

Informally, one can now interpret $\omega$ in the above presentation
of the braided Weyl algebra $\rd A_\qsub$  as a $\CC W$\dash valued form. 
This leads to a braided
version of Drinfeld's degenerate affine Hecke algebra \cite{D} and
Etingof\dash Ginzburg symplectic reflection algebra \cite{EG}, and is
a natural way to introduce braided Cherednik algebras.
Their formal definition is as follows.

\begin{definition}
\label{def:bcher}
A \emph{braided Cherednik algebra} associated to the matrix $\qq$ is
an algebra with triangular decomposition $S_\qsub(V) \tensor \CC W
\tensor S_\qsub(V^*)$ where
$$
     wxw^{-1}=w(x), \qquad wyw^{-1} = w(y), 
\qquad
     [y,x]_\qsub \in \CC W
$$
for $x\in V$, $y\in V^*$, $w\in W$, such that the braided
commutator $[\cdot,\cdot]_\qsub\colon V^*\times V \to \CC W$ has zero
kernels in $V^*$ or $V$.  
\end{definition}

To establish the connection to $\qq$\dash Cherednik algebras, recall  
that if there exists a $\qq$\dash Cherednik algebra 
of a group $\widetilde W\le \GL(V)$ such that $\widetilde W$ contains
$\Gamma_\qsub$, then $V$ is a Yetter\dash Drinfeld module over $W$ by 
Corollary~\ref{cor:yd}. This means that such $\qq$\dash Cherednik
algebra has braided reduction. We have

\begin{proposition}
\label{prop:brred}
Braided Cherednik algebras associated to the matrix $\qq$ are the same
as braided reductions of $\qq$\dash Cherednik algebras. 
\end{proposition}
\begin{proof}
First, assume that $A=S_\qsub(V) \tensor \CC \widetilde W \tensor S_\qtran(V^*)$
is a $\qq$\dash Cherednik algebra where $\widetilde W\supset
\Gamma_\qsub$. 
Then it is easily deduced from 
Proposition~\ref{prop:twist} that $\rd A$, which is the same algebra
as $A$ but
with generators $x_i$, $w$ and 
$$
 \uy_i= \gamma_i^{-1} y_i, 
$$
has triangular
decomposition $S_\qsub(V) \tensor \CC W \tensor S_\qsub(V^*)$  and
fulfils the relations in Definition~\ref{def:bcher}. 
To prove that $\rd A$ is a braided Cherednik algebra, we show that $A$
is a non\dash degenerate quadratic double if and only if the $\qq$\dash
commutator $[\cdot,\cdot]_\qsub$ between $V^*$ and $V$ in $\rd A$ 
has zero kernels in $V^*$, $V$. 

Indeed, we have $[\uy_j,x_i]_\qsub = \gamma_j^{-1} \cdot [y_j,x_i]$ for $x\in
V$. Hence it is enough to show that the kernels of $[\cdot,\cdot]$ and 
$[\cdot,\cdot]_\qsub$ are spanned by $\Gamma_\qsub$\dash homogeneous
elements. But observe that $\Gamma_\qsub$\dash homogeneous
elements in $V^*$ and in $V$ are precisely the simultaneous eigenvectors for the
action of $\Gamma_\qsub$. Furthermore, the kernels of $[\cdot,\cdot]$
and of $[\cdot,\cdot]_\qsub$  in $V^*$ and $V$ are $W$\dash submodules,
therefore $\Gamma_\qsub$\dash submodules and thus spanned by
eigenvectors for the action of $\Gamma_\qsub$, as required.

Second, assume that there is a braided Cherednik algebra of the form
$S_\qsub(V) \tensor \CC W \tensor S_\qsub(V^*)$. Then the group $W$ preserves 
the $\qq$\dash exterior squares $\wedge_\qsub^2 V$ and $\wedge_\qsub^2
V^*$. Hence $W$ preserves $\wedge_\qtran^2 V^* = \tau(\wedge_\qsub^2
V^*)$ and, by Proposition~\ref{prop:autom}, $W$ normalises $\Gamma_\qsub$. 
It follows that $\widetilde W := W\cdot \Gamma_\qsub$ is a
group which preserves $\wedge_\qsub^2 V$ and $\wedge_\qsub^2
V^*$. By Corollary~\ref{cor:yd}, $V$ is a Yetter\dash Drinfeld module with
respect to the action of $\widetilde W$ and the grading by elements of
$\Gamma_\qsub \le\widetilde  W$. Put 
$
               y_i = \gamma_i
%^{-1}
\uy_i
$;
then $x_1,\dots,x_n$, $y_1,\dots,y_n$ and $\widetilde W$ generate 
a quadratic double $A$, as shown in Proposition~\ref{prop:twist},
of the form $S_\qsub(V) \tensor \CC \widetilde W \tensor
S_\qtran(V^*)$. Our braided Cherednik algebra is the braided reduction
of $A$. Moreover, by what we have already proved, $A$ is a non\dash
degenerate quadratric double, i.e. the commutator $[\cdot,\cdot]$ 
between $V^*$ and $V$ has zero kernels, because this is true for
$[\cdot,\cdot]_\qsub$.  
\end{proof}

The Proposition and its proof imply a $W$\dash equivariance condition
for the braided commutator:

\begin{corollary}
\label{cor:brequiv}
Let $\rd{\mathcal H}(W)= S_\qsub(V)\tensor \CC W \tensor S_\qsub(V^*)$ 
be a braided Cherednik algebra of a group $W\subset \GL(V)$,
associated to a matrix $\qq$. Then the braided commutator 
$[\cdot,\cdot]_\qsub\colon V^*\tensor V\to \CC W$ is $W$\dash
equivariant and $\Gamma_\qsub$\dash equivariant. (The action of both
$W$ and $\Gamma_\qsub$ on $V^*\tensor V$ is standard diagonal, and on
$\CC W$ is by conjugation inside $\GL(V)$.)
\end{corollary}
\begin{proof}
As in the proof of Proposition~\ref{prop:brred}, put $\widetilde
W=W\cdot \Gamma_\qsub$ and view $\rd{\mathcal H}(W)$ as the braided
reduction of $\mathcal H(\widetilde W) = S_\qsub(V)\tensor \CC W
\tensor S_\qtran(V^*)$. To compute the braided commutator of $\uy\in
V^*$ and $x\in V$ in $\rd{\mathcal H}(W)$, we assume $\uy$ to be
$\Gamma_\qsub$\dash homogeneous, put $\uy=|y|y$ and write 
$
   [\uy, x]_\qsub = [ |y|y,x]_\qsub =|y| [y,x]
$
precisely as in Proposition~\ref{prop:twist}. Now for any $w\in
\widetilde W$ we have
$$
  [ w(\uy),w(x)]_\qsub = |w(y)| [w(y),w(x)]
  = w|y|w^{-1}\cdot w[y,x]w^{-1} = w|y|[y,x]w^{-1} = w[\uy,x]_\qsub w^{-1} 
$$
because $V^*$ is a Yetter\dash Drinfeld module for $\widetilde W$ and
the commutator $[\cdot,\cdot]$ in $\mathcal H(\widetilde W)$ is
$\widetilde W$\dash invariant. Extending  to arbitrary
$\uy$ by linearity, we obtain $\widetilde W=W\cdot \Gamma_\qsub$\dash 
equivariance of $[\cdot,\cdot]_\qsub$.
\end{proof}

Note that the group $W$ may not be stable under conjugation by
$\Gamma_\qsub$, but the braided commutator must still be
$\Gamma_\qsub$\dash equivariant.

\subsection{Negative braided Cherednik algebras}
\label{subsect:sigma}

Clearly, if $q_{ij}=1$ for all $i,j$, braided Cherednik algebras
associated to $\qq$ are ordinary rational Cherednik algebras. We will
now construct a family of braided Cherednik algebras of finite groups with 
$$
         q_{ij}=-1 \quad \text{for all}\ i,j=1,\dots,n,\ i\ne j.
$$
The matrix with such entries was denoted $\negone$ in the Introduction.
Recall, also from the Introduction, 
\begin{itemize}
\item[---] the elements  $\sigma_{ij}^{(\varepsilon)}$ of order $4$ in
  $\GL(V)$, 
  defined for indices $i\ne j$ and for $\varepsilon\in \CC^\times$;
\item[---] the finite group 
$
     W_{\mathcal C, \mathcal C'} = 
     \langle\{\sigma_{ij}^{(\varepsilon)} \mid \varepsilon\in \mathcal
     C\} 
     \, \cup \, 
     \{ \refl_i^{(\varepsilon')} \mid \varepsilon'\in \mathcal C'\}
     \rangle \le \GL(V)         
$,
where $\mathcal C'\subset \mathcal C$ are finite subgroups of
$\CC^\times$ such that $|\mathcal C|$ is even.
\end{itemize}
We will write $W_{\mathcal C, \mathcal C'}(n)$ to emphasise that
there is a separate group $W_{\mathcal C, \mathcal C'}$ in each rank
$n=\dim V$.  
We note that $ W_{\mathcal C, \mathcal C'}$ is an irreducible linear
group, i.e., it irreducibly acts on $V$, and keep in mind that 
$ W_{\mathcal C, \mathcal C'}(n)$ is one of the groups $G(m,p,n)$ or
$G(m, p,n)_+$ with such $m$ and $p$ as described in the Introduction. 
\begin{definition}
\label{def:negbraided}
Fix  a scalar function $c\colon \mathcal C' \to \CC$. 
The negative braided Cherednik algebra $\rd{\mathcal
  H}_{c}(W_{\mathcal C, \mathcal C'})$ is the algebra generated
  by $V=\mathrm{span}(x_1,\dots,x_n)$, $W_{\mathcal C, \mathcal C'}$ and $V^*=\mathrm{span}(\uy_1,\dots,\uy_n)$ subject to relations
\begin{itemize}
\item[$(i)$] $x_ix_j+x_jx_i=\uy_i\uy_j+\uy_j\uy_i=0$ for all $i\ne j$; 
\item[$(ii)$] $wx_i w^{-1} = w(x_i)$, $w\uy_i w^{-1} = w(\uy_i)$;
\item[$(iii)$] $\uy_j x_i+x_i\uy_j=c_1\sum\limits_{\varepsilon\in
  {\mathcal C}}
  \varepsilon \sigma_{ij}^{(\varepsilon)}$
  for all $i\ne j$, 
\newline 
$
\uy_i x_i-x_i\uy_i=1 + c_1\sum\limits_{j\ne i,\ \varepsilon\in {\mathcal C}}
 \sigma_{ij}^{(\varepsilon)}
+\sum\limits_{\varepsilon'\in {\mathcal C}'\setminus\{1\}}
c_{\varepsilon'}t_i^{(\varepsilon')}$,
\end{itemize}
where $i,j=1,\dots,n$ and $w\in W_{\mathcal C, \mathcal C'}$. 
\end{definition}

\begin{proposition}
\label{prop:isbca}
$\rd{\mathcal H}_c(W_{\mathcal C, \mathcal C'})$ is a braided Cherednik
  algebra.
\end{proposition}
\begin{proof}
The matrix  $\qq$ is given by $\qq=\negone$. 
Note that all indices form a single negative
block with respect to $\negone$.

Identify the group $\GL(V)$ with $\GL_n(\CC)$ via
the basis $\{x_i\}$ of $V$.
Let $m=|\mathcal C|$. 
Take $\widetilde W$ to be the complex reflection group
$G(m,1,n)$ of matrices in $\GL_n(\CC)$  with precisely $n$ non\dash
zero entries, all of which are $m$th roots of unity. 
Note that
$$
s_{ij}^{(\varepsilon)}:=(ij)\refl_i^{(\varepsilon)}\refl_j^{(\varepsilon^{-1})}
$$
and $\refl_i^{(\varepsilon)}$, $\varepsilon\in \mathcal C$,  
are complex reflections in $\widetilde W$. 
Let $\mathcal H(\widetilde W)$ be the algebra generated by 
$V$, $\widetilde W$ and $V^*$ subject to the relations in
Proposition~\ref{prop:qcomm} (with $q_{ij}=-1$ for all $i\ne j$ !) and
the commutation relation 
\begin{align*}
\label{eq:x_iy_j commutation}
\tag{$*$}
y_j x_i - x_i y_j 
&= (-\id) \cdot c_1 \sum_{\varepsilon\in \mathcal C} -\varepsilon  
  s_{ij}^{(\varepsilon)} \quad\text{if}\ i\ne j,
\\
\label{eq:x_iy_i commutation}
\tag{$**$}
y_i x_i - x_i y_i 
&= (-\id) \cdot\Bigl( \refl_i^{(-1)} + c_{-1}\cdot 1 
- c_1 \sum_{j\ne i,\ \varepsilon\in \mathcal C} s_{ij}^{(\varepsilon)}   
+ \sum_{\varepsilon'\in \mathcal C'\setminus\{\pm 1\}} 
c_{\varepsilon'} \refl_i^{(-\varepsilon')}\Bigr),
\end{align*}
where $(-\id)$ is the negative identity matrix in $\GL_n(\CC)$ 
(it is an element of $\widetilde W$ since $m$ is even).
The coefficient $c_{-1}$ is assumed to be zero if $-1\not\in\mathcal C'$. 
To observe that $\mathcal H(\widetilde W)$ is a $\negone$\dash
Cherednik algebra,  rewrite the commutation relation as 
\begin{align*}
yx-xy = &(-\id)\Bigl( c_{-1}\langle x,y\rangle + \frac{1}{2}\sum_i \langle x,2y_i\rangle \langle
y,x_i\rangle \refl_i^{(-1)} 
\\
&+
c_1 \sum_{i\ne j}\langle x,y_i-\varepsilon^{-1}y_j \rangle \langle x_i
-\varepsilon x_j, y\rangle s_{ij}^{(\varepsilon)}
+\sum_{i, \ \varepsilon'\in\mathcal C'\setminus\{\pm 1\}}
\frac{c_{\varepsilon'}}{1+\varepsilon'}
\langle x, (1+\varepsilon')y_i\rangle \langle x_i, y\rangle
\refl_i^{(-\varepsilon')} 
\Bigr).  
\end{align*}
This is the same as the commutator in Theorem~\ref{thm:structure}: 
given that there is only one block $B=\{1,\dots,n\}$ of indices which
is negative, one has $\gamma_\Bsub=(-\id)$. This commutator is non\dash
degenerate because of the coefficient in front of $\refl_i^{(-1)}$,
and is $\widetilde W$\dash equivariant, since 
$\refl_i^{(\varepsilon)}$ and $\refl_j^{(\delta)}$
are not conjugate in $\GL(V)$ if $\varepsilon\ne \delta$, and 
$s_{ij}^{(\varepsilon)}$ is
never conjugate to  $\refl_k^{(\varepsilon')}$ in $G(m,p,n)$ (see
\cite[Section 3]{DO}). 

Hence $\mathcal H(\widetilde W)$ is a $\negone$\dash Cherednik algebra by
Corollary~\ref{cor:structure}. It remains to observe that
$\rd{\mathcal H}_{c}(W_{\mathcal C, \mathcal C'})$ is a braided
reduction of $\mathcal H(\widetilde W)$. Indeed, let $(-\id)\beta_{ij}$ be
the commutator of $y_j$ and $x_i$ in $\mathcal H(\widetilde W)$, 
defined above.  
By Proposition~\ref{prop:twist}, the braided commutator 
$[\uy_j,x_i]_\minusone$ in the braided reduction of $\mathcal H(\widetilde
W)$ is equal to
$|\uy_j|(-\id)\beta_{ij}=\gamma_j^{-1}(-\id)\beta_{ij}$. 
Note that 
$\gamma_j^{-1}=\gamma_j$ acts on $x_i$, $i\ne j$, by $-1$ and on $x_j$
by $1$. Thus $\gamma_j^{-1}\cdot(-\id) = \refl_j^{(-1)}$. 
We are left to note that $\refl_j^{(-1)}s_{ij}^{(\varepsilon)} =
\sigma_{ij}^{(-\varepsilon)}$ and  $\refl_i^{(-1)}\refl_i^{(\varepsilon')}
=\refl_i^{(-\varepsilon')}$, therefore $ \gamma_j^{-1}(-\id)\beta_{ij}$ is
precisely the braided commutator of $\uy_j$ and $x_i$ in
Definition~\ref{def:negbraided}.   
\end{proof}

Using the notation from
the proof of Proposition~\ref{prop:isbca}, we can make another ``change of
variables'' in the $\negone$\dash 
Cherednik algebra $\mathcal H(\widetilde W)$ in the case when $-\id\in \widetilde W$. 
Namely,  $V$ becomes a Yetter\dash 
Drinfeld module for $\widetilde W$ via the grading $|v|=-\id$ for all
$v\in V$. By  Proposition~\ref{prop:twist}, the elements
$z_i = (-\id)\cdot y_i \in  \mathcal H(\widetilde W)$,  
together with the $x_i$ and $w\in W$, generate an algebra $\overline
H(\widetilde W)$ with relations
\begin{itemize}
\item[$(i)$] $x_ix_j+x_jx_i=z_iz_j+z_j z_i=0$ for all $i\ne j$;
\item[$(ii)$] $wx_iw^{-1}=w(x_i), w z_iw^{-1}=w(z_i)$ for all $w\in W_{{\mathcal C}, {\mathcal C}'}$, $i=1,\ldots,n$;
\item[$(iii)$] 
$z_jx_i+x_iz_j=c_1\sum\limits_{\varepsilon\in {\mathcal C}}
 - \varepsilon s_{ij}^{(\varepsilon)}$ for all $i\ne j$, and 
\item[]    
                $z_ix_i+x_iz_i= \refl_i^{(-1)} +
                c_{-1}\cdot 1  - c_1 \sum\limits_{j\ne i,\ \varepsilon\in
                \mathcal C} s_{ij}^{(\varepsilon)}    
+ \sum\limits_{\varepsilon'\in \mathcal C'\setminus\{\pm 1\}} 
c_{\varepsilon'} \refl_i^{(-\varepsilon')}$
for $i=1,\ldots,n$,
\end{itemize}      
obtained directly from the relations (\ref{eq:x_iy_j commutation}) and (\ref{eq:x_iy_i commutation}) in 
 the proof of Proposition~\ref{prop:isbca}.
We thus obtain
\begin{corollary}
\label{cor:anticommutator}
The algebra with the above presentation $(i)$--$(iii)$ has  
triangular decomposition $S_\minusone(V)\tensor \CC \widetilde W
\tensor S_{\minusone}(V^*)$.
\end{corollary}
% and the cross\dash commutator relation
%$$
%      z_j x_i + x_i z_j = \beta_{ij}  \qquad\forall i,j.
%$$
%Now, looking at the explicit form of $\beta_{ij}$ in the proof of the
%Proposition,  we conclude that there is a rational Cherednik algebra
%$H(\widetilde W)$, of
%the form $S(V)\tensor \CC \widetilde W \tensor S(V^*)$, given by
%relations 
%$[x_j,x_i]=[z_j,z_i]=0$, $[z_j,x_i]=(-1)^{\delta_{ij}} \beta_{ij}$. 
%The braided reduction has allowed us to easily discover a
%``twist'' of $H(\widetilde W)$ where all commutators in the defining
%relations are replaced by anticommutators. This
%construction works for rational Cherednik algebras of all  
%imprimitive complex reflection groups $G(m,p,n)$. 

\begin{remark}[The degenerate version]
\label{rem:negH0}
We introduce 
the ``degenerate'' negative braided Cherednik algebra $\rd{\mathcal
  H}_{0,c}(W_{\mathcal C,\mathcal C'})$ 
by omitting $1$ from the commutator $\uy_i x_i - x_i \uy_i$ in
Definition~\ref{def:negbraided}:
$$
   \uy_i x_i - x_i \uy_i = 
 c_1\sum\limits_{j\ne i,\ \varepsilon\in {\mathcal C}}
 \sigma_{ij}^{(\varepsilon)}
+\sum\limits_{\varepsilon'\in {\mathcal C}'\setminus\{1\}}
c_{\varepsilon'}t_i^{(\varepsilon')}.
$$
This is a braided Cherednik algebra, provided that the function $c$ is
not identically zero. The proof is the same as for 
 $\rd{\mathcal H}_{c}(W_{\mathcal C,\mathcal C'})$.
\end{remark}
\begin{remark}[The rank $2$ case]
\label{rem:rank2}
It turns out that when $\dim V=2$, the definition of $\rd{\mathcal
  H}_c(W_{\mathcal C,\mathcal C'}(2))$ and $\rd{\mathcal H}_c(W_{\mathcal C,\mathcal C'}(2))$ can be modified to allow one extra
  degree of freedom in choosing the parameter $c$. 
We hereby modify Definition~\ref{def:negbraided} to say that if $\dim
  V=2$, the algebra depends on $|\mathcal C'|+1$ parameters 
$c_1$, $c_1'$, $c_{\varepsilon'}$ ($\varepsilon'\in \mathcal C'\setminus
  \{1\}$), and the commutation relations in $\rd{\mathcal H}_c(W_{\mathcal C,\mathcal C'}(2))$  will
  be
\begin{itemize}
\item[$(iii)$]
$\uy_j x_i + x_i \uy_j = c_1 \sum\limits_{\varepsilon\in\mathcal C^2} 
\varepsilon \sigma_{ij}^{(\varepsilon)} + c_1 \sum\limits_{\varepsilon\in\mathcal C\setminus\mathcal C^2} 
\varepsilon \sigma_{ij}^{(\varepsilon)}$ when $\{i,j\}=\{1,2\}$,
\newline 
$
\uy_i x_i-x_i\uy_i=1 + c_1\sum\limits_{j\ne i,\ \varepsilon\in {\mathcal C}^2}
 \sigma_{ij}^{(\varepsilon)}
+ c_1'\sum\limits_{j\ne i,\ \varepsilon\in {\mathcal C}\setminus{\mathcal C}^2}
 \sigma_{ij}^{(\varepsilon)}
+\sum\limits_{\varepsilon'\in {\mathcal C}'\setminus\{1\}}
c_{\varepsilon'}t_i^{(\varepsilon')}$.
\end{itemize}
Here ${\mathcal C}^2$ denotes the set of squares of elements of
${\mathcal C}$ (the only subgroup of index $2$ in ${\mathcal C}$). 
The proof that $\rd{\mathcal H}_c(W_{\mathcal C,\mathcal C'}(2))$ is a
braided Cherednik algebra is the same as in
Proposition~\ref{prop:isbca}, but 
taking into account that in  the complex reflection group $G(m,1,2)$
the complex reflections 
$s_{12}^{(\varepsilon)}$ and $s_{12}^{(\varepsilon')}$ are conjugate if and
only if $\varepsilon'=\varepsilon \delta^2$ for some $\delta\in\mathcal C^2$.
\end{remark}

It turns out that the algebras $\rd{\mathcal H}_{c}(W_{\mathcal
  C,\mathcal C'})$ and $\rd{\mathcal H}_{0,c}(W_{\mathcal C,\mathcal C'})$
exhaust all possible ``negative braided'' Cherednik algebra structures
  over the group $W_{\mathcal C,\mathcal C'}$:

\begin{proposition}
\label{prop:unique}
Any braided Cherednik algebra $\rd{\mathcal H}$ of the form
$S_\minusone(V)\tensor \CC W_{\mathcal C,\mathcal C'} \tensor
S_\minusone(V^*)$
is isomorphic to $\rd{\mathcal H}_{c}(W_{\mathcal
  C,\mathcal C'})$ or to $\rd{\mathcal H}_{0,c}(W_{\mathcal C,\mathcal
  C'})$ for some choice of the parameter $c$. 
\end{proposition}
\begin{proof}
Let $\widetilde W = W_{\mathcal C,\mathcal C'}\Gamma_\minusone$, and
consider a $\minusone$\dash Cherednik algebra   $\mathcal H(\widetilde W)$
such that $\rd{\mathcal H}$   is its braided reduction (as in the
proof of Proposition~\ref{prop:brred}). The braided commutator 
$[\uy_j,x_i]_\minusone$ in $\rd{\mathcal H}$  rewrites 
as $\gamma_j^{-1}[y_j,x_i]$, where $\gamma_j^{-1}=\gamma_j=t_j^{(-1)}$ 
and $[y_j,x_i]$ is the commutator in
$\mathcal H(\widetilde W)$, necessarily given by 
$[y_j,x_i]=(-\id)(\mathrm{scalar}+\sum_s c_s \langle
x_i,\alpha_s^\vee\rangle \langle \alpha_s,y_i\rangle s)$. Here $s$
runs over some complex reflections in the group $\widetilde W$, and
$c_s$ are some scalars. Now
observe that $\widetilde W$ is contained in the complex reflection
group $G(m,1,n)$ where $m=|\mathcal C|$. We know what are the complex
reflections in $G(m,1,n)$; it follows that for $i\ne j$, the
only possible complex reflections appearing in the commutator $[y_j,x_i]$ 
are of the form $s_{ij}^{(\varepsilon)}$, and if $i=j$,
then they can be of the form 
 $t_i^{(\eta)}$ or $s_{ik}^{(\varepsilon)}$ 
for some $k\ne i$. We do not know what are the linear conditions on the scalars
$c_s$, because this depends on how the complex reflections split into
conjugacy classes in $\widetilde W$;
but we certainly know that the coefficients of the same complex reflection
$s_{ij}^{(\varepsilon)}$ in $[y_j,x_i]$ and in
$[y_i,x_i]$ differ by the factor of $\varepsilon$.

All this is sufficient to determine that the cross\dash commutation
relations in $\rd{\mathcal H}$ must be of the form 
\begin{itemize}
\item[]
$\uy_j x_i + x_i \uy_j = (\uy_j,x_i) + \sum_{\varepsilon\in\mathcal C} 
\varepsilon a(i,j,\varepsilon)\sigma_{ij}^{(\varepsilon)}$,
\item[]
$\uy_i x_i - x_i \uy_i = (\uy_i,x_i) + 
\sum_{j\ne i, \ \varepsilon\in\mathcal C} 
a(i,j,\varepsilon)\sigma_{ij}^{(\varepsilon)} 
+ \sum_{\varepsilon'\in\mathcal C'\setminus\{1\}} b(i,\varepsilon')
t_i^{(\varepsilon')}$
\end{itemize}
for some bilinear form $(\cdot,\cdot)\colon V^*\tensor V \to \CC$ and
some coefficients $a(i,j,\varepsilon)$ ($i\ne j$), 
$b(i,\varepsilon)$. 
Now we are going to use the $W_{\mathcal C,\mathcal C'}$\dash
equivariance of the braided commutator (Corollary~\ref{cor:brequiv}). 
The form $(\cdot,\cdot)$ must be $W_{\mathcal C,\mathcal C'}$\dash
invariant, and as $W_{\mathcal C,\mathcal C'}$ is an irreducible
group, $(\cdot,\cdot)=\lambda\langle\cdot,\cdot\rangle$ is proportional
to the  evaluation pairing. 

Equivariance of the second commutator formula
with respect to $\sigma_{i1}^{(1)}$ implies that 
$a(i,j,\varepsilon)=a(1,j,\varepsilon)$ and $b(i,\varepsilon')=b(1,\varepsilon')$,
and then equivariance under
$\sigma^{(1)}_{2j}$ implies that $a(1,j,\varepsilon)=a(1,2,\varepsilon)$.
Finally, equivariance under
$\sigma_{31}^{(1)}\sigma_{13}^{(\varepsilon)}=t_1^{(\varepsilon^{-1})} 
t_3^{(\varepsilon)}$ implies the equation 
$
    a(1,2,\varepsilon) = a(1,2,1) 
$,  
because $t_1^{(\varepsilon^{-1})} t_3^{(\varepsilon)} \sigma_{12}^{(1)} 
[t_1^{(\varepsilon^{-1})}
  t_3^{(\varepsilon)}]^{-1}=\sigma_{12}^{(\varepsilon)}$. 
The same result can be obtained by using equivariance under  
$t_1^{(\sqrt{\varepsilon^{-1}})}$. Thus, 
$a(i,j,\varepsilon)$ ($\varepsilon\in\mathcal C$) are all equal to some
constant $c_1$, and  $b(i,\varepsilon')=c_{\varepsilon'}$
($\varepsilon'\in\mathcal C'$). One concludes that 
$\rd{\mathcal H}\cong \rd{\mathcal H}_{0,c}(W_{\mathcal C,\mathcal
  C'})$ 
if $\lambda=0$, or 
$\rd{\mathcal H}\cong \rd{\mathcal H}_{c}(W_{\mathcal C,\mathcal C'})$
if $\lambda\ne0$, where $c$ 
is the function $\varepsilon'\mapsto c_{\varepsilon'}$ on $\mathcal C'$. 

The above argument only fails if the group 
$W_{\mathcal C,\mathcal C'}$ does not contain  $t_1^{(\varepsilon^{-1})} 
t_3^{(\varepsilon)}$ and $\mathcal C'$ does not contain 
$\sqrt{\varepsilon}$, for $\varepsilon\in\mathcal C$. This happens precisely
when $\dim V=2$ (the rank $2$ case). In this case, one may use
equivariance of the braided commutator under $\sigma_{12}^{(\delta)}$,
$\delta\in\mathcal C$, to establish
$a(1,2,\varepsilon)=a(1,2,\varepsilon^{-1}\delta^2)$ by observing that 
$(\sigma_{12}^{(\delta)})^{-1}\sigma_{12}^{(\varepsilon)}\sigma_{12}^{(\delta)}=\sigma_{12}^{(\varepsilon^{-1}\delta)}$. 
In this case, the algebra will depend not on $|\mathcal C'|$ but on
$|\mathcal C'|+1$ scalar parameters, as described in
Remark~\ref{rem:rank2}. 
\end{proof}

\begin{example}[Braided Cherednik algebra of type $B_n^+$]
\label{ex:bn+}
The smallest possible example of a nontrivial group $W_{\mathcal C,
  \mathcal C'}$ in rank $n$ corresponds to  $|\mathcal C|=2$ and
  $|\mathcal C'|=1$.  
The group $G(2,1,n)$ is the Coxeter group of type $B_n$, and 
$W_{\{\pm 1\},\{1\}}$ is the group of even elements in $B_n$. Denote
this group by $B_n^+$. It is generated by $\sigma_{ij}$,
$i,j=1,\dots,n$, $i\ne j$, so that $\sigma_{ij}^{(1)}=\sigma_{ij}$ and
$\sigma^{(-1)}_{ij}=(\sigma_{ij})^{-1}=\sigma_{ji}$. 

The following is the list of relations in the negative braided
Cherednik algebra of type $B_n^+$:
\begin{itemize}
\item
$x_i x_j +x_j x_i = \uy_i \uy_j + \uy_j \uy_i = 0$ for $i\ne j$;
\item
$\sigma_{ij} x_i = x_j \sigma_{ij}$,\quad 
$\sigma_{ij} x_j = -x_i \sigma_{ij}$,\quad
$\sigma_{ij} x_k = x_k \sigma_{ij}$ for $k\ne i,j$, and same with
  $\uy_i$ in lieu of $x_i$;
\item
$\uy_j x_i + x_i \uy_j = c (\sigma_{ij}-\sigma_{ji})$ for $i\ne j$;
\item
$\uy_i x_i - x_i \uy_i = 1+c \sum\limits_{j\ne i} (\sigma_{ij}+\sigma_{ji})$.
\end{itemize}
\end{example}

\section{Classification of braided Cherednik algebras}
\label{sect:Classification}

In this Section, we classify braided Cherednik algebras of finite
groups (under a
natural minimality assumption on the group $W$). We do this
 by showing that they
are braided tensor products of rational Cherednik algebras of
irreducible complex reflection groups and negative braided Cherednik
algebras of groups $G(m,p,n)$ and $G(m,p,n)_+$, introduced in the
previous section. 

\subsection{Braided tensor  product of algebras}

For $k=1,\dots,m$, let $\rd{\mathcal H}_k$ be a braided Cherednik algebra of a
finite group $W_k\subset \GL(V_k)$, associated to a matrix $\qq_k$ of size
$n_k \times n_k$ where $n_k = \dim V_k$. 
We would like to turn the 
vector space $\rd{\mathcal H}_1\tensor \dots \tensor \rd{\mathcal
  H}_m$ into a braided Cherednik algebra associated to a matrix $\qq$
of size $n=\sum_k n_k$, with submatrices $\qq_k$ along the main
diagonal. However, the standard tensor product
$A\tensor B$ of
algebras where $a\in A$ and $b\in B$ commute, is not general enough
because it would only give matrix $\qq$ with all entries outside the
submatrices  $\qq_k$ equal to $1$. 

It turns out that the appropriate tensor multiplication here is the
\emph{braided tensor product} of algebras, well known in the theory of
braided monoidal categories; see \cite{Mcompanion}. Let us recall this
notion without going into too much detail.
Let $\mathbf C$ be a braided tensor category, i.e., for each pair
$X$, $Y$ of 
objects there is a braiding $\Psi_{X,Y}\colon X\tensor Y \to Y
\tensor X$ which is a morphism in $\mathbf C$; these morphisms
satisfy axioms of the categorical braiding. An algebra in $\mathbf C$
is an object $A$ of $\mathbf C$ equipped with associative
multiplication $m_A\colon A\tensor A \to A$ and the unit map $1_A
\colon \mathbb{I}\to A$ that are morphisms in $\mathbf C$, where
$\mathbb I$ is the unit object in the category. 
The braided tensor product of algebras $A$ and $B$ in $\mathbf C$ is 
\begin{gather*}
A \utensor B := A\tensor B\ \text{as an object of}\ \mathbf C\,;
\\
m_{A\utensor B} = (m_A\tensor m_B)(\id_A \tensor \Psi_{B,A} \tensor
\id_B) \colon A\tensor B \tensor A \tensor B \to A \tensor B,
\\
1_{A\utensor B} = 1_A \tensor 1_B.
\end{gather*}
The categorical braiding axioms ensure that $m_{A\utensor B}$ is an
associative multiplication. 

\subsection{The braided tensor category $\mathcal M_{\Gamma,\mathcal R}$}

The category $\YD{\Gamma}$ of
Yetter\dash Drinfeld modules (as introduced in Section~\ref{sect:qd})
over a group $\Gamma$ is a braided category, with braiding 
$$
   X,Y\in\Ob\ \YD{\Gamma}  \qquad \mapsto \qquad 
   \Psi_{X,Y} \colon X\tensor Y \to Y \tensor X, \quad
\Psi_{X,Y}(x\tensor y)=|x|(y)\tensor x.             
$$
Our main example of a braided category will, however, be slightly
different. Let $\Gamma$ be an Abelian group. 
Fix a map $\mathcal R \colon \Gamma\times \Gamma \to
\CC^\times$ which is a \emph{unitary bicharacter}, i.e., 
$\mathcal R(\gamma,\cdot)$ is a homomorphism from $\Gamma$ to
$\CC^\times$ for fixed $\gamma\in \Gamma$, and $\mathcal
R(\gamma,\gamma')=\mathcal R(\gamma',\gamma)^{-1}$. 
Assume that $\mathcal R$ is \emph{non\dash degenerate}, that is,
$\gamma\mapsto \mathcal R(\gamma,\cdot)$ is an embedding 
of $\Gamma$ in the group $\widehat \Gamma$ of its multiplicative characters. 
Denote by $\widehat \Gamma_{\mathcal R}$ the image of this embedding.
Elements of $\widehat \Gamma_{\mathcal R}$ are viewed as $1$\dash
dimensional $\Gamma$\dash modules. 

\begin{definition}[The category $\mathcal M_{\Gamma,\mathcal R}$]
\label{def:braidcat}
Define the category $\mathcal M_{\Gamma,\mathcal R}$ as a full
subcategory of $\Gamma$\dash modules consisting of objects 
isomorphic to direct sums of modules in
$\widehat \Gamma_{\mathcal R}$. 
Each module $X$ in $\mathcal M_{\Gamma,\mathcal R}$ is  $\Gamma$\dash
graded by 
$$
      x\in X,\quad g(x)=\mathcal R(\gamma,\cdot)x \quad 
      \forall g\in \Gamma\qquad \Rightarrow \qquad
      |x| = \gamma\,.
$$ 
It is clear that such grading makes $X$ a Yetter\dash Drinfeld
module so that $\mathcal M_{\Gamma,\mathcal R}$ is a full subcategory
of $\YD{\Gamma}$ and defines the braiding
$\Psi_{X,Y}$ between $X$ and $Y$ in $\mathcal M_{\Gamma,\mathcal R}$. 
\end{definition}

In what follows, $\Gamma$ will be an Abelian group with fixed 
unitary non\dash degenerate bicharacter $\mathcal R$ on $\Gamma$. 
For $X\in\mathrm{Ob}\ M_{\Gamma,\mathcal R}$, denote
$$
     |X| = \{\gamma\in \Gamma : \text{there exists}\ x\in X,\ x\ne 0,\
     |x|=\gamma\}. 
$$
\subsection{Braided tensor product of braided Cherednik algebras}

We will now observe that if $\Gamma$ acts on a
braided Cherednik algebra  $\rd {\mathcal H}(W)$ in a certain
standard way,   $\rd {\mathcal H}(W)$ is guaranteed to be a
$\Gamma$\dash module algebra. 

\begin{definition}
\label{def:gammastr}
Let $\rd {\mathcal H}(W) \cong S_\qsub(V) \tensor \CC W \tensor
S_\qsub(V^*)$ be a braided Cherednik algebra, where $V$, as usual, is
spanned by $x_1,\dots,x_n$. 
A \emph{$\Gamma$\dash structure} on $\rd {\mathcal H}(W)$
is a representation $\rho\colon \Gamma\to \GL(V)$ such that:
\begin{itemize} 
\item[---] $V$ becomes an object of $\mathcal M_{\Gamma,\mathcal R}$,
  and $x_i\in V$ are simultaneous eigenvectors for $\rho(\Gamma)$;
\item[---]
$g^{-1}w^{-1}gw \in W \cap \rho(\Gamma)$ for all $g\in
\rho(\Gamma)$, $w\in W$. 
In particular, $(\gamma,w)  \mapsto \rho(\gamma) w \rho(\gamma)^{-1}$
is a $\Gamma$\dash action on $W$;
\item[---]
the braided commutator $[\cdot,\cdot]_\qsub\colon V^*\tensor V \to \CC W$ is
equivariant with respect to the $\Gamma$\dash action. 
\end{itemize}
\end{definition}
\begin{lemma}
\label{lem:degrees}
Let $\rd {\mathcal H}(W)$ be a  braided Cherednik algebra  
with $\Gamma$\dash structure. Then:

$(a)$ $W$ acts on the set $|V|\subset \Gamma$; 

$(b)$ $\CC W$ is a $\Gamma$\dash submodule of  $\rd {\mathcal H}(W)$,
and $|\CC W| = \{\gamma^{-1}\cdot w(\gamma) : \gamma\in |V|\}$;
 
$(c)$ $\rd {\mathcal H}(W)$ is an algebra in the category $\mathcal
M_{\Gamma,\mathcal R}$.
\end{lemma}
\begin{proof}
$(a)$  If $x\in V$ is a simultaneous eigenvector for $\rho(\Gamma)$,
  and $w$ is in $W$,
  then $w(x)$ is a simultaneous eigenvector for $w\rho(\Gamma)w^{-1}$
  corresponding to the same eigencharacter. But
  $w\rho(\Gamma)w^{-1}=\rho(\Gamma)$ by definition of $\Gamma$\dash
  structure. Thus, the action of $W$ permutes $\rho(\Gamma)$\dash
  simultaneous eigenspaces in 
  $V$, which are $\Gamma$\dash homogeneous components of $V$, hence
  $W$ permutes $\Gamma$\dash degrees of elements of $V$. 
  Note that the action of 
  $W$ on the set $|V|$ is such that $w(|x_i|)=|w(x_i)|$. 

$(b)$ Consider the $\Gamma$\dash action on $\End(V)$ given by
  $(\gamma,m) \mapsto \rho(\gamma) m \rho(\gamma)^{-1}$ for $m\in
  \End(V)$. It is then easy to see that the canonical isomorphism
  $\End(V)\cong V\tensor V^*$ 
  is an isomorphism of $\Gamma$\dash modules. Let $\{\uy_i\}$, as usual, be
  the basis of $V^*$ dual to $\{x_i\}$. Then $x_i \tensor \uy_j\in
  V\tensor V^*$ is a simultaneous eigenvector for $\Gamma$ of $\Gamma$\dash degree
  $|x_i ||x_j|^{-1}$. An element $w\in W$ is written as $\sum_{i=1}^n
  w(x_i)\tensor \uy_i\in V\tensor V^*$, and the $\Gamma$\dash degree
  of $  w(x_i)\tensor \uy_i$ is $w(|x_i|)\cdot |x_i|^{-1}$. Thus,
  $\Gamma$\dash degrees that appear in the $\Gamma$\dash submodule of $\CC
  W$  generated by $w$ are of the form $w(|x_i|)\cdot |x_i|^{-1}$, and
  the linear 
  independence of $ w(x_i)\tensor \uy_i$ in the expansion of $w$  
  implies that all these $\Gamma$\dash degrees actually appear in this
  submodule. 

$(c)$ Thus, $\Gamma$ acts on generators of $\rd {\mathcal H}(W)$, and 
we check that this action preserves the relations in
$\rd{\mathcal H}(W)$. 
The relation $wxw^{-1}=w(x)$ where $x\in V$ and  $w\in W$ becomes 
$\rho(\gamma)w\rho(\gamma)^{-1} \cdot \rho(\gamma)(x) \cdot 
\rho(\gamma)w^{-1}\rho(\gamma)^{-1} = (\rho(\gamma)w)(x)$, i.e., 
$w' x' w'^{-1} = w'(x')$ where $w'=\rho(\gamma)w\rho(\gamma)^{-1}\in W$ and
$x'=\rho(\gamma)(x)$. This is also a relation in $\rd{\mathcal
  H}(W)$. The relations $x_i x_j = q_{ij} x_j x_i$ are preserved since
the $x_i$ are simultaneous eigenvectors of $\Gamma$. Similarly, the relations
$w\uy w^{-1}=w(\uy)$ for $\uy\in V^*$, and $\uy_i \uy_j = q_{ij} \uy_j
\uy_i$, are preserved. 
Finally, the braided commutation relations between $\uy_j$
and $x_i$ are preserved because the braided commutator between $V^*$
and $V$ is $\Gamma$\dash equivariant.  
\end{proof}

The $\Gamma$\dash structure 
paves the way for introducing braided tensor product
$\utensor$ of
braided Cherednik algebras.
(Obviously, the usual tensor product
where the two tensorands commute is a particular case of this,
corresponding to the trivial ``$\{1\}$\dash structure'' on any braided
Cherednik algebra.) Let us write down the triangular
decomposition property of the braided tensor product:
\begin{lemma}
\label{lem:brprod}
Let $\rd {\mathcal H}(W)\cong S_\qsub(V) \tensor \CC W \tensor
S_\qsub(V^*)$ and 
$\rd {\mathcal H}(W')\cong S_{\qq'}(V') \tensor \CC W' \tensor
S_{\qq'}(V'^*)$ be braided Cherednik algebras with
$\Gamma$\dash structure, where $V$ is spanned by variables
$x_1,\dots,x_n$ and $V'$ is spanned by variables
$x_{n+1},\dots,x_{n+m}$. 
Let $\mathcal R$ be a non\dash
degenerate unitary bicharacter on $\Gamma$. The braided tensor product
of  $\rd {\mathcal H}(W)$ and 
$\rd {\mathcal H}(W')$ in the category $\mathcal M_{\Gamma,\mathcal R}$ 
has triangular decomposition
$$
\rd {\mathcal H}(W)\utensor \rd {\mathcal H}(W') \cong
S_{\widetilde\qq}(V\oplus V') \tensor (\CC W \utensor \CC W') \tensor
S_{\widetilde \qq}(V^*\oplus V'^*).
$$
The $(n+m)\times (n+m)$ 
matrix $\widetilde \qq=(\widetilde q_{ij})$ is given by 
$$
\widetilde q_{ij}=q_{ij} \ (i,j\le n), \qquad
\widetilde q_{ij}= q'_{ij}\ (n<i,j),   \qquad
\widetilde q_{ij}= \mathcal R(g_i,g_j)\ (\text{otherwise}),
$$
where $g_i=|x_i|$ is the $\Gamma$\dash degree of $x_i$; in particular, 
$\widetilde q_{ij}= \widetilde q_{ij}^{-1}$.
\end{lemma}
\begin{remark}
Warning: $\CC W \utensor \CC W'$ is not necessarily a group algebra!
\end{remark}
\begin{proof}[Proof of the Lemma]
It is clear that we may write $\rd {\mathcal H}(W)\utensor \rd
{\mathcal H}(W') $ as a tensor product 
$$
  (S_\qsub(V)\tensor S_{\qq'}(V)) \tensor (\CC W \tensor \CC W') \tensor
  (S_\qsub(V^*) \tensor S_{\qq'}(V'^*)) =: U^- \tensor U^0 \tensor U^+
$$
of subalgebras. 
The subalgebra $U^-$ is generated by $x_1,\dots,x_{n+m}$ modulo the
relations 
$$
   x_i x_j = q_{ij} x_j x_i \ (i,j\le n), 
\qquad
   x_i x_j = q'_{ij} x_j x_i \ (n<i,j), 
\qquad
   x_i x_j = \mathcal R(g_i,g_j) x_j x_i \ (\text{otherwise}),   
$$
the latter relation being dictated by the braided tensor product. 
Immediately $U^-=S_{\widetilde \qq}(V\oplus V')$ as required. 
In the same way $U^+=S_{\widetilde \qq}(V^*\oplus V'^*)$. 
Moreover, $U^- U^0$ is a  subalgebra. This follows from the
commutation relations
$w x_i = w(x_i) w$, $w' x_k = w'(x_k)w$ where $i\le n<k$, $w\in W$,
$w'\in W'$, and some way (provided by the braided tensor product) 
to re\dash express the product 
$w' x_i$ as an element in $V \tensor \CC W'$; and a way 
to re\dash express $w x_k$ as an element in $V'\tensor \CC W$.   
Similarly, $U^0 U^+$ is a subalgebra. 
\end{proof}

In general, however, $\utensor$ applied to braided Cherednik algebras  
$\rd {\mathcal H}(W)$ and $\rd{\mathcal H}(W')$ will not produce a
braided Cherednik algebra, at least because the associative algebra 
$\CC W \utensor \CC W'$ may not be the group algebra of $W\times
W'$. This generalisation of braided (and in particular, rational)
Cherednik algebras may deserve to be studied elsewhere. 
For the purposes of the present paper, we would like to 
force $\rd {\mathcal H}(W) \utensor \rd{\mathcal H}(W')$
to be a braided Cherednik algebra by some extra condition on the
bicharacter $\mathcal R$ on $\Gamma$. 
Here is the criterion for the braided product of 
two braided Cherednik algebras to be a braided Cherednik algebra of the
direct product of groups:

\begin{proposition}
\label{prop:braidedprod}
In the notation of  Lemma~\ref{lem:degrees} and  Lemma~\ref{lem:brprod}, 
$\rd {\mathcal H}(W)\utensor \rd {\mathcal H}(W')$ is a braided
Cherednik algebra of the group $W\times W'$  acting on the space
$V\oplus V'$, if and only if 
$\mathcal R( w(\gamma), w'(\gamma'))=\mathcal R(\gamma,\gamma')$ for
all $w\in W$, $w'\in W'$, $\gamma\in |V|$, $\gamma'\in |V'|$. 
\end{proposition}
\begin{proof}
Clearly, $\rd {\mathcal H}(W)\utensor \rd {\mathcal H}(W')$  is a
braided Cherednik algebra of the group $W\times W'$  acting on  
$V\oplus V'$, only if the following relations hold in
$\rd {\mathcal H}(W)\utensor \rd {\mathcal H}(W')$:
\begin{itemize}
\item[1.] 
$w w' = w'w$ for $w\in W$ and $w'\in W'$, equivalent to $\CC
W\utensor \CC W'$ being the group algebra of $W\times W'$;
\item[2.] 
$x_k w = w x_k$ ($w\in W$,
    $n<k\le n+m$),
$w' x_i = x_i w'$ ($w'\in W'$,    $1\le k\le n$), 
which are equivalent to the correct smash product relations between
$W\times W'$ and $x_1,\dots,x_{n+m}$; 
\item[3.] same as 2., but with $\uy_i$ instead of $x_i$.
\end{itemize}
Let us observe that conditions 1.--3.\ are not only necessary 
but also sufficient. Indeed, the commutation relation 
in $\rd {\mathcal H}(W)\utensor \rd {\mathcal H}(W')$
between $\uy_j$ and $x_i$ where,
say, $i\le n<j$, is
$$
   \uy_j x_i = \mathcal R(g_j^{-1},g_i) x_i \uy_j = \widetilde q_{ij}
   x_i \uy_j \qquad \Leftrightarrow \qquad [\uy_j,x_i]_{\widetilde \qq}=0
$$
where the matrix $\widetilde \qq$ is given in Lemma~\ref{lem:brprod},
and the same holds for $j\le n<i$. Thus, $[\uy_j, x_i]_{\widetilde
  \qq}\in \CC W \tensor \CC W'$ for all indices $i,j$. 
Besides that, the $\CC W\tensor \CC W'$\dash valued braided commutator 
$[\cdot,\cdot]_{\widetilde \qq}$ on $(V^*\oplus V'^*)\times
(V\oplus V')$ has no left or right kernel, because it coincides with
the non\dash degenerate commutator $[\cdot,\cdot]_\qsub$ (respectively 
$[\cdot,\cdot]_{\qq'}$) when restricted to $V^*\times V$ (respectively 
$V'^*\times V'$) and has zero restriction to $V^*\times V'$ and to
$V'^*\times V$. 
Thus, 1.--3.\ imply that 
$\rd {\mathcal H}(W)\utensor \rd {\mathcal
  H}(W')$ is a braided Cherednik algebra associated to the matrix
$\widetilde \qq$. 

But 1.--3.\ can clearly be rewritten as relations
$$
    w' a = a w' \quad (a\in \rd {\mathcal H}(W),\ w'\in W'),
\qquad
    bw  = wb \quad (b\in \rd {\mathcal H}(W'),\ w\in W).
$$
in the algebra $\rd {\mathcal H}(W)\utensor \rd {\mathcal
  H}(W')$. 
In terms of the bicharacter $\mathcal R$ on $\Gamma$ these are
equivalent to
$\mathcal R(g,h)=1$ for all $g\in|\CC W| $, $h\in|\rd {\mathcal
  H}(W')| $ and for all $g\in|\rd {\mathcal H}(W)|$,  
$h\in |\CC W'|$. 
This is precisely equivalent to the $(W,W')$\dash invariance of 
the restriction of $\mathcal R$ to $|V|\times |V'|$, stated in the
Proposition, because $|\CC W|$ consists of $\gamma^{-1} w(\gamma)$ 
where $w$ runs
over $W$ and $\gamma$ runs over $|V|$, and $|\rd {\mathcal
  H}(W')| $ lies in the subgroup of $\Gamma$ generated by $|V'|$;
similarly for $|\CC W'|$ and $|\rd {\mathcal
  H}(W)| $.
\end{proof}

\subsection{Braided Cherednik algebras: the main structural theorem}
\label{subsect:Braided Cherednik algebras: the main structural theorem}

Braided tensor multiplication is a powerful method of constructing
new associative algebras. Even restricted by the orthogonality condition in
Proposition~\ref{prop:braidedprod}, braided
tensor multiplication is sufficient for obtaining \emph{essentially} any braided
Cherednik algebra of a finite group as a product of algebras from the
following list: 
\begin{itemize}
\item[(1)]
$H_c(W)$ where $W$ is an irreducible complex reflection group
$G(m,p,n)$ or $G_4,\dots,G_{37}$ in the Shephard-Todd classification
  \cite{ST}; 
\item[(2)]
$H_{0,c}(W)$, $W$ same as in (1);
\item[(3)]
 $\rd{\mathcal H}_c(W)$, where $W$ is $G(m,p,n)$ with $m$ even,
 $\frac{m}{p}$ even, $n\ge 2$, or the subgroup $G(m,p,n)_+$ of even
 elements in $G(m,p,n)$ with $m$ even,
 $\frac{m}{2p}$ odd, $n\ge 2$;
\item[(4)]
$\rd{\mathcal H}_{0,c}(W)$, $W$ same as in (3).
\end{itemize}
We gave the definition of algebras $H_c(W)$ and $H_{0,c}(W)$ in
Example~\ref{ex:cherEG}.
These are rational Cherednik
algebras, whereas $\rd{\mathcal H}_c(W)$
(Definition~\ref{def:negbraided}, Remark~\ref{rem:rank2}) 
and  $\rd{\mathcal H}_{0,c}(W)$
(Remark~\ref{rem:negH0}) are negative braided Cherednik algebras. 
In other words, the  $(1)$--$(4)$ is the list of all rational and
negative braided Cherednik algebras of irreducible groups.

``Essentially any braided Cherednik algebra'' refers to $\rd{\mathcal
  H}(W)$ that satisfies the condition in

\begin{definition}
We say that the group $W$ is minimal for a braided Cherednik algebra 
$\rd{\mathcal H}(W) \cong S_\qsub(V) \tensor \CC W \tensor
S_\qsub(V^*)$, if the image of the braided commutator $[\cdot,\cdot
]_\qsub \colon V^* \times V \to \CC W$ does not lie in $\CC W_1$ for
any proper subgroup $W_1 \le W$.
\end{definition}
Every braided Cherednik algebra $\rd{\mathcal H}(W) \cong S_\qsub(V) \tensor \CC W \tensor
S_\qsub(V^*)$ contains 
a subalgebra given as $S_\qsub(V) \tensor \CC W_{\mathrm{min}} \tensor
S_\qsub(V^*)$ where $W_{\mathrm{min}}$ is the subgroup 
 generated by elements of $W$ that appear in
braided commutators between $V^*$ and $V$ with non\dash zero
coefficients. Clearly, all relevant information about a braided
Cherednik algebra $\rd{\mathcal H}(W)$ is contained in this
subalgebra, the braided Cherednik algebra of $ W_{\mathrm{min}}$. This
is the reason why rational Cherednik algebras are often considered
only over complex reflection groups. We apply the same principle to
arbitrary braided Cherednik algebras: 

\begin{theorem}
\label{thm:main}
Any braided Cherednik algebra $\rd{\mathcal H}(W)$ where $W$ is
minimal and finite, decomposes as a braided tensor product of algebras from
$(1)$--$(4)$ above in the category
$\mathcal M_{\Gamma,\mathcal R}$ 
for some  Abelian group $\Gamma$ and a non\dash degenerate
unitary bicharacter $\mathcal R$. 
\end{theorem}
\begin{remark}
The Theorem does not hold for $\qq$\dash Cherednik algebras.  
\end{remark}
\begin{proof}[Proof of Theorem~\ref{thm:main}]
Fix the triangular decomposition $\rd{\mathcal H}(W) \cong S_\qsub (V)
\tensor \CC W \tensor S_\qsub(V^*)$ where $V$ is spanned by the
variables $x_1,\dots,x_n$ and $V^*$ is spanned by
$\uy_1,\dots,\uy_n$. We put $\Gamma=\Gamma_\qsub$; recall
from Section~\ref{sect:q_cher}  that
$\Gamma_\qsub$ is the  subgroup of $\GL(V)$ generated by
$\gamma_1,\dots,\gamma_n$ where $\gamma_i(x_j)=q_{ij} x_j$. 
We let $\mathcal R = \mathcal R_\qsub$ be determined, via the extension
to the whole of $\Gamma_\qsub$ by the bicharacter property, 
by $\mathcal R_\qsub (\gamma_i,\gamma_j)=q_{ij}$. 
It is easy to check that 
$\mathcal R_\qsub$ is a well\dash defined non\dash degenerate unitary
bicharacter. 
In the course of the proof it will become apparent that the natural
representation of $\Gamma_\qsub$ on $V$ is a $\Gamma_\qsub$\dash structure on
$\rd{\mathcal H}(W)$.

\subsection*{Step 1: Factorisation into algebras $\mathcal H_B$
  indexed by blocks $B$}

Recall from Corollary~\ref{cor:block} that the index set splits into
disjoint blocks so that  
for $i\ne j$ one has $q_{ij}=q_{B,C}$ where $B$, $C$ are blocks, $i\in
B$, $j\in C$. Each 
block $B$ is either positive ($q_{B,B}=1$) or negative ($|B|>1$,
$q_{B,B}=-1$). 
Recall from Proposition~\ref{prop:brred} that $\rd{\mathcal H}(W)$ is
a braided reduction of a  
$\qq$\dash Cherednik algebra $\mathcal H(\widetilde W) \cong S_\qsub (V)
\tensor \CC \widetilde W \tensor S_\qtran(V^*)$. Here $\widetilde W$ is
a group containing $W$ (we assume $\widetilde W = W\cdot \Gamma_\qsub$ as
in the proof of Proposition~\ref{prop:brred}), and
the basis of $V^*$ is now given by $y_i =
\gamma_i^{-1} \uy_i$, $i=1,\dots,n$, so that
$$
    C_{ij} := [\uy_j,x_i]_\qsub = \gamma_j^{-1} [y_j,x_i].
$$
By Corollary~\ref{cor:zerocomm}, $y_j$ commutes with $x_i$ unless $i$,
$j$ are in the same block. Equivalently, 
$$
   \uy_j x_i = q_{ij} x_i \uy_j = \mathcal R(\gamma_j^{-1},\gamma_i)
   x_i \uy_j
$$
if $i,j$ are not in the same block. 
If $i,j$ are in the same block (say $B$), Theorem~\ref{thm:structure} implies
$$
   C_{ij} = \gamma_j^{-1} \gamma_\Bsub ( a_{ij} + \sum_s b_{i,j,s} s) 
$$
for certain constants $a_{ij}$ and $b_{i,j,s}$, where $s$ runs over a
set of complex reflections in $\GL(V_B)$ with $V_B=\oplus_{k\in B}\CC
x_k$. The group $\GL(V)$ contains the subgroup 
$\GL(V_{B_1})\times \dots\times \GL(V_{B_N})$ where $B_1,\dots,B_N$
are all blocks of indices, $B_1\cup \dots \cup B_N=\{1,\dots,n\}$. 
This corresponds to the direct sum
decomposition $V=V_{B_1}\oplus \dots \oplus V_{B_N}$.
Note the crucial fact that
$$
 \gamma_j^{-1} \gamma_\Bsub = 1\quad \text{if the block $B$ is positive}, 
 \quad
 \gamma_j^{-1} \gamma_\Bsub = t_j^{(-1)}\quad \text{if $B$ is negative}, 
$$
where $t_j^{(-1)}$ acts by $-1$ on $x_j$ and by $1$ on the rest of the
variables $x_i$. Thus, $t_j^{(-1)}\in \GL(V_B)$, and  
all elements of $W$ that appear in $C_{ij}$ ($i,j\in B$) 
with non\dash zero coefficients are elements of $\GL(V_B)$. By
minimality, $W$ is generated by such elements, and hence has direct
product decomposition 
$$
    W = W_{B_1}\times \dots \times W_{B_N} \quad \subset \quad  
\GL(V_{B_1})\times \dots \times \GL(V_{B_N}), \qquad
    W_{B_k} := W\cap \GL(V_{B_k}). 
$$
It is now clear that $\rd{\mathcal H}(W)$ (as a vector space) has
factorisation 
$$
  \rd{\mathcal H}(W) =  \mathcal H_{B_1} \tensor \dots \tensor
  \mathcal H_{B_N}, 
$$  
where $\mathcal H_{B}$ is the subalgebra of $ \rd{\mathcal H}(W)$
generated by $x_i$, $\uy_i$ ($i\in B$) and $W_B$. 
If $B$ is a positive block, 
the algebra $\mathcal H_{B}$ has triangular decomposition 
$\mathcal H_{B} \cong S(V_B)\tensor \CC W_B \tensor S(V_B^*)$ and is a
rational Cherednik algebra of a complex reflection group $V_B$. 
If $B$ is negative block, $\mathcal H_{B}$ has triangular
decomposition $S_\minusone(V_B)\tensor  \CC W_B \tensor S_\minusone(V_B^*)$.

Let us show that the algebra $\rd{\mathcal H}(W)$ and 
all subalgebras $\mathcal H_{B}$ have $\Gamma_\qsub$\dash structure
given by the action of $\Gamma_\qsub$ on $V$ (and hence on all $V_B$
that are $\Gamma_\qsub$\dash submodules of $V$).
According to Definition~\ref{def:gammastr}, we have to check that
\begin{itemize}
\item[1.] $\gamma w \gamma^{-1}w^{-1} \in W_B \cap \Gamma_\qsub$ 
for $\gamma\in \Gamma_\qsub$, $w\in W_B$;
\item[2.] the braided commutator $[\cdot,\cdot]_\qsub \colon V_B^*
  \times V_B\to \CC W_B$ is $\Gamma_\qsub$\dash equivariant.
\end{itemize}  
In 1., we already know that $W$ normalises $\Gamma_\qsub$ by
Proposition~\ref{prop:autom}, hence it is enough to check that 
$\Gamma_\qsub$ normalises $W_B$ in $\GL(V)$. But this follows from 2.,
because by minimality of $W$,  $W_B$ is generated by elements of
$\GL(V)$ that appear in the braided commutator $[\cdot,\cdot]_\qsub
\colon V_B^*   \times V_B\to \CC W$. Of course, 2.\ is true by
Corollary~\ref{cor:brequiv}. 

Let us now show that 
$\rd{\mathcal H}(W)$ is a tensor product of the  $\mathcal H_{B}$ not
only as a vector space but as an algebra in the category 
$\mathcal M_{\Gamma_\qsub,\mathcal R_\qsub}$. 
Since we already have the $\Gamma_\qsub$\dash structure on $\mathcal
H_B$ and tensor factorisation of $\rd{\mathcal H}(W)$
into the algebras $\mathcal H_B$, it is enough to check that
the commutation relations in $\mathcal H_B\utensor \mathcal H_{B'}$ 
 between $\mathcal H_B$ and $\mathcal H_{B'}$
for blocks $B\ne B'$ hold also in $\rd{\mathcal H}(W)$.
By Lemma~\ref{lem:brprod},   the $x_i$ $\qq$\dash commute in 
$\mathcal H_{B}\utensor\mathcal H_{B'}$, as well as the $\uy_i$; the
same holds in $\rd{\mathcal H}(W)$. Furthermore, the definition of
$\utensor$ prescribes the relations $\uy_j x_i = q_{ij} x_i \uy_j$ in 
$\mathcal H_{B}\utensor\mathcal H_{B'}$ where $i\in B$, $j\in B'$;
we have already shown in this proof that the same holds in
$\rd{\mathcal H}(W)$. Similarly for $i\in B'$, $j\in B$. 
Finally, let us find the relations between $W_B$ and $\mathcal H_{B'}$
and between $\mathcal H_B$ and $W_{B'}$ 
in $\mathcal H_{B}\utensor\mathcal H_{B'}$. 
The group $W_B$ acts on $\{\gamma_i\mid i\in B\}$ by
conjugation, and for two blocks $B\ne B'$ we have 
$$
 \mathcal R_\qsub(w\gamma_i w^{-1},w' \gamma_{i'} w'^{-1}) = \mathcal
 R_\qsub(\gamma_i, \gamma_{i'}) = q_{ii'}
$$ 
for all $w\in W_B$, $w'\in W_{B'}$, $i\in B$, $i'\in B'$ simply
because  $q_{jj'}=q_{ii'}$ for all $j\in B$,
$j'\in B'$ by definition of a block. Therefore, by
Proposition~\ref{prop:braidedprod} $W_B$ commutes with 
$\mathcal H_{B'}$ and $\mathcal H_B$ commutes with $W_{B'}$ in
$\mathcal H_{B}\utensor\mathcal H_{B'}$. But again, the same happens
in $\rd{\mathcal H}(W)$. Hence the braided tensor product 
$\mathcal H_{B_1}\utensor \dots \utensor \mathcal H_{B_N}$ indeed 
coincides with $\rd{\mathcal H}(W)$. 

\subsection*{Step 2: Factorisation of $\mathcal H_B$, $B$ a positive
  block, into rational Cherednik algebras of irreducible groups}

It remains to break up each of the algebras $\mathcal H_B$ into a
braided tensor product, in the category $\mathcal
M_{\Gamma_\qsub,\mathcal R_\qsub}$, of ``elementary'' braided Cherednik
algebras from the list $(1)$--$(4)$.  The more familiar case is that
of a positive block, where $\mathcal H_B$ has triangular
decomposition $S(V_B)\tensor \CC W_B \tensor S(V_B^*)$ with commutation
relation 
$$
   yx-xy = (x,y)_B\cdot 1 + \sum_s c_s \langle x,\alpha_s^\vee\rangle
   \langle \alpha_s,y\rangle s,
$$
where $s$ runs over complex reflections in $W_B$, and $(\cdot,\cdot)_B$
is some $W_B$\dash invariant bilinear form on $V_B\times V_B^*$. 
By a known result
on complex reflection groups, $W_B\subset \GL(V_B)$ is
a direct product $W^1\times \dots\times W^l$ of irreducible complex
reflection groups corresponding to a direct sum decomposition
$V_B=V^1\oplus \dots \oplus V^l$.  Denote by $\pi^k\colon V_B \to V^k$
the projection of $V$ onto its direct summand $V^k$. 
The dual space $V_B^*$ has
the dual direct sum decomposition $V^{1*}\oplus \dots\oplus V^{l*}$
with $V^{k*}=\im \pi^{k*}$. Since $V^k$ are irreducible $W$\dash
submodules of $V$, the $W$\dash invariant pairing $(\cdot,\cdot)_B$
between $V$ and $V^*$ is of the form $\sum_{k=1}^l \lambda_k\langle
\pi^k(\cdot), \cdot\rangle$ for some $\lambda_k\in \CC$.
Moreover, any complex reflection $s\in W$ belongs to one of the $W^k$,
thus $\alpha_s\in V^k$ and $\alpha_s^\vee\in V^{k*}$. 
It follows that $x\in V^k$ and $y\in V^{k'}$ commute in $\mathcal H_B$
for $k\ne k'$, and $\mathcal H_B$ decomposes as the tensor product
$$
\mathcal H_B = H(W^1) \tensor \dots \tensor H(W^l)
$$
of commuting subalgebras.
Here $H(W^k)=S(V^k)\tensor\CC W^k \tensor S(V^{k*})$ with the main
commutation relation 
$$
    yx-xy = \lambda_k \langle x,y\rangle + \sum_{s\in W^k}
 c_s \langle x,\alpha_s^\vee\rangle
   \langle \alpha_s,y\rangle s 
$$
between $y\in V^{k*}$ and $x\in V^k$, 
thus is a rational Cherednik algebra isomorphic to either
$H_{0,c}(W^k)$ or $H_c(W^k)$, depending on whether $\lambda_k$ is zero
or not.  

It remains to note that the standard tensor product of commuting
subalgebras is in this case the same as braided tensor product in the
category  $\mathcal M_{\Gamma_\qsub,\mathcal
  R_\qsub}$. First of all, $\Gamma_\qsub$ acts by scalars on $V_B$ and
$V^*_B$, hence trivially on $W_B$. Thus $V^k$, $\CC W^k$ and $V^{k*}$ are
$\Gamma_\qsub$\dash submodules of $V_B$, $\CC W_B$ and $V_B^*$,
respectively.  We now only need to check that the commutation
relations between $H(W^k)$ and $H(W^{k'})$ inside $\mathcal H_B$
(where these two subalgebras commute) are the
same as in the braided tensor product $H(W^k)\utensor H(W^{k'})$. 
Note the $\Gamma_\qsub $\dash degrees that arise in
the $\Gamma_\qsub $\dash module $ \mathcal H_B$ lie in the subgroup of
$\Gamma_\qsub$ generated by $\{\gamma_i \mid i\in B\}$; therefore, 
the value of $\mathcal R_\qsub$ on any two such degrees is $1$,
because $\mathcal R_\qsub(\gamma_i,\gamma_j)=q_{ij}=1$ for any $i,j$ in the
positive block $B$. Thus, $H(W^k)$ and $H(W^{k'})$ indeed commute in
$H(W^k)\utensor H(W^{k'})$.

\subsection*{Step 3: Factorisation of $\mathcal H_B$, $B$ a negative
  block, into braided Cherednik algebras $\protect\rd{\mathcal H}(W^k)$}

Now assume that $B$ is a negative block.
The group $W_B$ may not be a complex reflection group. 
By an observation at Step~1 of this proof, Theorem~\ref{thm:structure} and
Lemma~\ref{lem:compref}(3), $W_B$ is generated by some elements of the
form
$$
 t_j^{(-1)}, \quad 
 t_j^{(-1)}t_j^{(\eta)}
 \ (\eta\ne1\ \text{root of unity}), 
 \quad\text{or}\quad
 t_j^{(-1)}(ij)t_i^{(\varepsilon)}t_j^{(\varepsilon^{-1})}
 \ (\varepsilon\in\CC^\times),
 \qquad i,j\in B,
$$
besause only such elements may appear in $C_{ij}$ ($i,j\in B$) with
nonzero coefficients. 
We rewrite the list of possible
generators of $W_B$ as 
$$
  t_j^{(\eta)}\quad(\eta\ \text{any root of unity}), 
  \qquad
  \sigma_{ij}^{(\varepsilon)}\quad (\varepsilon\in\CC^\times),
 \qquad i,j\in B.  
$$
Call two indices $i,j\in B$ linked, if $i=j$ or $W_B$ contains
an element $\sigma_{ij}^{(\varepsilon)}$ for some $\varepsilon\in
\CC^\times$. The relation ``linked'' is symmetric and transitive,
because
$$
   \sigma_{ij}^{(\varepsilon)} = \sigma_{ji}^{(-\varepsilon^{-1})}, 
\qquad
   (\sigma_{ij}^{(\varepsilon)})^{-1} \sigma_{jk}^{(\delta)} 
   \sigma_{ij}^{(\varepsilon)} = \sigma_{ik}^{(\varepsilon\delta)}.
$$
Let $B=\mathcal O^1 \cup \dots \cup \mathcal O^l$ be the partition of $B$
into equivalence classes, and denote $V^k =\oplus\{\CC x_i \mid i\in
\mathcal O^k\}$. The generating set for $W_B$ is 
partitioned into pairwise commuting subsets 
$\{t_j^{(\eta)},\sigma_{ij}^{(\varepsilon)} \in W_B
\mid  i,j\in \mathcal O^k\}$, $k=1,\dots,l$, thus $W_B$ is a direct
product $W^1\times \dots \times W^l$ of groups acting on the direct
sum $V_B=V^1\oplus \dots\oplus V^l$ of spaces.  The algebra $\mathcal
H_B$ is then a tensor product 
$$
   \mathcal H_B = (S_\minusone(V^1) \tensor \CC W^1 \tensor S_\minusone(V^{1*}))
   \tensor \dots \tensor 
(S_\minusone(V^l) \tensor \CC W^l \tensor S_\minusone(V^{l*})),
$$
of vector spaces,
where $V^{k*}=\oplus\{\CC \uy_i \mid i\in \mathcal O^k\}$.
Observe that each  $S_\minusone(V^k) \tensor \CC W^k \tensor
S_\minusone(V^{k*})$ is a subalgebra, because the braided 
commutator $C_{ij}$ of $\uy_j$ and
$x_i$ ($i,j\in \mathcal O^k$) may only contain generators
$\sigma_{ij}^{(\varepsilon)}$, $t_j^{(\eta)}$ of $W_B$
that lie in $W^k$. Thus, 
$\rd{\mathcal H}(W^k):= S_\minusone(V^k) \tensor \CC W^k \tensor
S_\minusone(V^{k*})$ is a braided Cherednik algebra. 

Let us show that 
$$
   \mathcal H_B =  \rd{\mathcal H}(W^1) \tensor \dots \tensor
   \rd{\mathcal H}(W^l) 
$$
is a braided tensor product of algebras in the category 
$\mathcal M_{\Gamma_\qsub,\mathcal R_\qsub}$. First of all, 
$V^k$ and $V^{*k}$ are $\Gamma_\qsub$\dash submodules of $V_B$ and
$V_B^*$, respectively, because they are spanned by simultaneous 
eigenvectors of $\Gamma_\qsub$. Next, since $W^k$ is generated by its
elements that 
appear with nonzero coefficients in braided commutators of $V^k$ and
$V^{k*}$, and the braided commutator map is $\Gamma_\qsub$\dash
equivariant (proved in Step 1), the group algebra $\CC W^k$ is a
$\Gamma_\qsub$\dash submodule of $\CC W_B$. 
This gives the $\Gamma_\qsub$\dash structure on $\rd{\mathcal
  H}(W^k)$. 
It remains to check that
the commutation relations between $\rd{\mathcal H}(W^k)$ and
$\rd{\mathcal H}(W^{k'})$  inside $\mathcal H_B$ are the same as in
the braided tensor product 
$\rd{\mathcal H}(W^k)\utensor \rd{\mathcal H}(W^{k'})$. 

In $\rd{\mathcal H}(W^k)\utensor \rd{\mathcal H}(W^{k'})$, the
variables $x_i$ and $x_j$ ($i\in \mathcal O^k$, $j\in \mathcal
O^{k'}$) $q_{ij}$\dash commute, and the same happens in $\mathcal
H_B$. Similarly for $\uy_i$ and $\uy_j$. Furthermore, $\uy_j$ and
$x_i$ $q_{ij}$\dash commute in $\rd{\mathcal H}(W^k)\utensor
\rd{\mathcal H}(W^{k'})$, and the same happens in $\mathcal H_B$,
since $C_{ij}=0$: $i$, $j$ are not linked,  hence there is no element
$\sigma_{ij}^{(\varepsilon)}$ in $W_B$. Finally, for any $w\in W^k$ and
$w'\in W^{k'}$ one has $\mathcal R_\qsub(w(\gamma_i),w'(\gamma_j))=
\mathcal R_\qsub(\gamma_i,\gamma_j)=-1$, simply because $w(\gamma_i)\in |V^k|\subset |V_B|$, $w'(\gamma_j)\in
|V^{k'}|\subset |V_B|$, and the value of $\mathcal R_\qsub$ at any
pair of distinct elements of $|V_B|$ is $-1$ as $B$ is a negative
block. Hence by Proposition~\ref{prop:braidedprod}, $W^{k'}$ commutes
with $\rd{\mathcal H}(W^k)$ and $\rd{\mathcal H}(W^{k'})$ commutes
with $W^k$ in $\rd{\mathcal H}(W^k)\utensor \rd{\mathcal
  H}(W^{k'})$. But the same relations hold in $\mathcal H_B$, thus
$\mathcal H_B =  \rd{\mathcal H}(W^1)\utensor \dots \utensor
\rd{\mathcal H}(W^l)$.

\subsection*{Step 4: Proof that each of the algebras $\protect\rd{\mathcal
    H}(W^k)$ is an ``elementary'' braided Cherednik algebra}

Recall that we are working with a negative block $B$ of indices, and
have already factorised $\mathcal H_B$ into braided Cherednik algebras 
$\rd{\mathcal H}(W^k)$, $k=1,\dots,l$, where $W^k\subset \GL(\mathcal
V^k)$ and $V^k=\oplus \{\CC x_i \mid i\in\mathcal O^k\}$. 
We fix an index $k$ and will
show that $\rd{\mathcal H}(W^k)$ is isomorphic to one of the
``elementary'' braided
Cherednik algebras, listed in $(1)$--$(4)$ before the Theorem. 
Without the loss of generality we may assume that
the set $\mathcal O^k$ of indices is $1,2,\dots,d$. 
If $d=1$, then $W^k$ is a cyclic group (an irreducible complex
reflection group of rank $1$),
and $\rd{\mathcal H}(W^k)$ is a rational Cherednik algebra
isomorphic to $H_c(W^k)$ or $H_{0,c}(W^k)$.

Assume $d\ge 2$. All indices in  $\{1,\dots,d\}$ are pairwise linked, that is 
for each pair $i,j\in\{1,\dots,d\}$ of distinct indices, there is at
least one nonzero number --- call it $\varepsilon_{ij}$ --- such that 
$\sigma_{ij}^{(\varepsilon_{ij})}\in W^k$. 

We may assume that $\varepsilon_{ij}=1$ for all $i\ne j$ in
$\{1,\dots,d\}$. Indeed, we may change the basis of $V^k$ 
by rescaling the variable $x_{i}$ by a factor of 
$\varepsilon_{12}\varepsilon_{23}\dots \varepsilon_{i-1,i}$ 
and denote the new basis again by $\{x_i\}$. The braided Cherednik
algebra structure of $\rd{\mathcal H}(W^k)$ obviously does not change
under such rescaling, nor does the action of $\Gamma_\qsub$. We apply
rescaling to the dual basis in $V^{k*}$ 
so that $\{x_i\}$ and $\{\uy_i\}$, $i=1,\dots,d$, 
remain a pair of dual bases. Now with respect to the new basis,
$W^k$ contains $\sigma_{i-1,i}^{(1)}$ 
and hence also contains $\sigma_{i,i-1}^{(1)}=(\sigma_{i-1,i}^{(1)})^{-1}$
for each $i=2,\dots,d$.
It then follows from the relation 
$\sigma_{ba}^{(1)} \sigma_{bc}^{(1)} \sigma_{ab}^{(1)} =
\sigma_{ac}^{(1)}$ that $W^k$ contains $\sigma_{ij}^{(1)}$ for any
pair $i\ne j$, $i,j=1,\dots,d$. 

Besides $\sigma_{ij}^{(1)}$, the group $W^k$ may have some other
generators, namely some of $t_j^{(\eta)}$ and
$\sigma_{ij}^{(\varepsilon)}$. 
We replace each generator $t_j^{(\eta)}$ by
$$
     t_1^{(\eta)} = \sigma_{j1}^{(1)} t_j^{(\eta)} \sigma_{1j}^{(1)}
     \quad \in W^k,
$$  
and each generator $\sigma_{ij}^{(\varepsilon)}$ by 
$$
   t_1^{(\varepsilon)} t_2^{(\varepsilon^{-1})}  = \sigma_{21}^{(1)} 
   \sigma_{j2}^{(1)} \sigma_{i1}^{(1)} \sigma_{ij}^{(\varepsilon)}
   \sigma_{1i}^{(1)} \sigma_{2j}^{(1)}.
$$
Thus the new set of generators for the same group $W^k$ contains 
$\sigma_{ij}^{(1)}$ for all $i\ne j$, $i,j=1,\dots,d$, and also 
$t_1^{(\eta)}$ and $t_1^{(\varepsilon)} t_2^{(\varepsilon^{-1})}$ for
some unknown choice of the $\eta$'s and $\varepsilon$'s. 
Let
$$
\mathcal C = 
\{\varepsilon \in \CC^\times \mid
t_1^{(\varepsilon)}t_2^{(\varepsilon^{-1})}\in W^k\}, 
\qquad
\mathcal C' = \{\eta \in \CC^\times \mid t_1^{(\eta)}\in W^k\}. 
$$ 
Then $\mathcal C$ (respectively $\mathcal C'$) 
is a finite subgroup of $\CC^\times$ because it
is the inverse image of $W^k$ under a group monomorphism 
$\varepsilon \mapsto t_1^{(\varepsilon)}t_2^{(\varepsilon^{-1})}$ 
(respectively $\eta\mapsto t_1^{(\eta)}$) from $\CC^\times$ to $\GL(V^k)$. 
Moreover, $\mathcal C'\subseteq \mathcal C$ because if 
$t_1^{(\eta)}\in W^k$, then 
$t_1^{(\eta)}t_2^{(\eta^{-1})} =
t_1^{(\eta)}[\sigma_{12}^{(1)}t_1^{(\eta)}\sigma_{21}^{(1)} ]^{-1}$ is
also in $W^k$. Besides that, $\mathcal C$ contains $-1$, as
$t_1^{(-1)}t_2^{(-1)} = (\sigma_{12}^{(1)})^2\in W^k$, hence $\mathcal
C$ is of even order. We have proved that 
$$
   W^k = W_{\mathcal C, \mathcal C'}(d).
$$
By Proposition~\ref{prop:unique}, $\rd{\mathcal H}(W^k)$ is isomorphic
to $\rd{\mathcal H}_c(W_{\mathcal C, \mathcal C'}(d))$ or to 
$\rd{\mathcal H}_{0,c}(W_{\mathcal C, \mathcal C'}(d))$ for some
parameter~$c$. 
\end{proof}
\begin{remark}
Note that to form a braided tensor product 
$\rd{\mathcal H}=\rd{\mathcal H}_1\utensor \dots \utensor \rd{\mathcal H}_m$ 
of braided Cherednik algebras
of irreducible groups (algebras listed in $(1)$--$(4)$ above), one
needs $m(m-1)/2$ extra nonzero complex parameters $r_{ab}$, $1\le
a<b\le m$.  
The matrix $\qq$ for $\rd{\mathcal H}$ can be written 
as a block matrix with $m^2$ blocks $M_{ab}$, $a,b=1,\dots,m$,
such that:
\begin{itemize}
\item the size of  $M_{ab}$ is $(\mathrm{rank}\ \rd{\mathcal H}_a
  \times \mathrm{rank}\ \rd{\mathcal H}_b)$;
\item in a diagonal block $M_{aa}$, all entries are $1$ or else all
  entries outside the main diagonal are~$-1$;
\item in an off\dash diagonal block $M_{ab}$ where $a<b$ (respectively $a>b$),
  all entries are   equal to $r_{ab}$ (respectively $r_{ab}^{-1}$).
\end{itemize}
The commutation relations in the braided tensor product include
$xx'=r_{ab} x' x$ whenever $x$ is one of the $x_i$ variables in
$\rd{\mathcal H}_a$ and $x'$ is one of the $x_i$ variables in
$\rd{\mathcal H}_b$.
\end{remark}

\section{Universal embeddings and braided Dunkl operators}
\label{sect:Braided Dunkl operators}

In the last Section of the paper, we embed braided Cherednik algebras
in  modified Heisenberg quadratic doubles, introduced here.
We use this result to arrive at the explicit
formulae for braided Dunkl operators. 

\subsection{Degenerate $\qq$-Cherednik algebras and Heisenberg
  quadratic doubles}

%By Proposition~\ref{prop:brred}, 
%a braided Cherednik algebra is a braided reduction (in particular, a
%subalgebra) of some $\qq$\dash Cherednik algebra $\mathcal
%H_c(W)$. The latter is a quadratic double, hence we may apply
%the results of Section~\ref{sect:qd} to embed $\mathcal
%H_c(W)$ in a Heisenberg
%quadratic double; this will provide an embedding for its braided
%reduction as well. 

We say that a $\qq$-Cherednik algebra of the form $\mathcal H_{0,c}(W)
= S_\qsub(V) \tensor \CC W \tensor S_\qtran(V^*)$ is degenerate, if 
the commutator of $y\in V^*$ and $x\in V$ in $\mathcal H_{0,c}(W)$ 
is 
$
   \beta'(y\tensor x) = yx-xy = \sum_B \gamma_\Bsub \sum_s c_s \langle x,
   \alpha_s^\vee \rangle \langle \alpha_s, y\rangle s 
$
(compare with Theorem~\ref{thm:structure}). Here and below $B$ are blocks of indices with respect to the matrix $\qq$, and
$s$ runs over complex reflections in $\GL(V_B)\subset \GL(V)$; we
continue to use the notation introduced 
in Section~\ref{sect:q_cher}.
We would like to construct a morphism (not necessarily injective)
from a degenerate $\qq$\dash Cherednik algebra
to a Heisenberg quadratic double over $W$. 
This is done via Theorem~\ref{thm:perfect}; the crucial step is to represent 
$\mathcal H_{0,c}(W)$ as a $\star$\dash product
(Definition~\ref{def:stardiamond}) of two quadratic
doubles. This is done as follows:
% The formula for $\beta'$ can be written
%as
%$\sum_{B,s} \langle L_{\gamma_\Bsub s}(x),y\rangle
%  \gamma_\Bsub s$, 
%%\qquad
%%\text{where}\  
%where $L_{\gamma_\Bsub s}(x)= c_s \langle x, \alpha_s^\vee \rangle
%  \alpha_s
%$. We have
%%In the same way as for rational Cherednik algebra we may represent 
%%$\mathcal H_{0,c}(W)$ as a result of $\star$\dash operation
%%(Definition~\ref{def:stardiamond}):
$$
   \mathcal H_{0,c}(W) = \mathcal H_{0,c}(W) \star \mathcal H_{0,c_0}(W)\,,
$$  
where the value of the fixed parameter $c_0$ at $s$ is 
$\langle \alpha_s,\alpha_s^\vee\rangle^{-1}$. 
%We refer to elements $\gamma_\Bsub s\in W$ as \emph{$\qq$\dash
%  reflections} 
%and 
Application of  Theorem~\ref{thm:perfect} now yields 
a Yetter\dash Drinfeld module 
$$
           Y_\qsub = \bigoplus_{B,s} \CC \cdot [\gamma_\Bsub s].       
$$
%Here and below the sum is taken over the pairs $B,s$ where $B$ is a
%           block of indices 
%           and $s$ is a complex reflection in $\GL(V_B)$, such that
%           $\gamma_\Bsub s\in W$. 
We refer to elements $\gamma_\Bsub s\in W$ as \emph{$\qq$\dash
  reflections}. The YD module structure on  $Y_\qsub$ 
is induced by its embedding in the YD module $\CC W \tensor V$ 
%(as in
%the proof of Theorem~\ref{thm:perfect}) 
via 
$ [\gamma_\Bsub s]=\gamma_\Bsub s \tensor \alpha_s$.
The dual YD module 
$
           Y_\qsub^* = \bigoplus_{B,s} \CC \cdot [\gamma_\Bsub s]^*   $
embeds in $\CC W \tensor V^*$ via $ [\gamma_\Bsub s]^* = (\gamma_\Bsub s)^{-1}
           \tensor \langle \alpha_s,\alpha_s^\vee\rangle^{-1}
           \alpha_s^\vee$. 
%Theorem~\ref{thm:emb0} implies the following
By Theorem~\ref{thm:perfect}, 
%\begin{theorem}
%\label{thm:emb-0}
the $W$\dash equivariant maps $\mu_c\colon V \to
Y_\qsub$, $\nu\colon V^* \to Y_\qsub^*$
given by
$$
   \mu_c(x) = \sum_{B,s} c_s \langle x,\alpha_s^\vee\rangle [\gamma_\Bsub
   s],
   \qquad 
   \nu(y) = \sum_{B,s} \langle \alpha_s, y \rangle [\gamma_\Bsub s]^*,
$$
extend to an algebra morphism
$
  \mathcal H_{0,c}(W) \to A_{Y_\qsub} = 
   S(Y_\qsub,\Psi_{Y_\qsub}) \tensor 
   \CC W \tensor
  S(Y_\qsub^*,\Psi_{Y_\qsub}^*)
%\qed
$.
%\end{theorem}
\startcmm
%%%%%%%%%%%%%%%%%%%%%%%%%%%%%%%%%%%%%%%%%%%%%%%%%%%%%%%%%%
\subsection{The algebra $S(Y_\qsub,\Psi_{Y_\qsub})$}  
Recall from Section~\ref{sect:yd} that the algebra
$U^-:= S(Y_\qsub,\Psi_{Y_\qsub})$ that we have
just shown to contain $\qq$\dash commuting elements $\mu_c(x_i)$, has
presentation 
$$
 U^- = T(Y_\qsub) / \lgen \ker (\id +
 \Psi_{Y_\qsub})\rgen. 
$$
Let $B$, $C$ be distinct blocks   of indices. 
Note that any $s\in \GL(V_B)$ acts on $[\gamma_\Csub t]\in
Y_\qsub$ trivially, and $\gamma_\Bsub$ acts by the scalar
$q_{B,C}$. This is because $[\gamma_\Csub t]=\gamma_\Csub t \tensor \alpha_t$
in the module $\CC W \tensor V$, so that $\gamma_\Bsub$ conjugates
$\gamma_\Csub t $  trivially and acts on $\alpha_t\in V_C$ by $q_{B,C}$.
Thus the
braiding between the generators 
$[\gamma_\Bsub s]$ and $[\gamma_\Csub t]$ of $U^-$ is 
$$
\Psi_{Y_\qsub}( [\gamma_\Bsub s]\tensor [\gamma_\Csub t])
= \gamma_\Bsub (s( [\gamma_\Csub t] )) \tensor [\gamma_\Bsub s]
= q_{B,C}  [\gamma_\Csub t] \tensor [\gamma_\Bsub s]. 
$$
As $q_{C,B}=q_{B,C}^{-1}$, this means that the generators of $U^-$
corresponding to these two blocks $q_{B,C}$\dash commute. The whole
algebra $U^-$ factorises as 
$$
   U^- = U^-_{B_1} \utensor \dots \utensor U^-_{B_N}
$$
where $B_1,\dots,B_N$ are all blocks of indices, and 
and $U^-_{B_k}$ is the
subalgebra generated by generators $[\gamma_{B_k} s]$ of $U^-$. 
This is a braided tensor product in the braided category 
$\mathcal M_{\bar \Gamma,\bar{\mathcal R}}$
where the group
$$
\bar \Gamma = \langle \gamma_{B_1},\dots,\gamma_{B_N} \rangle
\subset \GL(V)
$$
is  Abelian, with bicharacter defined by 
$\bar{\mathcal R}(\gamma_\Bsub,\gamma_\Csub)=q_{B,C}$.  

The principal part of the structure of $U^-$ is therefore the relations
in each of the subalgebras $U^-_{B_k}$. We thus assume that all
indices are in a single block $B$. If $B$ is positive, we are dealing
with a rational Cherednik algebra, and the algebra $U^-$ was already
considered in Section~\ref{sect:yd}; see Remark~\ref{rem:BK}.

\subsection{Classical and twisted classical complex reflection groups}

If the block $B=\{1,\dots,n\}$ is negative, one has $\gamma_\Bsub=-\id$. 
The group $W$ is generated by $-\id \cdot s$, where $s$ is a complex
reflection in an $n$\dash dimensional space which acts imprimitively,
preserving $\{\CC x_1,\dots,\CC x_n\}$. 
It is not too difficult to list finite subgroups of $\GL_n(\CC)$ which are so
generated. Recall the group $\Symm_n\le \GL_n(\CC)$ of permutation
matrices and denote 
$$
(\Symm_n)_\sgn = \{ \det(\pi)\cdot \pi \mid \pi \in \Symm_n \}.
$$
\begin{proposition}
\label{prop:ABCD}
Subgroups of $\GL_n(\CC)$ that are generated by $-\id \cdot s$ with $s$ 
imprimitive complex reflections, and are not decomposable into direct products
of smaller such subgroups, are
\begin{itemize}
\item[$(i)$]
$G(m,p,n)$, $m$ even, $mn/p$ even;
\item[$(ii)$] $G(m,p,n)\times \{\pm \id\}$, $m$ odd;
\item[$(iii)$] $G(m,p,n)_\sgn:= T(m,p,n)\lcprod (\Symm_n)_\sgn$, $m$ even,
  $m/p$ odd;
\item[$(iv)$] $G(m,m,n)_\sgn:= T(m,m,n)\lcprod (\Symm_n)_\sgn$, $m$ odd.\qed
\end{itemize}\end{proposition}
To each of the
groups $(i)$--$(iv)$ there is associated a $\negone$\dash Cherednik
algebra and the corresponding algebra $U^-=
S(Y_\qsub,\Psi_{Y_\qsub})$  containing
$n$ anticommuting elements of degree $1$.
Some of  these quantum (anti)symmetric algebras $U^-$ were
found earlier:

\begin{example}[Algebras of flat connections]
The smallest group in the family $(iv)$ in Proposition~\ref{prop:ABCD} is
$G(1,1,n)_\sgn=(\Symm_n)_\sgn$. The corresponding algebra 
$$
U^- = \frac{\CC\langle e_{ij} \mid  i,j=1,\dots,n, \ i\ne j \rangle}{
\lgen e_{ij}=e_{ji}, \ 
e_{ij}^2 = e_{ij}e_{kl}+e_{kl}e_{ij}=e_{ij}e_{jk}+e_{jk}e_{ki}+e_{ki}e_{ij}=0,\
i,j,k,l\ \text{distinct}\rgen}
$$ 
was introduced by Majid in \cite{Mnoncomm}. The elements 
$\theta_i = \sum_{j\ne i}e_{ij}$ anticommute and
generate a subalgebra termed algebra of flat connections with constant
coefficients for the
noncommutative differential cohomology of $\Symm_n$. 

The Weyl groups $B_n$ and $D_n$ correspond to $G(2,1,n)$ (the smallest
in the family $(i)$)
and twisted $G(2,2,n)$ (the smallest in the family $(iii)$) of
Proposition~\ref{prop:ABCD}, and 
the relevant algebras of flat connections were found by Kirillov and
Maeno \cite{KMexterior}. 
\end{example}
%%%%%%%%%%%%%%%%%%%%%%%%%%%%%%%%%%%%%%%%%%%%%%%%%%%%%%%%%%%%%%%%%%%%%%%%%
\endcmm

%\subsection{Modified Heisenberg quadratic doubles}
\subsection{$\qq$-Cherednik
  algebras are subalgebras in modified Heisenberg quadratic doubles}
\label{subsect:embq}

One can obtain a version of the above morphism $\mathcal H_{0,c}(W)
\to A_{Y_\qsub}$ for non\dash degenerate $\qq$\dash Cherednik
algebras. A new ingredient for this is the 
operation $\diamond$, introduced in Section~\ref{sect:qd}. 

Suppose that $Y$ is a module over a group $W$, and
  that $Y$ has two different Yetter\dash Drinfeld structures over $W$;
  that is, two $W$\dash gradings satisfying the Yetter\dash Drinfeld
  axiom with respect to the same $W$\dash action on $Y$. 
  These two YD structures give rise to
  two braidings $\Psi_1$, $\Psi_2$ on $Y$ and 
  two Heisenberg quadratic doubles 
$$
   A_Y^i = S(Y,\Psi_i)\tensor \CC W \tensor S(Y^*,\Psi_i^*), \qquad
   i=1,2. 
$$
\begin{definition}
The quadratic double $A_Y^1\diamond A_Y^2$ is called a \emph{modified
Heisenberg quadratic double} of the two Yetter\dash Drinfeld
structures on $Y$. 
\end{definition}
By definition of $\diamond$, the triangular decomposition 
of $A_Y^1\diamond A_Y^2$ is 
$$
  \frac{T(Y)}{\lgen \ker(\id+\Psi_1)\cap \ker(\id+\Psi_2) \rgen }
\tensor \CC W \tensor  \frac{T(Y^*)}{\lgen \ker(\id+\Psi_1^*)\cap
  \ker(\id+\Psi_2^*) \rgen}\,.
$$
Now, by  Theorem~\ref{thm:structure} an arbitrary $\qq$\dash Cherednik algebra can be written as 
$S_\qsub(V)\tensor \CC W \tensor S_\qtran(V^*)$ 
%$\mathcal H_c(W)$ 
with the commutator 
%$\beta\colon V^*\tensor V \to
%\CC W$ 
$
     \beta(y\tensor x) = \sum_B \gamma_\Bsub (x,y)_B + \beta'(y\tensor x)
$
%
%+\sum_s
%     c_s \langle x, 
%   \alpha_s^\vee \rangle \langle \alpha_s, y\rangle s\bigr)
%$
between $V^*$ and $V$, with $\beta'(y\tensor x)$ as above. 
Denote this algebra by $\mathcal H_{(\cdot,\cdot),c}(W)$, where
$(\cdot, \cdot)=\sum_B \gamma_B(\cdot,\cdot)_B$ is the $\CC
\bar\Gamma$\dash valued pairing between $V$ and $V^*$. 
% to show the
%dependence on the pairings $ (x,y)_B$.
Clearly, 
$$
   \mathcal H_{(\cdot,\cdot),c}(W) = 
   \mathcal H_{0,c}(W) \diamond
   \mathcal H_{(\cdot,\cdot),0}(W),
$$
where $\mathcal H_{0,c}(W)$ is the degenerate $\qq$\dash Cherednik
algebra with commutator $\beta'$ as above, and 
$\mathcal H_{(\cdot,\cdot),0}(W)$
is the $\qq$\dash Cherednik algebra with commutator $\beta-\beta'$. 
We have already constructed a morphism $\mathcal H_{0,c}(W)\to
A_{Y_\qsub}$, and will now turn to the algebra $ \mathcal
H_{(\cdot,\cdot),0}(W)$.

Recall that the module $Y_\qsub$ has $W$\dash grading given by
$$
  |[\gamma_\Bsub s]| = \gamma_\Bsub s.
$$
Assume that the group $W$ contains the Abelian group
$$
\bar \Gamma = \langle \gamma_{B_1},\dots,\gamma_{B_N} \rangle
\subset \GL(V)\,;
$$ 
since $W$
permutes subspaces $V_B$, the group $\bar \Gamma$ is normal in
$W$. It follows that we can introduce the second,  $\bar\Gamma$\dash
valued grading 
$$
    | [\gamma_\Bsub s]|_{\bar \Gamma} := \gamma_\Bsub
$$
on $Y_\qsub$, which  also makes $Y_\qsub$ a Yetter\dash
Drinfeld module over $W$. 
This second YD structure leads to a
Heisenberg quadratic double
$$
       A_{Y_\qsub}^{\bar \Gamma} \cong 
       S(Y_\qsub,\tau_{\bar \Gamma}) \tensor 
       \CC W \tensor  S(Y_\qsub^*,\tau_{\bar \Gamma}^*),
$$
where the braiding $\tau_{\bar \Gamma}$ is given by 
$\tau_{\bar \Gamma}([\gamma_\Bsub s]\tensor [\gamma_\Csub t]) = 
q_{B,C} [\gamma_\Csub t] \tensor [\gamma_\Bsub s]$. Clearly,
$S(Y_\qsub,\tau_{\bar \Gamma}^*)$ is nothing but 
a $(q_{B,C})$\dash polynomial algebra of $Y_\qsub$.

%We are free to introduce $W$\dash invariant parameters $t_{B,s}\in
%\CC$ and to impose the commutation relation 
%$$
%[\gamma_\Bsub s]^*[\gamma_\Csub t] - [\gamma_\Csub t][\gamma_\Bsub s]^* =
%\delta_{B,C}\delta_{s,t}t_{B,s} \gamma_\Bsub
%$$
%in $A_{Y_\qsub}^{\bar \Gamma}$.

\begin{lemma}
\label{lem:val-}
%If $c_s$ are generic, for some values of parameters $t_{B,s}$ in
%$A_{Y_\qsub}^{\bar \Gamma}$
For some pairings $(\cdot,\cdot)_B$,  
the maps $\mu_c\colon V \to Y_\qsub$, $\nu\colon V^* \to
Y_\qsub^*$ defined above extend to a morphism
$\mathcal H_{(\cdot,\cdot),0}(W) \to
%\hookrightarrow
A_{Y_\qsub}^{\bar \Gamma}$. This morphism is injective if the roots $\{\alpha_s
\mid c_s\ne 0\}$ span $V$.  
\end{lemma}
\begin{proof}
Let $B,C\subset \{1,\dots,n\}$ be blocks of indices. For any $i\in B$,
$j\in C$, $i\ne j$ the relation $x_i x_j = q_{B,C}x_j x_i$ holds in
$S_\qsub(V)$. Note that $\mu_c(x_i)$ is a linear combination of 
$[\gamma_B s]$ where $s\in W$ runs over complex reflections in $V_B$,
and similarly $\mu_c(x_j)$ is a combination of 
$[\gamma_C t]$ where $t$ runs over complex reflections in $V_C$. 
The relation $[\gamma_B s][\gamma_C t] = q_{B,C}
[\gamma_C t][\gamma_B s]$ holds in $ S(Y_\qsub,\tau_{\bar \Gamma})$
for any such $s,t$, hence $\mu_c$ extends to a morphism $S_\qsub(V)
\to  S(Y_\qsub,\tau_{\bar \Gamma})$.
If $\{\alpha_s\mid c_s\ne 0\}$ span $V$, $\mu_c\colon V \to Y_\qsub$
is injective; then $\mu_c(x_i)$ are linearly independent vectors in
$Y_\qsub$ that generate a subalgebra in $S(Y_\qsub,\tau_{\bar
  \Gamma})$ isomorphic to the $\qq$\dash polynomial algebra, therefore
the extension of $\mu_c$ to $S_\qsub(V)$ is injective. 
Similarly, $\nu$ extends to a morphism $S_\qtran(V^*)
\to  S(Y_\qsub^*,\tau_{\bar \Gamma}^*)$, injective if all
$\alpha^\vee_s$ span $V^*$ (that is, all $\alpha_s$ span $V$). 

The $\CC \bar \Gamma$\dash valued 
commutator of $\nu(y_j)$ and $\mu_c(x_i)$ in $A_{Y_\qsub}^{\bar
  \Gamma}$   is equal to $\delta_{B,C}\gamma_B 
\sum_s c_s \langle x,\alpha_s^\vee\rangle \langle \alpha_s,
  y\rangle$ where $s\in W$ runs over complex reflections in $V_B$.
Let the same formula define the pairing $\gamma_B (x_i,y_j)_B$ between
$V$ and $V^*$. Then
we have a morphism $\mathcal H_{(\cdot,\cdot),0}(W) \to
A_{Y_\qsub}^{\bar \Gamma}$ that extends the maps $\mu_c$, $\nu$. 
It follows from the triangular decomposition that if this morphism is
injective on $S_\qsub(V)$ and $S_\qtran(V^*)$, it 
is injective on $\mathcal H_{(\cdot,\cdot),0}(W)$.
%
%, whereas the commutator of $y_j$
%  and $x_i$ in $\mathcal H_{(\cdot,\cdot),0}(W) $ is 
%$\delta_{B,C}\gamma_B (x_i,y_j)_B$. We have to ensure that these two 
%$\CC\bar\Gamma$\dash valued pairings between $V^*$ and $V$ 
%can be made equal by choosing the parameters $t_{B,s}$. 
%As far as the genericity of $c_s$ is concerned, it is well enough to 
%assume that all $c_s$ are nonzero.
%
%Note that both pairings, viewed as elements of $V\tensor V^* \tensor
%\CC\bar \Gamma$, lie in the image of the $W$\dash equivariant map 
%$V\tensor V^* \to V\tensor V^*\tensor \CC \bar \Gamma$ given by
%$x\tensor y_i \mapsto x\tensor y_i \tensor \gamma_B$ where $B$ is the
%block containing the index $i$. It is thus enough to show that 
%$\alpha_s\tensor\alpha_s^\vee$
\end{proof}
Applying Proposition~\ref{prop:diamond}, we obtain 
%consider the modified Heisenberg quadratic double
%$$
%A'_{Y_\qsub} := A_{Y_\qsub}\diamond
%A_{Y_\qsub}^{\bar \Gamma} .
%$$
\begin{theorem}
\label{th:embdg modified}
If the parameter $c$ is such that the roots $\{\alpha_s\mid c_s\ne 0\}$ 
span $V$, then for some pairings $(\cdot,\cdot)_B$, 
the above maps $\mu_c\colon V \to Y_\qsub$, $\nu\colon V^* \to
Y_\qsub^*$ extend to an embedding
$\mathcal H_{(\cdot,\cdot),c}(W) \hookrightarrow A_{Y_\qsub}\diamond A^{\bar\Gamma}_{Y_\qsub}$.
\qed
\end{theorem}
%\begin{remark}
%\label{rem:fine}
%The commutation relation 
%$[\gamma_\Bsub s][\gamma_\Csub t] - [\gamma_\Csub t][\gamma_\Bsub s] =
%\delta_{B,C}\delta{s,t}\gamma_\Bsub$ in $A_{Y_\qsub}^{\bar \Gamma}$
%can be fine\dash tuned by introducing a $W$\dash invariant factor
%$t_{B,s}$ on the right. For any given pairings $(\cdot,\cdot)_B$, the
%values of these extra factors can be so chosen that
%$A'_{Y_\qsub} = A_{Y_\qsub}^{\bar \Gamma} \diamond
%A_{Y_\qsub}$ contains 
%$\mathcal H_{(\cdot,\cdot),c}(W)$. The extra factors do not affect the
%tensorands in the triangular decomposition 
%$$
%A'_{Y_\qsub}
%= S'(Y_\qsub, \Psi_{Y_\qsub}) \tensor 
%\CC W \tensor S'(Y_\qsub^*, \Psi_{Y_\qsub}^*),
%$$
%where
%$$
%S'(Y_\qsub, \Psi_{Y_\qsub}) = T(Y_\qsub) / 
%\lgen\ker(\id + \Psi_{Y_\qsub})\cap
%\wedge^2_{\bar\Gamma}Y_\qsub\rgen. 
%$$
%Note that the algebra $S'(Y_\qsub, \Psi_{Y_\qsub})$ will
%  still be a braided tensor product of components indexed by blocks
%  in the same way as $U^-$ above, since the relations 
%$
%[\gamma_\Bsub s][\gamma_\Csub t]=q_{B,C}
%[\gamma_\Csub t][\gamma_\Bsub s]
%$, present both in $U^-$ and in $S(Y_\qsub,\tau_{\bar \Gamma})$,
%make their way into the algebra $S'(Y_\qsub,
%  \Psi_{Y_\qsub})$. 
%\end{remark}
\begin{remark}
The $W$\dash invariant $\CC \bar \Gamma$\dash valued pairing
$(\cdot,\cdot)$ between $V$ and $V^*$ in the Theorem depends on the
parameter $c$. As in any $\qq$\dash Cherednik algebra, this pairing is
of a special kind: namely, it is obtained from a scalar $W$\dash
invariant pairing by the change of variables as in 
Corollary~\ref{cor:anticommutator}. If the group $W$ is irreducible, such
pairing is unique up to a scalar factor. One deduces that any  $\qq$\dash
Cherednik algebra of an irreducible group, with at least one nonzero
parameter $c_s$, embeds in a modified Heisenberg double.
\end{remark}

The following is left as an exercise to the reader:
\begin{exercise}
Describe an embedding of the twist of a rational Cherednik algebra,
introduced just before Corollary~\ref{cor:anticommutator}, in an appropriate
version of a modified Heisenberg double. 
\end{exercise} 
\begin{remark}
Embedding of a braided Cherednik algebra $\rd{\mathcal H}(W)$ 
in a modified braided Heisenberg double, described in
Theorem~\ref{th:emb0}, is obtained by first embedding $\rd{\mathcal
  H}(W)$  in a $\qq$\dash Cherednik algebra $\mathcal H(\widetilde W)$ 
with $\widetilde W = W \cdot \Gamma_\qsub \cdot \bar \Gamma$, then
applying Theorem~\ref{th:embdg modified} 
%(and Remark~\ref{rem:fine})
to $\mathcal H(\widetilde W)$.
\end{remark}

\subsection{Braided Dunkl operators}

\label{subsect:Braided Dunkl operators}

We will now consider the braided Cherednik algebra $\rd{\mathcal
  H}_c(W)$ of the irreducible group $W=W_{\mathcal
  C,\mathcal C'})$, as introduced in Definition~\ref{def:negbraided}. 
The parameter $c$ is a function $c\colon \mathcal C'\to \CC$ 
(with the exception of rank $n=2$, see Remark~\ref{rem:rank2}; we are
  going to ignore this exception and claim that in rank $2$, the proof
  may easily be modified as appropriate). 
The algebra has triangular decomposition 
$\rd{\mathcal   H}_c(W) = S_\minusone(V)\tensor \CC W \tensor
  S_\minusone(V^*)$ where $V$ is spanned by $x_1,\dots,x_n$ and $V^*$ is
  spanned by $\uy_1,\dots,\uy_n$. 
The 
group $\Gamma_\minusone$ is generated by $n$ commuting involutions
  $\gamma_i$, $\gamma_i(x_j)=-1$ ($i\ne j$), $\gamma_i(x_i)=1$. 

Denote by $\mathit{pr}$ the projection 
$\id_{S_\minusone(V)} \tensor \varepsilon_W \tensor \varepsilon_{V^*} \colon  S_\minusone(V)\tensor \CC W \tensor
  S_\minusone(V^*)\to S_\minusone(V)$ onto $S_\minusone(V)$, where 
$\varepsilon_W\colon \CC W \to \CC$
is the algebra
  morphism such that $\varepsilon_W(w)=1$ for $w\in W$, and 
$\varepsilon_{V^*}\colon S_\minusone(V^*)\to \CC$ is the algebra morphism
such that $\varepsilon_{V^*}(V^*)=0$.
The \emph{braided Dunkl operators} attached to the group $W$ are
$$
 \rd\nabla_i\colon  S_\minusone(V) \to  S_\minusone(V), 
\qquad
  \rd\nabla_i(a) = \mathit{pr}(\uy_i a).
$$
Our last goal is to prove  formula (\ref{eq:braided Dunkl anticommutative})
for $\rd\nabla_i$, given in the Introduction. 

We put $\widetilde W=W\cdot \Gamma_\minusone\cdot\{\pm\id\}$ 
and embed $\rd{\mathcal H}_c(W)$, 
as a braided reduction, in the $\qq$\dash Cherednik algebra
  $\mathcal H_{c}(\widetilde W)$. Explicitly, we have $\uy_i=\gamma_i^{-1}
y_i = \gamma_i y_i$ for $i=1,\dots,n$, where $y_1,\dots,y_n$ span
$V^*$ in $\mathcal H_c(\widetilde W)$. 
The $\qq$\dash complex reflections in $\widetilde W$ are 
$$
    -\id\cdot s_{ij}^{(\varepsilon)}, \ \varepsilon\in\mathcal C; \quad 
    -\id\cdot t_i^{(\varepsilon')},\ \varepsilon'\in\pm\mathcal C', 
$$
where $-\id\cdot s_{ij}^{(\varepsilon)}= \gamma_i
\sigma_{ij}^{(\varepsilon)}$. 
(Formally, the parameters $c$ in $\rd{\mathcal H}_c(W)$ and in
$\mathcal H_c(\widetilde W)$ are not the same, but they are identified
in a rather obvious way.) 
We can find the generalised Dunkl
operators $\nabla_i$ for the algebra $\mathcal H(\widetilde W)$ and
then put $\rd\nabla_i = \gamma_i \nabla_i$. 

The generalised Dunkl operators $\nabla_i$ will be computed using the
procedure described in~\ref{subsect:revisited}. 
First, we embed 
$\mathcal H_c(\widetilde W)$ in
a modified Heisenberg quadratic double as in~\ref{subsect:embq}. 
We then have the following ``generalised root system''
of~\ref{subsect:revisited}: 
$$
   \alpha_{-\id\cdot s_{ij}^{(\varepsilon)}} = x_i - \varepsilon x_j, 
\quad
   \alpha^\vee_{-\id\cdot s_{ij}^{(\varepsilon)}} = c_1(y_i - \varepsilon^{-1} y_j), 
\quad
   \alpha_{-\id\cdot t_i^{(\varepsilon')}} = x_i,
\quad
   \alpha^\vee_{-\id\cdot t_i^{(\varepsilon')}} =
   c_{\varepsilon'}(1-\varepsilon')y_i.     
$$
It follows from \ref{subsect:revisited} that
$$
   \nabla_i = \partial_i+ \sum_w \langle \alpha_w, y_i \rangle \bar\partial_w,
$$
where 
$\partial_i$ is defined by $\partial_i(x_1^{a_1}\dots x_n ^{a_n}) = 
a_i x_1^{a_1}\dots x_i^{a_i-1}\dots x_n^{a_n}$ and 
$\bar\partial_w$ are uniquely defined operators on $S_\minusone(V)$
satisfying
$$
\bar\partial_w(x)=\langle x, \alpha^\vee_w\rangle, 
\qquad
\bar\partial_w(ab) = \bar\partial_w(a)w(b)+a\bar\partial_w(b),
$$
for $x\in V$, $a,b\in S_\minusone(V)$. 
If we know that $\bar\partial_w$ lowers the degree in $S_\minusone(V)$ by
$1$, both rules are equivalent to the equation
$$
     [\bar\partial_w, x]=\langle x,\alpha_w^\vee\rangle w 
$$ 
in $\End_\CC(S_\minusone(V))$. Put $\sigma_{ij}:=\sigma_{ij}^{(1)}$ and define
$$
D_{ij} =\frac{1}{{x_i^2-x_j^2}}\bigl((x_i+x_j)(1-\sigma_{ij})+(x_i-
x_j)(1-\sigma_{ji})\bigl) \ .
$$
We claim that $c_1 \gamma_i D_{ij} = \bar\partial_{-\id\cdot s_{ij}^{(1)}}
+\bar\partial_{-\id\cdot s_{ij}^{(-1)}}$. Indeed, 
it is not difficult to check, using the anticommutativity of the $x_i$
and the fact that $x_i^2 -x_j^2$ is central in $S_\minusone(V)$, that
$$
 [\gamma_i D_{ij}, x_i] = (-\id) \cdot s_{ij}^{(1)} + (-\id)\cdot
 s_{ij}^{(-1)},
\qquad
 [\gamma_i D_{ij}, x_j] = (-\id) \cdot s_{ij}^{(1)} - (-\id)\cdot
 s_{ij}^{(-1)},
\qquad
 [\gamma_i D_{ij}, x_k] = 0.
$$
Conjugating evertyhing with $t_j^{(\varepsilon)}$ shows 
  that $c_1 \gamma_i t_j^{(\varepsilon)} D_{ij} t_j^{(\varepsilon)-1} 
= \bar\partial_{-\id\cdot s_{ij}^{(\varepsilon)}}
+\bar\partial_{-\id\cdot s_{ij}^{(-\varepsilon)}}$. 

In the same way it is shown that
if $ D_i^{(\varepsilon')}=\dfrac{1}{x_i}(1-t_i^{(\varepsilon')})$,
then $c_{\varepsilon'} \gamma_i D_i =
\bar\partial_{t_i^{(\varepsilon')}}$ because 
$$
[\gamma_i D_i^{(\varepsilon')}, x_k] =
\delta_{ik}(1-\varepsilon')t_i^{(\varepsilon')}.
$$
(verified directly). 
We thus have the following expansion for $\nabla_i$:
$$
\nabla_i = \partial_i + \gamma_i 
c_1\sum_{j\ne i,\ \varepsilon \in 
\widetilde {\mathcal C}}
{t_j^{(\varepsilon)}}D_{ij} t_j^{(\varepsilon)-1}+ 
\gamma_i \sum_{\varepsilon' \in {\mathcal C}'\setminus\{1\}}
\frac{c_{\varepsilon'}}{1-\varepsilon'}D_i^{(\varepsilon')}\ ,
$$
where $\widetilde {\mathcal C}$ is a set of $|\mathcal C|/2$ elements
of $\mathcal C$ distinct modulo the subgroup $\{\pm 1\}$. 
Multiplying by $\gamma_i$ on the left and observing that
$\gamma_i\partial_i=\rd\partial_i$, we obtain 
formula~(\ref{eq:braided Dunkl anticommutative})
for $\rd\nabla_i$ as given in the Introduction.

\end{document}